\documentclass[11pt, twoside]{article} 
\usepackage{amsmath}    %obligatory package
\usepackage{amssymb}   	%obligatory package
\usepackage{amsthm}     %obligatory package
\usepackage{graphicx}    %obligatory package

\usepackage{textcomp}    %obligatory package
\usepackage[T1]{fontenc} %obligatory package
\usepackage{marvosym}  %obligatory package (for the \Letter symbol for indicating corresponding author)

\usepackage{authblk} 	%obligatory package
\usepackage[usenames]{xcolor} %obligatory package
\usepackage{lastpage} 	%obligatory package

\usepackage[latin1]{inputenc} % allows accented characters									
\usepackage{indentfirst}
\usepackage{amsfonts,amscd}	% symbols and special characters
\usepackage{empheq}
\usepackage{latexsym}		% special characters 
\usepackage{graphicx,graphics,epsfig}								
\usepackage{graphpap}										
\usepackage{float}
\usepackage{xcolor}
\usepackage{makeidx}      % remissive index
\usepackage{appendix}	  % appendix
\usepackage{multicol}										
\usepackage{color}

\usepackage[normalem]{ulem}
\usepackage{changepage}
\usepackage{hyperref} % optional arguments: [dvipdfm], [dvips]
\hypersetup{
	colorlinks=true,
	linkcolor=blue,         
   citecolor=blue,          
   urlcolor=blue          
}

\usepackage[dvips=false,pdftex=false,vtex=false]{geometry}
\usepackage{dsfont} % fonts for sets of nymbers, reals, complex, integers, etc.

\numberwithin{equation}{section}
\numberwithin{table}{section}
\numberwithin{figure}{section}

\newcommand{\orcid}[1]{\href{https://orcid.org/#1}{\includegraphics{logo_orcid.jpg}}} % orcid symbol with link to orchid page

\theoremstyle{definition}

\theoremstyle{definition}

\theoremstyle{definition}

\title{From Varadhan's Limit to Eigenmaps: A Guide to the Geometric Analysis behind Manifold Learning}

\author[1]{\textbf{Chen-Yun Lin %\orcid{0000-0002-6096-9716}
} \footnote{The first author is partially supported by a PSC-CUNY Grant, e-mail: chenyunlin.math@gmail.com}}
\author[2]{\textbf{Christina Sormani %\orcid{0000-0002-2295-2585}
}\footnote{The second author is partially supported by partially supported by NSF DMS 1612049 and a PSC-CUNY Grant, e-mail: sormanic@gmail.com}}
\affil[1]{Lehman College, 250 Bedford Park Boulevard West, Bronx, NY 10468, USA}
\affil[2]{CUNY Graduate Center, 365 Fifth Avenue, NY, NY 10468, USA}

\newcommand\shorttitle{From Varadhan's Limit to Eigenmaps}

\makeatletter
\def\volume#1{\def\@volume{#1}}

\def\shortauthor#1{\def\@shortauthor{#1}}

\def\shorttitle#1{\def\@shorttitle{#1}}

\makeatother

\date{2023} % year of publication 
\volume{\#} % number of the volume

\makeatletter
\let\thetitle\@title
\let\thedate\@date
\let\thevol\@volume
\let\theauthor\@author
\let\theshortauthor\@shortauthor
\let\theshorttitle\@shorttitle 
\makeatother

\makeatletter
\renewcommand{\maketitle}{\bgroup\setlength{\parindent}{0pt}

\parindent=1em
\renewcommand{\thefootnote}

\phantom{ }

\vspace{-1.5cm} \noindent\includegraphics[scale=.35]{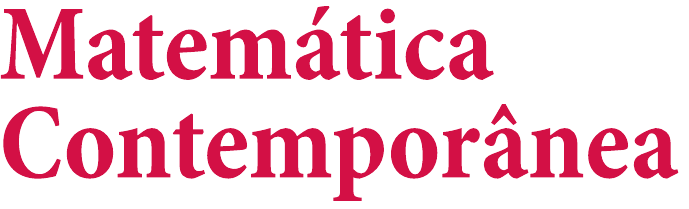}
 \hfill
 \includegraphics[scale=0.07]{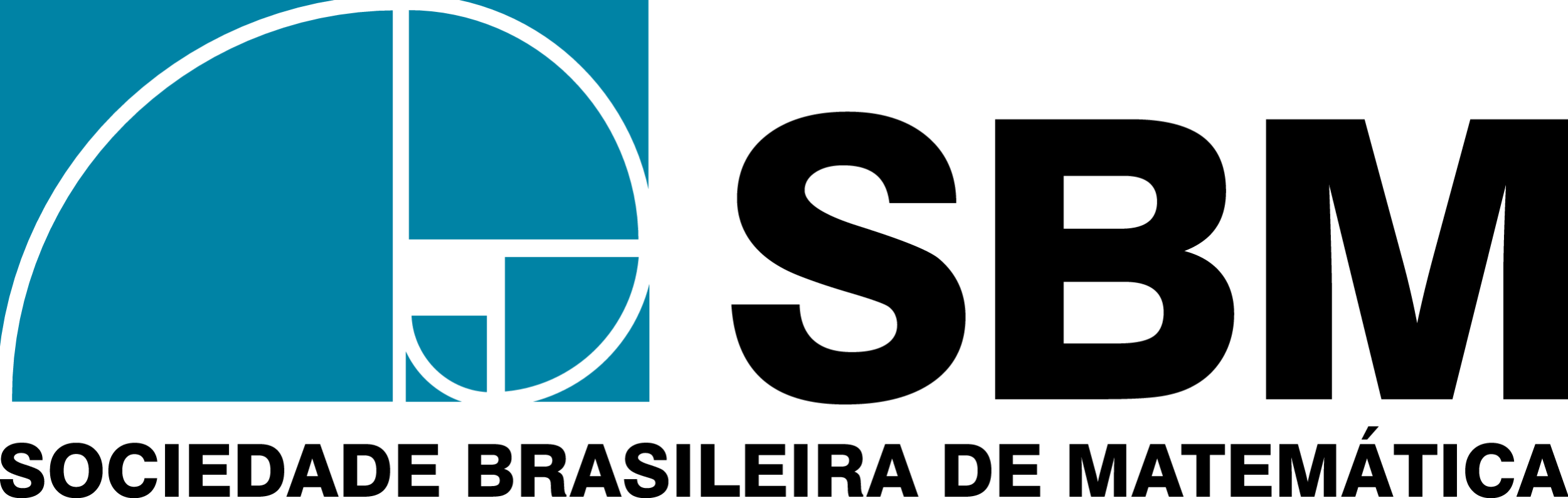}

\vspace*{-.2cm}

\noindent{\scriptsize Vol. \thevol,\  \pageref{FirstPage}--\pageref{LastPage}} \hfill \copyright \thedate \qquad \qquad

\vspace*{-.2cm}

\noindent{\scriptsize\url{http://doi.org/will\_be\_added\_later}}
\vspace{1truecm}
\begin{center}{\vbox{\titlefont\@title}}\end{center}
\vspace{0.5truecm}
\begin{center}{\@author} \end{center}

\egroup
}

\renewcommand{\thefootnote}{\fnsymbol{footnote}}
\renewcommand{\@fnsymbol}[1]{%
    \ifcase#1 \or {\,\Letter\!} \or\textasteriskcentered\or \textasteriskcentered\textasteriskcentered 
    \else\@ctrerr\fi}
\makeatother

\newcommand*{\titlefont}{\fontsize{18}{21.6}\selectfont\textbf}

\makeatletter
\renewcommand\@author{\ifx\AB@affillist\AB@empty\AB@author\else
      \ifnum\value{affil}>\value{Maxaffil}\def\rlap##1{##1}%
    \AB@authlist\\[\affilsep]\vbox{\AB@affillist}
    \else  \AB@authors\fi\fi}
\makeatother

\usepackage{amsmath, amsthm, amssymb}
\usepackage{amsfonts}
\usepackage{graphicx}
\usepackage{xcolor}
\usepackage{mathbbol}
\usepackage{hyperref}
\usepackage{txfonts}
\usepackage{soul}
\usepackage{mathtools}  
\usepackage{tabulary}
\usepackage{booktabs}

\newtheorem{thm}{Theorem}[section]
\newcommand{\bt}{\begin{thm}}
\newcommand{\et}{\end{thm}}

\newtheorem{conj}[thm]{Conjecture}

\newtheorem{cor}[thm]{Corollary}   %remember switch all {coro} to {cor}
\newcommand{\bc}{\begin{cor}}
\newcommand{\ec}{\end{cor}}

\newtheorem{lem}[thm]{Lemma}   %remember to switch all {lem} to {lem}
\newcommand{\bl}{\begin{lem}}
\newcommand{\el}{\end{lem}}

\newtheorem{prop}[thm]{Proposition}
\newcommand{\bp}{\begin{prop}}
\newcommand{\ep}{\end{prop}}

\newtheorem{defn}[thm]{Definition}
\newcommand{\bd}{\begin{defn}}    % This produces an error????    
\newcommand{\ed}{\end{defn}}

\newtheorem{rmrk}[thm]{Remark}   %remember to switch all {rmrk} to {rmrk}

\newcommand{\br}{\begin{rmrk}}
\newcommand{\er}{\end{rmrk}}

\newcommand{\mGHto}{\stackrel { \textrm{mGH}}{\longrightarrow} }
\newcommand{\GHto}{\stackrel { \textrm{GH}}{\longrightarrow} }

\newcommand{\be}{\begin{equation}}

\newcommand{\ee}{\end{equation}}

\newcommand{\injrad}{\textrm{injrad}}

\newcommand{\diam}{\operatorname{Diam}}

\newcommand{\Ricci}{\rm{Ricci}}
\newcommand{\Scal}{{\rm Scal}} % new Jan 2016
\newcommand{\disjointunion}{\sqcup}

\newcommand{\Lip}{\operatorname{Lip}}

\newcommand{\sect}{\operatorname{sect}}
\newcommand{\vol}{\operatorname{Vol}}

\begin{document}

\maketitle
\thispagestyle{plain}
\renewcommand{\thefootnote}{\arabic{footnote}}   % don't delete!!!
\setcounter{footnote}{0}     % don't delete!!!
\setcounter{page}{1} % The editors will insert the correct initial pagenumber
\label{FirstPage}	

\pagestyle{myheadings} \markboth{\hfil  $\hspace{1.5cm}$  
C-Y Lin and C Sormani
\hfil $\hspace{3cm}$ } {\hfil$\hspace{1.5cm}$
{From Varadhan's Limit to Eigenmaps:}
\hfil}

%%%%%%%%%%%
% Minipage with: dedication, Abstract, keywords, and MSC
%%%%%%%%%%%
\begin{center}
\noindent
\begin{minipage}{0.85\textwidth}\parindent=15.5pt

\smallskip
\begin{center}
\large{\textsl{Dedicated to Professor Varadhan,\\ Fellow of The World Academy of Sciences  }}
\end{center}

\smallskip

{\small{
\noindent {\bf Abstract.} We present an overview of the history of the heat kernel and eigenfunctions on Riemannian manifolds
and how the theory has lead to modern methods of analyzing high dimensional data.  }
\smallskip

% Please enter at most 6 keywords here with lowercase letters separated by commas.
\noindent {\bf{Keywords:}} heat kernel embeddings, eigenmaps, manifold learning.
\smallskip

% Please enter at most 5 Mathematics Subject Classification codes here. Please use 2010 classification codes, which can be found on the following link: http://www.ams.org/msc//msc2010.html.
\noindent{\bf{2020 Mathematics Subject Classification:}} 58J35, 35K05, 53C23.
}

\end{minipage}
\end{center}

\tableofcontents

\section{\bf Introduction}
Geometry is one of the most beautiful and ancient forms of mathematics.   Inspired by Professor Romain Murenzi's speech at the {\em 2022 TWAS TYAN Virtual Workshop on Differential Geometry}, we will present the Geometric Analysis that leads a very young and active branch of Geometry called Manifold Learning.  Manifold Learning is a form of machine learning that has applications in biology, computer vision, natural language processing, and data mining.   It is built upon half a century of Geometric Analysis.   Here we review the history from Varadhan's theorem relating the heat kernel and the distance on a Riemannian manifold, to the
study of eigenfunctions and eigenvalues in relation to the heat kernel, to convergence of Riemannian
manifolds, to embeddings of Riemannian manifolds using the heat kernel, to the modern notions of eigenmaps and dimension reduction through manifold learning.   Throughout we provide detailed examples demonstrating the ideas.  We assume only knowledge of Professor Manfredo P. do Carmo's {\bf{Riemannian Geometry}} and we provide a guide to both the original papers and various useful textbooks.

\section{\bf Heat Kernels and Distances on Riemannian Manifolds}\label{sect-heat}

We will begin with a 1967 theorem of TWAS Fellow and Abel Laureate, Srinivasan Varadhan, who earned his doctorate at the Indian Statistical Institute in 1963:

\begin{thm} \cite{Varadhan:1967} The distance between two points, $x$ and $y$ in a compact Riemannian manifold, $M$, satisfies
\be \label{V-lim}
d(x,y)^2=\lim_{t\to 0} -4t \log h_t(x,y)
\ee
where $h_t(x,y)$ is the heat kernel.  This convergence is uniform on bounded regions, $U$,
in the sense that: $\forall \varepsilon>0 \,\,\exists T_{U,\varepsilon}>0$ such that $\forall t\in (0, T_{U,\varepsilon})$ we have
\be
| d(x,y)^2+4t\log(h_t(x,y))|<\varepsilon 
\ee
which is equivalent to saying that 
\be
e^{-d^2(x,y)/(4t)}\cdot e^{-\varepsilon/(4t)} < h_t(x,y)< e^{-d^2(x,y)/(4t)}\cdot e^{\varepsilon/(4t)}.
\ee
\end{thm}

Recall that the heat kernel 
\be
h_t: M\times M \to (0,\infty) \textrm{ defined for } t>0
\ee
is defined to be the unique function such that
\be
u(x,t)=\int_M f(y) h_t(x,y) \, dvol_y,
\ee
for any solution, 
$u:M\times (0,\infty)\to [0,\infty)$ of the heat equation, 
\be
\Delta_x u(x,t)= \partial_t u(x,t) \textrm{ with initial data } \lim_{t\to 0^+} u(x,t)=f(x).
\ee 
As a consequence, for any $x\in M$,
\be\label{normalize}
\int_{M} h_t(x,y) \, dvol_y = 1.
\ee
Furthermore, it is easy to verify the heat semigroup property holds:
\be\label{semigroup}
\int_M h_t(x,y)h_s(y,z) \, dvol_y=h_{t+s}(x,z).
\ee

For a good resource with the details as to why the heat kernel exists on a compact Riemannian manifold and proofs of its various properties, we refer the reader to Chavel's textbook {\bf Eigenvalues in Riemannian Geometry} \cite{Chavel-text}.  For more theorems about the heat kernel see Schoen-Yau's text {\bf Lectures on Differential Geometry} \cite{Schoen-Yau-text}.   

Varadhan proved his theorem by estimating the heat kernel on normal neighborhoods 
using the Euclidean heat kernel
and then built his estimate on $h_t(x,y)$ for $x$ far from $y$ by applying the heat semigroup property as well
as the properties of the exponential map based at $x$ near the minimizing geodesic from $x$ to $y$.   This proof appeared in \cite{Varadhan:1967} and is quite readable for geometers .    He wrote a probabilistic proof in \cite{Varadhan:1967diffusion}. 

\subsection{The Heat Kernel on a Circle}\label{sect-circle-heat}

On a circle, ${\mathbb{S}}^1$ parametrized by $\theta\in [0,2\pi)$ with $g_{{\mathbb{S}}^1}=d\theta^2$,
then the Laplacian is $\Delta=\partial_\theta^2$ and
the heat kernel is the following sum which converges due to the exponential decay:
\be\label{circle-kernel}
h_t(\theta_1,\theta_2)=\sum_{m\in {\mathbb{Z}} }\tfrac{1}{\sqrt{4\pi t} } \, e^{-(\theta_1-\theta_2+2\pi m)^2/(4t)}.
\ee
Notice that as $t\to 0$, the term in the sum which dominates because it decays the least quickly, is the
term with 
\be
e^{-(\theta_1-\theta_2+2\pi m)^2/(4t)}= e^{ -d_{{\mathbb S}^1}^2(\theta_1,\theta_2)^2/(4t)}
\ee
where $d_{{\mathbb S}^1}(\theta_1,\theta_2)$ is the Riemannian distance between the points in ${\mathbb{S}}^1$:
\be
d_{{\mathbb S}^1}(\theta_1,\theta_2)=\min\{ |\theta_1-\theta_2+2\pi m|\,:\, m \in {\mathbb Z} \}
\ee
which is the length of the smaller arc between $\theta_1$ and $\theta_2$.  

\subsection{The Heat Kernel on a Rescaled Circle}\label{sect-recircle-heat}

On a circle, ${{\mathbb{S}}^1_R}$, we rescale the circle by a factor $R>0$ so that $x=R\theta\in [0,2\pi R)$ with 
$
g_R=dx^2=R^2d\theta^2.
$
Then the Laplacian rescales $\Delta_{R}=\partial_x^2=R^{-2} \partial_\theta^2$ and the heat kernel
is 
\be\label{circle-kernel_R1}
h^R_t(\theta_1,\theta_2)=\sum_{m\in {\mathbb{Z}} }\tfrac{1}{\sqrt{4\pi t}\,\, } \, e^{-(\theta_1-\theta_2+2\pi m)^2 R^2/(4t)}.
\ee

Note that the heat kernel of ${\mathbb S}_R^1$ can be related to the heat kernel of  ${\mathbb S}^1$
with dilated time:
\be\label{rescaled-heat}
h^R_t(\theta_1,\theta_2)=  R^{-1} \,h_{t/R^2}(\theta_1,\theta_2)
\ee
because
\be
\sum_{m\in {\mathbb{Z}} }\tfrac{1}{\sqrt{4\pi t} \,\,} \, e^{-(\theta_1-\theta_2+2\pi m)^2 R^2/(4t)}=
R^{-1} \sum_{m\in {\mathbb{Z}} }\tfrac{1}{\sqrt{4\pi (t/R^2)} \,\,} \, e^{-(\theta_1-\theta_2+2\pi m)^2/(4(t/R^2))}
\ee
%Varadhan time too hard to compute

\subsection{The Heat Kernel on Tori} \label{sect-tori-heat}
Next let us consider the flat torus formed by taking the isometric product of rescaled circles,
\be
{\mathbb{S}}_A^1\times {\mathbb{S}}_B^1=({\mathbb{S}}^1\times {\mathbb{S}}^1, \, A^2\,d\theta^2 + B^2\,d\varphi^2).
\ee
Its heat kernel is
\be \label{heat-AB}
h^{A,B}_t((\theta_1,\varphi_1),(\theta_2,\varphi_2))=
\sum_{n,m\in {\mathbb{Z}} }\tfrac{1}{(4\pi t) } 
\, e^{-\left( (\theta_1-\theta_2+2\pi n)^2 A^2+(\varphi_1-\varphi_2+2\pi m)^2 B^2\right) /(4t)}.
\ee
Notice that
\be
4t \log h^{A,B}_t((\theta_1,\varphi_1),(\theta_2,\varphi_2))=
4t \log h^{A}_t(\theta_1,\theta_2)+4t \log h^{B}_t(\varphi_1,\varphi_2).
\ee
Since the distances squared also add,  we have
\be
d^2_{{\mathbb{S}}_A^1\times {\mathbb{S}}_B^1}((\theta_1,\varphi_1),(\theta_2,\varphi_2))=d^2_{{\mathbb{S}}_A^1}(\theta_1,\theta_2)+d^2_{{\mathbb{S}}_B^1}(\varphi_1,\varphi_2)
\ee
which matches what Varadhan's limit predicts in this setting.

\subsection{The Heat Kernel on Spaces of Constant Curvature}
In general on a homogeneous manifold, $M^n$, the heat kernel will be a function only of distance and time:
\be
h^M_t(x,y)= h^M_t(s) \textrm{ where } s=d(x,y).
\ee

The heat kernel on a sphere was first computed by Minakshisundaram \cite{Minak-sphere}.  
%\textcolor{red}{I CANNOT FIND THIS ARTICLE }
%{\color{blue} Neither do I. }
The heat
kernel on hyperbolic space was first computed  in dimension 2 to be
\be
h^{{\mathbb H}^{2}}_t(s)=\frac{\sqrt{2}}{(4\pi t)^{3/2}} e^{-t/4}\int_s^\infty  
\frac{u e^{-u^2/(4t)} }{ (\cosh(u)-\cosh(s))^{1/2} } \, du
\ee
by McKean \cite{McKean:1970}.   In three dimensions it was shown to be
\be
h^{{\mathbb H}^{3}}_t(s)=\frac{(-1)}{2\pi}\frac{1}{\sqrt{4\pi t\,}}
\left(\frac{1}{\sinh(s)}\frac{\partial}{\partial s}\right) e^{-t}e^{-s^2/4t}
\ee
by Davies-Mandouvalos in \cite{Davies-Mandouvalos:1988}.  They applied an
unpublished recurrence relation proved by Millson:
\be
h^{{\mathbb H}^{m+2}}_t(s)= -\frac{e^{-mt}}{2\pi \sinh(s)} \frac{\partial}{\partial s}   h^{{\mathbb H}^{m}}_t(s)
\ee
to provide the formula for the heat kernel in all dimensions.
 See also the
alternative proof by Grigor'yan-Noguchi in \cite{Grigoryan-Noguchi:1998}.

When $M_\kappa^n$ is the simply connected Riemannian manifold of constant curvature $\kappa$, we write
\be \label{h-kappa}
h^{\kappa,n}_t(s) = h^{M_\kappa^n}_t(s).
\ee
Since $M_\kappa^n$ is just a rescaled sphere ${\mathbb{S}}^n_{\kappa^{-1/2}}$ when $\kappa>0$
and is just rescaled hyperbolic space ${\mathbb{H}}^n_{\kappa^{-1/2\,}}$ when $\kappa<0$, we know
%{\color{blue} \href{https://mathoverflow.net/questions/409193/heat-kernel-on-hyperbolic-space-of-variable-curvature}{rescaling}   CHECK WITH BETTER SOURCE
\be
h^{\kappa,n}_t(s) = \kappa^{n/2} h^{{\mathbb S}^{n}}_{\kappa t}(\kappa^{1/2} s) \textrm{ for } \kappa>0
\ee
and
\be
h^{\kappa,n}_t(s) = |\kappa|^{n/2} h^{{\mathbb H}^{n}}_{|\kappa| t}(|\kappa|^{1/2} s) \textrm{ for } \kappa<0.
\ee

\subsection{Estimating the Heat Kernel on Families of Manifolds}\label{sect-est-heat}

In 1976 Debiard-Gaveau-Mazet \cite{Debiard-et-al:1976} estimated the heat kernel uniformly from above and below.

\begin{thm}
 For simply connected manifolds, $M$, with negative sectional curvature between bounds 
\be\label{kappaMkappa'}
\kappa\le \sect_M \le \kappa'<0,
\ee
the heat kernel satisfies
\be \label{DGM}
h^{\kappa',n}_t(s) \ge  h^M_t(x,y) \ge h^{\kappa,n}_t(s). 
\ee
where $h^{\kappa,n}(t)(s)$ is defined as in (\ref{h-kappa}) and where $s=d_M(x,y)$. 
\end{thm}

Applying this theorem, we thus have
\be 
s^2+4t \log h^{\kappa,n}_t(s) \ge  d_M(x,y)^2+4t \log h^M_t(x,y) \ge s^2+4t \log h^{\kappa',n}_t(s).
\ee
So $\left|d_M(x,y)^2+4t \log h^M_t(x,y)\right| $ is bounded above by
\be
 \max\left\{|s^2+4t \log h^{\kappa,n}_t(s)| , |s^2+4t \log h^{\kappa',n}_t(s)| \right\}.
\ee
Thus the uniform convergence in Varadhan's Theorem is uniform even allowing the Riemannian manifold to vary within this class of Riemannian manifolds $M$ such that (\ref{kappaMkappa'}) holds.

In 1981, Cheeger and Yau \cite{Cheeger-Yau:1981} proved the following estimate that only requires uniform lower bounds on the Ricci curvature.

\begin{thm}
If $M^n$ has $\Ricci \ge (n-1)H$ for any $H\in (-\infty, \infty)$,
then 
we still have a lower bound
\be\label{Cheeger-Yau-heat}
h^M_t(x,y) \ge h^{H,n}_t(d_M(x,y)).
\ee
\end{thm}

So if $M$ is also simply connected with $\sect_M\le \kappa'$ we still have (\ref{DGM}) and its consequences.   

In 1982, Cheeger-Gromov-Taylor found an upper bound on the heat kernel for points $x,y \in B(p,D)$
depending on $D$, on the injectivity radius and the volumes of balls in the manifold (See Theorem 1.4 in \cite{Cheeger-Gromov-Taylor}).  In 1984 Chavel found an upper bound on the heat kernel which depends on isoperimetric constants of balls (See Theorems 8 and 9 in Chapter VIII of \cite{Chavel-text}) thus there is a uniform bound on the Varadhan time for manifolds with same isoperimetric profile functions and the same lower bound on Ricci curvature.   These bounds may also be used to provide uniform control on the heat kernel.

Cheng-Li-Yau studied complete manifolds and discussed some of the
difficulties proving upper bounds before finding an upper bound on the growth of the heat kernel \cite{Cheng-Li-Yau:1981}.   
Cheng and Li found upper bounds on the heat kernel depending on the Sobolev constant of the manifold \cite{Cheng-Li:1981}.    Additional bounds on the heat kernel and also a Harnack Inequality for Riemannian manifolds is proven by Li and Yau in 1986 \cite{Li-Yau:1986}.   
This was then applied by Li in \cite{Li:1986} to estimate the behavior of the heat kernel as $t\to \infty$.   For more details see the text of Schoen and Yau \cite{Schoen-Yau-text}.   

It should be noted that Varadhan's original proof of his limit in \cite{Varadhan:1967diffusion} was completed using probabilistic methods and there has been much work done using these methods that we do not present here.  These
results would require an expertise in probability and would better be presented elsewhere.

\section{\bf Eigenfunctions and Eigenvalues}\label{sect-evalue}

One of the most beautiful properties of the heat kernel on a compact Riemannian manifold is its Sturm-Liouville decomposition:
\be \label{Sturm-Liouville}
h_t(x,y)=\sum_{j=0}^\infty e^{-\lambda_jt} \phi_j(x)\phi_j(y)
\ee 
where $0=\lambda_0<\lambda_1\le \lambda_2\le \cdots$ are eigenvalues and
$\phi_j$ are an $L^2$ orthonormal collection of eigenfunctions of the Laplacian:
\be
\Delta \phi_j(x)=-\lambda_j \phi_j(x) 
\ee
For a quick hint at why this expansion gives the heat kernel, note that each term
satisfies the heat equation:
\be \label{each-term}
\partial_t e^{-\lambda_jt} \phi_j(x)= -\lambda_j e^{-\lambda_jt} \phi_j(x)=\Delta e^{-\lambda_jt} \phi_j(x).
\ee
It is also symmetric $h_t(x,y)=h_t(y,x)$.   We will explain in more detail below.

Note that Sturm-Liouville theory is named for Jacques Charles Francois Sturm (1803-1855) and Joseph Liouville (1809-1882).  From the age of 16, Sturm supported his mother and siblings as a mathematics tutor: he moved to Paris to serve as a tutor to a student of his that continued on to study there.   He became a secondary school teacher after the Revolution of 1830 and by 1834 had published his work developing Sturm-Liouville Theory with Joseph Liouville.  Joseph Liouville was the son of an army officer in France and had studied engineering at Ecole Polytechnique, but was working as a mathematics assistant at various institutions at the time.  After the publication of this important work and other theorems, both became professors of mathematics at Ecole Polytechnique in 1838. 

See also work of Miranda \cite{Miranda:1955} applying Sturm-Liouville Theory to Riemannian manifolds.   This is 
well explained in Chavel's textbook \cite{Chavel-text}  and Schoen-Yau's text  \cite{Schoen-Yau-text}.   For a survey of earlier results, the text of Berger, P. Gauduchon and E. Mazet, \cite{Berger-et-al:1971} is often recommended.  This last book includes many explicit computations of eigenvalues and eigenfunctions for a variety of Riemannian manifolds.

\subsection{Finding Eigenfunctions using the Rayleigh Quotient}\label{sect-find-evalues}

Let us begin this section by reviewing how to find the eigenfunctions.
Recall that $L^2(M)$ is the Hilbert space of square Lebesgue integrable functions:
\be
L^2(M)=\left\{ \phi:M \to {\mathbb R}\,\, |\,\, \int_M \phi^2(x) \,dvol_x  < \infty \right\}
\ee
with the $L^2$-inner product and $L^2$-norm,
\be 
<\phi\,, \,\tilde{\phi}>_{L^2(M)}=\int_M  \phi(x) \tilde{\phi}(x) \, dvol_x 
\quad \textrm{ and } \quad ||\phi ||^2_{L^2}=\int_M  \phi^2(x)  \, dvol_x.
\ee
To find the eigenfunctions, we begin with $\lambda_0=0$ and its
constant eigenfunction, $\phi_0(x)=\phi_0$, which is normalized by 
\be \label{phi-0}
1=||\phi_0||^2_{L^2(M)}=\int_M \phi_0^2 \, dvol = \phi_0^2\vol(M) \quad \implies \quad \phi_0(x)=\vol(M)^{-\tfrac{1}{2}\,\,}.
\ee
We find the rest of the eigenvalues iteratively using Rayleigh's Quotient:
\be\label{Rayleigh}
\lambda_j=\inf  \int_M |\nabla \phi(y)|^2 \, dvol \, 
\ee
where the infimum is taken over all $\phi \in L^2(M)$ such that
\be \label{3.8}
||\phi ||_{L^2}=1 \textrm{ and }<\phi\,, \,\phi_k>_{L^2}=0 \textrm{ for } k=0,...,(j-1).
\ee

Let us review why the above Rayleigh's Quotient method works.
Observe first that the conditions on $\phi$ in (\ref{3.8}) are designed to guarantee an $L^2$ orthonormal collection.  Next recall Green's Theorem
\be
\int_M g(\nabla \phi(y), \nabla \tilde \phi) \, dvol = -\int_M\phi(y) \Delta \tilde \phi \, dvol.
\ee
Observe that any $\phi_j$ such that $\Delta \phi_j=-\lambda_j\phi_j$ with $\lambda_j>0$
has $<\phi_j\,, \,\phi_k>_{L^2}=0$ for $\lambda_k<\lambda_j$ because
\begin{eqnarray}
<\phi_j\,, \,\phi_k>_{L^2}
&=&\int_M \phi_j(x)\phi_k(x)\, dvol_x= \tfrac{1}{\lambda_j}\int_M \Delta\phi_j(x)\phi_k(x)\, dvol_x\quad \\
&=&\tfrac{1}{\lambda_j}\int_M \phi_j(x)\Delta\phi_k(x)\, dvol_x\\
&=&\tfrac{\lambda_k}{\lambda_j}\int_M \phi_j(x)\phi_k(x)\, dvol_x=  \tfrac{\lambda_k}{\lambda_j} <\phi_j\,, \,\phi_k>_{L^2}.
\end{eqnarray}
If we rescale $\phi_j$ so that $||\phi_j||_{L^2}=1$, then it achieves the infimum in (\ref{Rayleigh}):
\begin{eqnarray}
 \int_M |\nabla \phi_j(y)|^2 \, dvol_y&=&-\int_M \phi_j(y) \Delta \phi_j(y) \, dvol_y\\
 &=&\lambda_j \int_M  \phi_j(y) \phi_j(y) \, dvol_y=\lambda_j.
\end{eqnarray}
It is not difficult to prove that the collection of $\phi: M \to {\mathbb R}$ which achieve the infimum in (\ref{Rayleigh}) is
exactly the finite dimensional eigenspace for $\lambda_j$.  
We then choose an $L^2$ orthonormal basis of eigenfunctions for each eigenspace.   

Finally we see that a  complete collection of eigenfunctions found using this method forms a basis for $L^2(M)$
by proving that they span $L^2(M)$.  That is, every $f\in L^2(M)$
has a Fourier Expansion:
\be \label{trunc-f}
f(x) = \lim_{N\to \infty} f_N(x) \textrm{ where } f_N(x)=\sum_{j=0}^N <\phi_j \,, f>_{L^2(M)} \phi_j(x)
\ee
which converges in $L^2(M)$.   See for example the Chavel textbook \cite{Chavel-text} for more details.

\subsection{The Sturm-Liouville Decomposition of the Heat Kernel}

If we wish to solve the heat equation with initial data from (\ref{trunc-f}),
\be
\partial_t u_N(x,t)=\Delta u_N(x,t) \textrm{ and } \lim_{t\to 0^+} u_N(x,t)=f_N(x),
\ee
we can solve one term at a time as in (\ref{each-term}) to see that
\begin{eqnarray}
u_N(x,t)&=&\sum_{j=0}^N <\phi_j \,, f>_{L^2} e^{-\lambda_jt} \phi_j(x)\\
&=&\sum_{j=0}^N  \int_M \phi_j(y) f(y) \, dvol_y e^{-\lambda_jt} \phi_j(x)\\
&=& \int_M f(y) \sum_{j=0}^N e^{-\lambda_jt} \phi_j(x)\phi_j(y)  \, dvol_y\\
&=& \int_M f(y) h^N_t(x,y) \, dvol_y
\end{eqnarray}
where 
\be
h^N_t(x,y)=\sum_{j=0}^N e^{-\lambda_jt} \phi_j(x)\phi_j(y).
\ee
To complete the proof of the Sturm-Liouville Decomposition, we must
cautiously take $N\to \infty$, and show $u_N(x,t)$ converges
to a solution $u(x,t)$ of the heat equation
\be
\partial_t u(x,t)=\Delta u(x,t) \textrm{ and } u(x,0)=f(x).
\ee
and
\be
u(x,t)=\int_M f(y) (\lim_{N\to \infty} h^N_t(x,y)) \, dvol_y
\ee
so that 
\be
\lim_{N\to \infty} h^N_t(x,y)=h_t(x,y)
\ee
 with convergence in $L^2(M\times M)$.
See Chavel's text \cite{Chavel-text} for details.

The trace of the heat kernel is then defined naturally as follows:
\be\label{heat-trace}
\int_{x\in M} h_t(x,x) \, dvol_x = 
\sum_{j=0}^\infty e^{-\lambda_jt}.
\ee

%\textcolor{red}{MUST Mention the Weyl expansion here.   We have not mentioned it anywhere.}

%\textcolor{red}{Maybe show decay using Coifman Lafon type estimate.  Or check Chavel or Schoen-Yau}

\subsection{Truncating the Sturm-Liouville Expansion of the Heat Kernel}\label{sect-trunc-heat}

Let us consider how large we must take $N$ in order to estimate the heat kernel with the truncated sum
\be \label{truncated}
h^N_t(x,y)=\sum_{j=0}^N e^{-\lambda_jt} \phi_j(x)\phi_j(y).
\ee
We will define the {\em uniform Sturm-Liouville number}, $N_{t,\epsilon}$ to be the
smallest natural number such that
\be \label{S-L-1}
|h_t(x,y) -h^N_t(x,y)| < \epsilon \qquad \forall N\ge N_{t,\epsilon}.
\ee 
We also define the {\em $L^2$ Sturm-Liouville number}, $N_{t,\epsilon, L^2}$ to be the
smallest natural number such that $\forall N\ge N_{t,\epsilon, L^2}$ we have
\be \label{S-L-2}
||h_t(x,y) -h^N_t(x,y)||_{L^2(M\times M)}=\int_{M\times M} (h_t(x,y) -h^N_t(x,y))^2\,dvol_xdvol_y < \epsilon.
\ee 
The {\em $L^2$ Sturm-Liouville number} is easy to estimate because $||\phi_j(x)||_{L^2(M)}=1$ and
\begin{eqnarray*}
||\phi_j(x)\phi_j(y)||^2_{L^2(M\times M)}&=&\int_M\int_M (\phi_j(x)\phi_j(y))^2\,dvol_xdvol_y\\
&=&
\int_M\phi_j^2(x)\,dvol_x \,\, \cdot\,\,\int_M \phi^2_j(y) \,dvol_y\\
&=&
 ||\phi_j||^2_{L^2(M)} \, \cdot \, ||\phi_k||^2_{L^2(M)}=1.
\end{eqnarray*}
So, by the $L^2$ triangle inequality,
\be \label{tail-estimate}
||h_t(x,y) -h^N_t(x,y)||_{L^2(M\times M)} \le \sum_{j=N+1}^\infty e^{-\lambda_jt} ||\phi_j(x)\phi_j(y)||_{L^2(M\times M)}
= \sum_{j=N+1}^\infty e^{-\lambda_jt}. 
\ee
This last sum is the {\em ``tail of the trace of the heat kernel"} so it converges to $0$ as $N\to \infty$
and can even be estimated knowing only information about the eigenvalues of the manifold.
 
%\textcolor{red}{ADD ABOUT ITS DECAY AND N USED TO ESTIMATE IT.  Maybe a subsection on the tail which has exactly what Coifman-Lafon wrote later.  FORGET COIFMAN LAFON.}.

\subsection{Eigenvalues and Eigenfunctions on a Standard Circle}\label{sect-circle-evalue}
On a circle, ${\mathbb{S}}^1$ parametrized by $\theta\in [0,2\pi)$ with Laplacian is 
\be
\Delta=\partial_\theta^2.
\ee
The zeroeth eigenvalue is $\lambda_0=0$ with $\phi_0=(2\pi)^{-\tfrac{1}{2}\,\,}$ by (\ref{phi-0}).
The first eigenvalue is $\lambda_1=1$ and the eigenspace for this eigenvalue consists of periodic eigenfunctions of the form $A\sin(\theta-\theta_0)$ where $\theta_0\in [0,2\pi)$.   This eigenspace is 2 dimensional 
 $\lambda_2=\lambda_1=1$ and we say this eigenvalue has multiplicity 2.   We can choose two orthonormal eigenfunctions
 \be
 \phi_1(\theta)=A_1\sin(\theta) \quad \textrm{ and } \quad \phi_2(\theta)=A_1\sin(\theta-\pi/2)=A_1\cos(\theta)
 \ee
 where 
 \be
 1=\int_0^{2\pi} A_1^2 \sin^2(\theta)\, d\theta= (\tfrac{1}{2}\,) \int_0^{2\pi} A_1^2 \, d\theta=(\tfrac{1}{2}\,)2\pi A^2_1,\quad \textrm{ so } A_1=\pi^{-\tfrac{1}{2}\,\,}. 
 \ee
Continuing forward we find that $\lambda_{2j}=\lambda_{2j-1}=j^2$ with
\be
\phi_{2j-1}(\theta)=A_j\sin(j\theta) \quad \textrm{ and } \quad \phi_{2j}(\theta)=A_j\cos(j\theta)
 \ee
 where 
 \be
 1=\int_0^{2\pi} A_j^2 \sin^2(j\theta)\, d\theta= (\tfrac{1}{2}\,) \int_0^{2\pi} A_j^2 \, d\theta=\pi A_j^2 \textrm{ so }
 A_j=\pi^{-\tfrac{1}{2}\,\,} 
 \ee
So following (\ref{Sturm-Liouville}), the truncated Sturm-Liouville expansion for the heat kernel ${\mathbb S}^1$ is
 \be \label{ST-S1}
h^{2N+1}_t(\theta_1,\theta_2)=
 \tfrac{1}{2\pi} + \tfrac{1}{\pi}  \sum_{j=1}^N e^{-t j^2} \left( \sin(j\theta_1)\sin(j\theta_2) + \cos(j\theta_1)\cos(j\theta_2 ) \right).
\ee

The trace of the heat kernel on ${\mathbb S}^1$ is
\be
\int_{\theta=0}^{2\pi} h_t(\theta,\theta) \, d\theta =  1 +  2 \sum_{j=1}^\infty  e^{-j^2t} 
\ee
Since the eigenfunctions of ${\mathbb S}^1$ are uniformly bounded, $|\phi_j(x)| \le \pi^{-{1}/{2}\,\,}$, the Sturm-Liouville expansion
converges as quickly as the trace of heat kernel and we can use the same estimates to compute
the uniform and the $L^2$ Sturm-Liouville numbers.
%\be
%\textcolor{red}{ADD NUMBERS HERE? similar to Coifman Lafon}
%\ee

If we try to apply Varadhan's formula to the truncated Sturm-Liouville expansion for the heat kernel, we see that
as $t\to 0$, the truncated heat kernels, $h_t^{2N+1}(\theta_1,\theta_2) $,
converge to
\be
(2\pi )^{-1}  + (\pi )^{-1} \sum_{j=1}^N  \left( \sin(j\theta_1)\sin(j\theta_2) + \cos(j\theta_1)\cos(j\theta_2 ) \right).
\ee
So we end up with nothing:
\be
\lim_{t\to 0} -4t \log h_t^{2N+1}(\theta_1,\theta_2))=0
\ee
In order to estimate Varadhan's formula using the Sturm-Liouville expansion with an error less than $\varepsilon>0$, 
we would have to first fix $t_\varepsilon$ small enough that
\be
|(-4t \log h_t(\theta_1,\theta_2))-d^2_{{\mathbb S}^1}(\theta_1,\theta_2)|<\varepsilon/2
\ee
and then for fixed $t=t_\varepsilon$ take $N_t$ large enough that 
\be
|(-4t \log h^N_t(\theta_1,\theta_2))-(-4t \log h_t(\theta_1,\theta_2))|<\varepsilon/2
\ee
This would then give
us
\be
|d^2_{{\mathbb S}^1}(\theta_1,\theta_2)+4t \log h^{N}_{t}(\theta_1,\theta_2)|<\varepsilon 
\ee
for  $\theta_1, \theta_2 \in {\mathbb S}^1, \,\, t=t_{\varepsilon /2}, \textrm{ and } N=N{_{t_\varepsilon, \varepsilon /2}}$.
To get a better estimate, we need to take $t$ smaller first and then $N$ larger.   

\subsection{Eigenvalues and Eigenfunctions on a Rescaled Circle}\label{sect-recircle-evalue}

Let us consider a circle, ${{\mathbb{S}}^1_R}$, where $x=R\theta\in [0,2\pi R)$ with 
$
g_R=dx^2=R^2d\theta^2.
$
The Laplacian rescales 
\be
\Delta_{R}=\partial_ x^2= R^{-2} \partial_\theta^2
\ee
 and the volume is $2\pi R$.
We find that $\lambda_0=0$ with $\phi_0=(2\pi R)^{-{1}/{2}\,\,}$ by (\ref{phi-0}).
The rest of the eigenvalues come with multiplicity two:
\be
\lambda_{2j}=\lambda_{2j-1}=j^2/R^2
\ee
 with $L^2$ orthonormal eigenfunctions
\begin{eqnarray}
\phi_{2j-1}(\theta)&=&A_j\sin(j\theta) =A_j \sin(jx/R) \\
 \phi_{2j}(\theta)&=&A_j\cos(j\theta)=A_j\cos(jx/R)
 \end{eqnarray}
 where 
 \be
 1=\int_0^{2\pi} A_j^2 \sin^2(j\theta)\, R d\theta= \tfrac{1}{2}\, \int_0^{2\pi} A_j^2 \, R d\theta=2\pi A_j^2R/2 \textrm{ so }
 A_j=(\pi R)^{-\tfrac{1}{2}\,\,} 
 \ee
    So the Sturm-Liouville expansion for ${\mathbb S}^1_R$ is
 \be 
h^R_t(\theta_1,\theta_2)=\tfrac{1}{2\pi R} + \tfrac{1}{\pi R}
 \sum_{j=1}^\infty e^{-t j^2/R^2} \left( \sin(j\theta_1)\sin(j\theta_2) + \cos(j\theta_1)\cos(j\theta_2 ) \right).
\ee
Note that this rescales like (\ref{rescaled-heat}).
The trace of the heat kernel is then
\be
\int_{\theta=0}^{2\pi} h^R_t(\theta,\theta) \, R d\theta =  1 + 2\sum_{j=1}^\infty e^{-tj^2/R^2} .
\ee
Observe that as $R\to 0$ this converges to $0$ as expected.

\subsection{Eigenvalues and Eigenfunctions on Tori}\label{sect-tori-evalue}

Given any $A\le B$, the flat torus formed by taking the isometric product of rescaled circles,
\be
{\mathbb{S}}_A^1\times {\mathbb{S}}_B^1=({\mathbb{S}}^1\times {\mathbb{S}}^1, \, dx^2+dy^2),
\ee
where $x=A\theta$ and $y=B\varphi$ has the following Laplacian, 
\be
\Delta_{A,B}=\partial_ x^2 + \partial_ y^2 = A^{-2} \partial_\theta^2  + B^{-2} \partial_ \varphi^2.
\ee
Its volume is $4\pi^2AB$ so, by (\ref{phi-0}), for $\lambda_0=0$ the eigenfunction is 
$\phi_0=(4\pi^2AB)^{-\tfrac{1}{2}\,\,}$. 
The rest of the eigenfunctions come in sets of four as follows:
\be \label{torus-4}
\{\cos(j\theta)\cos(k\varphi), \cos(j\theta)\sin(k\varphi), 
\sin(j\theta)\cos(k\phi), \sin(j\theta)\sin(k\varphi)\}
\ee
which we can rescale by constants $A_{j,k}$ to be $L^2$ orthonormal.
Each set has the same eigenvalue
\be
\lambda_{j,k}=j^2/A^2+k^2/B^2.
\ee   
The ordering of these eigenvalues depends on the values of $A$ and $B$, so first lets examine the values of $A_{j,k}$.

When $j=0$ and $k=0$, we recover $\lambda_{0,0}=\lambda_0=0$ and see by (\ref{phi-0}) that there is
a single $L^2$-orthonormal eigenfunction,
\be 
\lambda_{0,0}=0 \textrm{ with }V_{0,0}=\{A_{0,0}\cos(0)\cos(0)\} 
\ee
$\textrm{ where } A_{0,0}= \tfrac{1}{\sqrt{(2\pi A)(2\pi B)}\,} =\tfrac{1}{2\pi\sqrt{AB}}$.

When $j>0$ and $k=0$, we also get $\lambda_{j,0}=j^2/A^2$ with two $L^2$-orthonormal eigenfunctions,
\be \label{Vj0}
\lambda_{j,0}=\tfrac{j^2}{A^2} \textrm{ with }
V_{j,0}=\{A_{j,0}\cos(j\theta)\cos(0), A_{j,0}\sin(j\theta)\cos(0)\} 
\ee
$\textrm{ where } A_{j,0}= \tfrac{1}{\pi \sqrt{2AB}\,}$
because
\be
1=\int_0^{2\pi}\int_0^{2\pi} A_{j,0}^2\sin^2(j\theta) AB\,d\theta d\varphi
=A_{j,0}^2(\tfrac{1}{2}\,) (2\pi A)(2\pi B).
\ee

When $j=0$ and $k>0$, we get $\lambda_{0,k}=k^2/B^2$ with two $L^2$-orthonormal eigenfunctions,
\be
\lambda_{0,k}=\tfrac{k^2}{B^2} \textrm{ with }
V_{0,k}=\{A_{0,k}\cos(0)\cos(k\varphi), A_{0,k}\cos(0)\sin(k\varphi)\} 
\ee
$\textrm{ where } A_{0,k}=  \tfrac{1}{\pi \sqrt{2AB}\,}.$

When both $j,k>0$, we get $\lambda_{j,k}=j^2/A^2+k^2/B^2$ with four $L^2$-orthonormal eigenfunctions as in (\ref{torus-4}),
\be
\lambda_{j,k}=\tfrac{j^2}{A^2}+\tfrac{k^2}{B^2} \textrm{ with }
V_{j,k}=\{A_{j,k}\cos(j\theta)\cos(k\varphi), ..., A_{j,k}\sin(j\theta)\sin(k\varphi)\}
\ee
$\textrm{ where } A_{j,k}= \tfrac{1}{\pi \sqrt{AB}\,}.$
because 
\begin{eqnarray}
1&=&\int_0^{2\pi}\int_0^{2\pi} A_{j,k}^2\sin^2(j\theta)\sin^2(k\varphi) AB\,d\theta d\varphi\\
&=&A_{j,k}^2 \,\cdot\, \int_0^{2\pi}\sin^2(j\theta) A\,d\theta\,\cdot\,   \int_0^{2\pi} \sin^2(k\varphi) B\, d\varphi\\
&=& A_{j,k}^2 \,\cdot\,(2\pi A (\tfrac{1}{2}\,))\,\cdot\,(2\pi B (\tfrac{1}{2}\,))=A_{j,k}^2\pi^2AB.
\end{eqnarray}

When $A=B=1$, the eigenvalues are ordered as follows with the eigenspaces: 
\begin{eqnarray}
\lambda_{0,0}=0^2+0^2=0 & V_{0,0} \textrm{ of multiplicity } 1\\
\lambda_{1,0}=\lambda_{0,1}=1^2+0^2=1 & V_{1,0}\cup V_{0,1} \textrm{ of multiplicity } 2+2=4\\
\lambda_{1,1}=1^2+1^2=2 & V_{1,1} \textrm{ of multiplicity } 4\\
\lambda_{2,0}=\lambda_{0,2}=2^2+0^2=4 & V_{2,0}\cup V_{0,2} \textrm{ of multiplicity } 2+2=4\\
\lambda_{2,1}=\lambda_{1,2}=2^2+1^2=5 & V_{2,1}\cup V_{1,2} \textrm{ of multiplicity } 4+4=8\\
\lambda_{2,2}=2^2+2^2=8 & V_{2,2} \textrm{ of multiplicity } 4\\
\lambda_{3,0}=\lambda_{0,3}=3^2+0^2=9 & V_{3,0}\cup V_{0,3} \textrm{ of multiplicity } 2+2=4 %\\
%\lambda_{3,1}=\lambda_{1,3}=3^2+1^2=10 & V_{3,1}\cup V_{1,3} \textrm{ of multiplicity } 4+4=8\\
%\lambda_{3,2}=\lambda_{2,3}=3^2+2^2=13 & V_{3,2}\cup V_{2,3} \textrm{ of multiplicity } 4+4=8\\
%\lambda_{4,0}=\lambda_{0,4}=4^2+0^2=16 & V_{4,0}\cup V_{0,4} \textrm{ of multiplicity } 2+2=4\\
%\lambda_{4,1}=\lambda_{1,4}=4^2+1^2=17 & V_{4,1}\cup V_{1,4} \textrm{ of multiplicity } 4+4=8\\
%\lambda_{3,3}=3^2+3^2=18 & V_{3,3} \textrm{ of multiplicity } 4
\end{eqnarray}

When $A=1$ and $B=\tfrac{1}{2}\,$, the eigenvalues $\lambda_{j,k}=j^2+4k^2$
are ordered quite differently with different multiplicities:
\begin{eqnarray}
\lambda_{0,0}=0^2+0^2=0 & V_{0,0} \textrm{ of multiplicity } 1\\
\lambda_{1,0}=1^2+4\cdot 0^2=1 & V_{1,0} \textrm{ of multiplicity } 2\\
\lambda_{0,1}=\lambda_{2,0}=0^2+4\cdot 1^2=4 & V_{0,1} \cup V_{2,0} \textrm{ of multiplicity } 2+2=4\\
\lambda_{1,1}=1^2+4\cdot 1^2=5 & V_{1,1} \textrm{ of multiplicity } 4\\
\lambda_{2,1}=2^2+4\cdot 1^2=8 & V_{2,1} \textrm{ of multiplicity } 4\\
\lambda_{3,0}=3^2+4\cdot 0^2=9 & V_{3,0} \textrm{ of multiplicity } 2 %\\
%\lambda_{3,1}=3^2+4\cdot 1^2=13 & V_{3,1} \textrm{ of multiplicity } 4\\
%\lambda_{0,2}=\lambda_{4,0}=4^2+4\cdot 0^2=16 & V_{0,2}=V_{4,0} \textrm{ of multiplicity } 2+2=4\\
%\lambda_{1,2}=1^2+4\cdot 2^2=17 &  V_{1,2} \textrm{ of multiplicity } 4\\
%\lambda_{2,2}=2^2+4\cdot 2^2=18 & V_{2,2} \textrm{ of multiplicity } 4\\
%\lambda_{4,1}=4^2+4\cdot 1^2=20 & V_{4,1} \textrm{ of multiplicity } 4\\
%\lambda_{5,0}=5^2+4\cdot 0^2=25 & V_{5,0} \textrm{ of multiplicity } 2\\
%\lambda_{3,2}=3^2+4\cdot 2^2=25 & V_{3,2} \textrm{ of multiplicity } 4\\
%\lambda_{0,3}=0^2+4\cdot 3^2=27 & V_{0,3} \textrm{ of multiplicity } 2
\end{eqnarray}

This ordering of the eigenfunctions becomes important when we try to apply these 
eigenfunctions to Manifold Learning.   

\subsection{Eigenvalues and Eigenfunctions on a Thin Torus}\label{sect-thin-evalue}

Let us consider a thin torus, ${\mathbb S}^1_A \times {\mathbb S}^1_B$ with $A=1$ and $B=1/10$.  Here  the eigenvalues are $\lambda_{j,k}=j^2+100k^2$.    Following our method from Section~\ref{sect-tori-evalue} we have the following
first 23 eigenvalues (counting multiplicity) and their eigenspaces:
\begin{eqnarray}
\lambda_{0,0}=0^2+0^2=0 & V_{0,0} \textrm{ of multiplicity } 1\\
\lambda_{j,0}=j^2+100\cdot 0^2=j^2&  V_{j,0} \textrm{ of multiplicity } 2 \textrm{ for }j=1,...,9\\
\lambda_{0,1}=\lambda_{10,0}=100 & V_{0,1}\cup V_{10,0} \textrm{ of multiplicity } 2+2=4.
\end{eqnarray}
 Note that  the eigenspaces $V_{j,0}$ 
as in (\ref{Vj0}) do
not have functions depending on $\varphi$.  Thus the truncated Sturm-Liouville expansion with $N=19$
for the
heat kernel of this thin torus does not depend on $\varphi$:
$$
h^{19}_t((\theta_1,\varphi_1),(\theta_2,\varphi_2))=\tfrac{1}{4\pi^2 B}\\  
+\, \tfrac{1}{2\pi^2 B}\sum_{j=1}^{9} e^{-j^2 t} 
(\sin(j\theta_1) \sin(j\theta_2)  +\cos(j\theta_1)\cos(j\theta_2)) 
$$
because $(1/(2\pi \sqrt{AB}))^2=1/(4\pi B^2)$ when $A=1$.

Note that this $N=19$ truncated Sturm-Liouville expansion for the
heat kernel of a thin torus is equal to the $N=19$ truncated Sturm-Liouville expansion for the
heat kernel of a circle, ${\mathbb S}^1$, that we found in (\ref{ST-S1}) rescaled by $1/(2\pi B)$.
{\em Intuitively the truncated heat kernel views this thin torus as almost the same
space as a circle. }   This becomes important when using eigenfunctions for Manifold Learning: sometimes it can be advantageous to view a thin torus as almost the same as a circle and other times it is better to distinguish them.

\subsection{Eigenvalues and Eigenfunctions on a Sphere}\label{sect-sphere-evalue}

The eigenfunctions of the sphere, ${\mathbb S}^2$ are called Laplace's spherical harmonics.  
They are explained beautifully in many sources including wikipedia:
\href{https://en.wikipedia.org/wiki/Spherical_harmonics#Laplace's_spherical_harmonics}{Spherical Harmonics}.
The eigenvalues take the form can be written as
\be
\ell(\ell+1) \textrm{ with eigenfunctions } Y_{\ell}^m=\cos( m \varphi)P_{\ell}^m(\cos(\theta))
\ee
where $m=-\ell, -\ell+1,...,0,...,\ell-1, \ell$ and $P_\ell^m$ is a Legendre Polynomial
where here we are using the coordinates $\varphi\in [0,\pi]$ and $\theta\in [0,2\pi]$ so that
\be
g_{{\mathbb S}^2}= d\varphi^2 + \sin^2(\varphi)d\theta^2.
\ee
The first eigenvalue $\lambda_0=0(1)=0$ with $\phi_0=(4\pi)^{-1/2}$ by (\ref{phi-0}).  The next eigenvalue is
$\lambda_1=1(2)=2$ with eigenfunction
\be \label{Z-sphere}
\phi_1(\theta,\varphi)=A_1 \cos(\varphi)=A_1 \,\, Z(\theta,\varphi)
\ee
where $Z(\theta,\varphi)$ is the $Z$ coordinate of the standard Riemannian isometric  
embedding of the sphere into Euclidean space,
\be \label{XYZ-sphere}
X(\theta,\varphi)=\sin(\varphi)\cos(\theta) \qquad Y(\theta,\varphi)=\sin(\varphi)\sin(\theta) 
\qquad Z(\theta,\varphi)=\cos(\varphi).
\ee
By symmetry we can see that $X(\theta,\varphi)$ and $Y(\theta,\varphi)$ also give eigenfunctions, so 
\be \label{XY-sphere}
\phi_2(\theta,\varphi)=A_1X(\theta,\varphi) \quad \textrm{ and } \quad
\phi_3(\theta,\varphi)=A_1Y(\theta,\varphi)
\ee
with $\lambda_2=\lambda_3=\lambda_1=2$.   To compute $A_1$, we can use the symmetry and the fact
that $X^2+Y^2+Z^2=1$ to see by (\ref{Z-sphere}) and (\ref{XY-sphere}) that
\begin{eqnarray}
1&=&\int_{{\mathbb S}^2 } \phi_1^2 \,dvol = \tfrac{1}{3} \int_{{\mathbb S}^2 } (\phi_1^2+\phi_2^2+\phi_3^2) \,dvol \\
&=&\tfrac{1}{3} \int_{{\mathbb S}^2 } A_1^2 (X^2+Y^2+Z^2) \,dvol
=\tfrac{1}{3} \int_{{\mathbb S}^2 } A_1^2 (1) \,dvol = \frac{4\pi A_1^2}{3}.
\end{eqnarray}
Thus $A_1=\sqrt{3/(4\pi)}$.    It is an exercise to check that these are orthonormal.

The $N=4$ truncated Sturm-Liouville expansion of the
heat kernel for ${\mathbb S}^2$ 
 can be viewed as a dot product of the embedding map (\ref{XYZ-sphere}):
\begin{eqnarray}
\qquad h_t^4(x_1,x_2) 
&=&\tfrac{1}{4\pi}+ \tfrac{3}{4\pi} e^{-2t } \left(X(x_1)X(x_2)+Y(x_1)Y(x_2)+Z(x_1)Z(x_2)\right)\qquad \quad\\
&=&\tfrac{1}{4\pi}+ \tfrac{3}{4\pi}  e^{-2t }
\cos\left(d_{{\mathbb S}^2}(x_1,x_2)\right).
\end{eqnarray}
This dependence on the distance is a consequence of
the symmetry of ${\mathbb S}^2$ and is unrelated Varadhan's formula because it
holds without taking $N\to \infty$ and $t\to 0$.

\subsection{Estimates on the First Eigenvalue}

In 1969, Cheeger found the following lower bound for the first nonzero eigenvalue of the Laplacian:

\begin{thm}
Let $M$ be a compact manifold with no boundary. Then,
\be \label{cheeger-ineq}
\lambda_1\ge \tfrac{1}{4} h^2(M).
\ee
Here $h(M)$ is Cheeger's Constant:
\be
h(M) = \inf\left\{ \frac{\vol_{n-1}(\partial U) }{ \min\{\vol_n(U), \vol_n(M\setminus U)\} }\right\}
\ee
where the infimum is taken over all domains $U \subset M$. 
\end{thm}
  
This was proven using Rayleigh's Quotient and 
coarea formulas in a clever way.   See, for example, the proof in \cite{Schoen-Yau-text}

In 1975 S-T Yau extended Cheeger's result and applied it to prove a lower bound on this eigenvalue in terms of the Ricci curvature on the manifold in \cite{Yau:1975}.  

\begin{thm}
When $\Ricci\ge (n-1)H$, then
\be
\lambda_1\ge \tfrac{1}{4} I^2(M) 
\ee
where 
\begin{eqnarray}
\quad I(M)^{-1} 
&\le& \omega_{n-1}\frac{\diam(M)}{\vol(M)} \int_0^{\diam(M)} \left(|H|^{-\tfrac{1}{2}\,\,}\sin(|H|^{\tfrac{1}{2}\,}r\right)^{n-1}\, dr \quad
\textrm{ when }H>0, \\
\quad I(M)^{-1} 
&\le& \omega_{n-1}\frac{\diam(M)}{\vol(M)} \int_0^{\diam(M)} (r)^{n-1}\, dr \qquad \qquad 
\textrm{ when }H=0,\\
\quad I(M)^{-1} 
&\le& \omega_{n-1}\frac{\diam(M)}{\vol(M)} \int_0^{\diam(M)} \left(|H|^{-\tfrac{1}{2}\,\,}\sinh(|H|^{\tfrac{1}{2}\,}r\right)^{n-1}\, dr \quad
\textrm{ when }H<0.
\end{eqnarray}
\end{thm}

In 1980 Buser proved Cheeger's inequality is sharp by showing that any manifold has a Riemannian metric tensor which achieves this inequality \cite{Buser:1980}.  
In 1982, Buser \cite{Buser:1982} proved the following uniform estimate:

\begin{thm}
If $M^n$ has $\Ricci\ge (n-1)H$ where $H<0$ then
\be
\lambda_1(M)\le c(n)\, \left(h(M)\,|H|^{1/2}+ h^2(M)\right).
\ee
\end{thm}

Buser presented two proofs in \cite{Buser:1982}, the second of which has a Jacobi field argument which is especially nice to read for those who enjoy do Carmo's textbook \cite{doCarmo:text}.   Buser also studied a variety of examples in these papers as well.

\subsection{Estimating Eigenvalues on Families of Manifolds}\label{sect-est-evalue}

Explicit computations of the eigenfunctions and eigenvalues for many Riemannian manifolds can be found in the French text of Berger, Gauduchon, and Mazet \cite{Berger-et-al:1971}.   See also Chavel for many examples with constant negative curvature \cite{Chavel-text}.

S.Y. Cheng gave upper estimates on all the eigenvalues of a Riemannian manifold with $Ricci \ge (n-1)H$ \cite{Cheng:1975}.  In 1980, Li and Yau improved this result, proving many results including  
\be
\lambda_k\le C_{H,D,n} \left(\tfrac{H+1}{V}\right)^{2/n} + C'_{H,n}
\ee
for any $n$ dimensional Riemannian manifold $M^n$ with $Ricci \ge (n-1)H$, $\diam(M)\le D$ and $\vol(M)\ge V$
\cite{Li-Yau:1980}.  

There are also lower bounds found by Cheng-Li in \cite{Cheng-Li:1981} depending on the Sobolev constant.   See also
the work of Debiard in \cite{Debiard-et-al:1976} and of Hess-Schrader-Uhlenbrock  in \cite{Hess-Schrader-Uhlenbrock:1980}.
For more details see the text of Schoen and Yau \cite{Schoen-Yau-text}.

\section{\bf Classes of Riemannian Manifolds and Convergence}

In Sections~\ref{sect-est-heat} and~\ref{sect-est-evalue} we have seen that entire classes of families of Riemannian manifolds have uniform bounds on their heat kernels and eigenvalues.   This is closely related to the various notions of convergence of Riemannian manifolds, their corresponding compactness theorems, and how the heat kernels and eigenvalues behave on converging sequences.   To be more clear, suppose we have a notion of convergence of Riemannian manifolds, $M_j \to M_\infty$, then we say a class, ${\mathcal{M}}$, of Riemannian manifolds is {\em compact} with respect to this notion of convergence if the following holds:
\be \label{class-compact}
M_j \subset {\mathcal{M}} \implies \exists \textrm{ a subseq }M_{j_k} \textrm{ s.t. } M_{j_k}\to M_\infty \in {\mathcal{M}}.
\ee
The class $\mathcal{M}$ is only {\em precompact} if the subsequence converges to $M_\infty$ that is not in ${\mathcal{M}}$.   

For many notions of convergence, $M_j \to M_\infty$, given any pair of points in the limit space, $p,q\in M_\infty$,
there are $p_j,q_j\in M_j$ that have $p_j\to p_\infty$ and $q_j\to q_\infty$.
We say a notion of convergence {\em preserves distances} if $M_j \to M_\infty$ implies
\be\label{preserves-distances}
 d_j(p_j,q_j)\to d_\infty(p,q).
\ee 
We say a notion of convergence
{\em preserves the heat kernel} if $M_j \to M_\infty$ implies
\be \label{preserves-heat-kernel}
 \forall t>0,\,\,h^{M_j}_{t}(p_j,q_j)\to h^{M_\infty}_{t}(p,q). 
\ee 
We say there is {\em continuity of the eigenvalues} if $M_j \to M_\infty$ implies
\be \label{continuity-evalues}
 \forall k \in {\mathbb N} \,\, \lambda_k(M_j) \to \lambda_k(M_\infty).
\ee
Sometimes we only have semicontinuity of the eigenvalues.

If we have a family of manifolds, ${\mathcal{M}}$, and a notion of convergence which makes this family compact and this notion of convergence implies continuity of the eigenvalues, then we
have uniform bounds on the eigenvalues.   Otherwise, we could take a sequence of $M_j\in \mathcal{M}$ whose eigenvalues are diverging, find the converging subsequence and use continuity of the eigenvalues on the subsequence to reach a contradiction.  Similarly, if the notion of convergence preserves the heat kernel, then it has uniform bounds on the heat kernel.

Below we will discuss various notions of convergence and their relationship with heat kernels and eigenvalues.

\subsection{Smooth, $C^0$, $C^{k,\alpha}$, and $C^k$  Convergence}

Recall that a function $h: U\subset {\mathbb R}^n \to {\mathbb R}$ is in $C^0(U)$ if it is  continuous  and the $C^0$ norm is
\be
|h|_{C^0}= \sup_{x\in U} |h(x)|.
\ee
We say that $h\in C^{0,\alpha}(U)$ if it's Holder norm is bounded:
\be
|h|_{C^{0,\alpha}(U)}=\sup_{x\neq y \in U} \frac{|h(x)-h(y)|}{|x-y|^\alpha}<\infty.
\ee
We say that $h\in C^k(U)$ if 
\be
|h|_{C^{k}(U)}=\sup_{x\in U} \sum |\partial_{x_1}^{\beta_1}\cdots \partial_{x_n}^{\beta_n} h(x)|<\infty,
\ee
where the sum is over $(\beta_1,...,\beta_n)$
such that $\beta_1+\cdots+\beta_n\le k$. 
We say that $h\in C^{k,\alpha}(U)$ if
\be
|h|_{C^{k,\alpha}(U)}=|h|_{C^{k}(U)}+  \sum |\partial_{x_1}^{\beta_1}\cdots \partial_{x_n}^{\beta_n} h(x)|_{C^{0,\alpha}}(U)<\infty,
\ee
where the sum is over $(\beta_1,...,\beta_n)$ such that $\beta_1+\cdots+\beta_n=k$.  

We say an atlas, $\{\varphi_i\}$, of charts on $M$ is $C^k$ or $C^{k,\alpha}$ if the transition
maps are $C^k$ or $C^{k,\alpha}$ respectively.
If we fix a $C^k$ or $C^{k,\alpha}$ atlas, $\{\varphi_i\}$, of charts on $M$, then with respect to that atlas we can define:
\be
|h|_{C^k(M),\{\varphi\} }=\sup |h\circ \varphi |_{C^k(U)}\quad 
\textrm{ or }\quad  |h|_{C^{k,\alpha}(M),\{\varphi\} }=\sup |h\circ \varphi |_{C^{k,\alpha}(U)}
\ee
where the sup is over all coordinate charts $\varphi: U \to M$ in the atlas.  

A function between compact manifolds, $H: M_1\to M_2$ is $C^k$, if there exists a $C^k$ atlas of charts on $M_1$ and on $M_2$ such that the function composed with these charts is $C^k$ respectively.    It is a $C^k$ diffeomorphism if it is invertible and both it and its inverse are $C^k$.  It is a $C^{k,\alpha}$ diffeomorphism if everything here is $C^{k,\alpha}$.
It is a smooth or $C^\infty$ diffeomorphism if it is $C^k$ for all $k \in {\mathbb N}$.  

We say $(M_j, g_j) \to (M_\infty,g_\infty)$ in the $C^k$ sense if $H_j: M_\infty \to M_j$ is a $C^{k+1}$ diffeomorphism (which means it has $(k+1)$ continuous derivatives) and if there exists a $C^k$ atlas on $M_\infty$ such that
for all $g_\infty$-unit length $ u,v \in TM_\infty$ we have
\be
|H_j^*g_j(u,v)- g_\infty(u,v)|_{C^k(M), \{\varphi_i\}}<\epsilon_j \to 0.
\ee
We say $(M_j, g_j) \to (M_\infty,g_\infty)$ in the $C^{k,\alpha}$ sense if $H_j: M_\infty \to M_j$ is a $C^{k+1,\alpha}$ diffeomorphism  and if there exists a $C^{k,\alpha}$ atlas on $M_\infty$ such that
for all $g_\infty$-unit length $ u,v \in TM_\infty$ we have
\be
|H_j^*g_j(u,v)- g_\infty(u,v)|_{C^{k,\alpha}(M),\{\varphi_i\}}<\epsilon_j \to 0.
\ee
It converges smoothly in the $C^\infty$ sense if all the diffeomorphisms are $C^\infty$ and if there exists a 
$C^\infty$ atlas on $M_\infty$ such that
for all $g_\infty$-unit length $ u,v \in TM_\infty$ we have
\be
|H_j^*g_j(u,v)- g_\infty(u,v)|_{C^{\infty}(M),\{\varphi_i\}}<\epsilon_j \to 0.
\ee

\subsection{Smooth Convergence of Tori and Warped Tori}

Note that if we consider tori ${\mathbb S}^1_{A_j}\times {\mathbb S}^1_{B_j}$ and $A_j\to A_\infty>0$
and $B_j \to B_\infty>0$ then we have $C^\infty$ convergence to ${\mathbb S}^1_{A_\infty}\times {\mathbb S}^1_{B_\infty}$.   This can be seen by taking $\theta, \varphi\in [0,2\pi)$ and writing the metric tensors,
\be
g_{{\mathbb S}^1_{A}\times {\mathbb S}^1_{B}}=A^2 d\theta^2 + B^2 d\varphi^2,
\ee
and taking $H_j: {\mathbb S}^1_{A_\infty}\times {\mathbb S}^1_{B_\infty}\to {\mathbb S}^1_{A_j}\times {\mathbb S}^1_{B_j}$ preserving $\theta$ and $\varphi$.   Then
\be
|H_j^*g_j(u,v)- g_\infty(u,v)|_{C^0}=|(A_j^2-A_\infty^2)u(\theta)v(\theta)+(B_j^2-B_\infty^2)u(\varphi)v(\varphi)|_{C^0}\to 0
\ee
It is easy to see there is $C^k$ convergence as well.

A sequence of warped product metric tensors on a torus converges in the $C^k$ sense
\be\label{warp-Ck}
g_j=d\varphi^2 +f_j(\varphi)^2 d\theta^2 \,\,\to\,\, g_\infty=d\varphi^2 +f_\infty(\varphi)^2 d\theta^2
\ee
if and only if  $|f_j-f_\infty|_{C^k} \to 0$.  On tori we must assume the warping factors are periodic,
$f_j(\varphi)=f_j(\varphi+2\pi)$, and positive,
$f_j>0$ and $f_\infty>0$, to have a diffeomorphism.   If we only have $f_j(0)=f_j(\pi)$ without requiring the
derivatives to match as implied by periodicity, then we only have a $C^0$ metric tensor and $C^0$ convergence.  To have
a $C^k$ metric tensor, we need $k$ derivatives to match:
\be
f_j(0)=f_j(\pi)\quad f_j'(0)=f_j'(\pi)\quad f_j''(0)=f_j''(\pi)\quad ...\quad f_j^{(k)}(0)=f_j^{(k)}(\pi).
\ee
To have $C^k$ convergence, we need all these derivatives to converge 
$$
f_j^{(i)}(0)\to f_\infty^{(i)}(0) \textrm{ for } i=0,...,k.
$$

A sequence of warped product metric tensors on a 
sphere will have (\ref{warp-Ck}) with $f_j(0)=f_j(\pi)=0$ and $f_\infty(0)=f_\infty(\pi)$.   To have
a $C^k$ metric tensor across the poles, we need $k$ derivatives to match that of the standard sphere:
\begin{eqnarray}
f_j(0)&=&f_j(\pi)=\sin(0)=0\\
 f_j'(0)&=&f_j'(\pi)=\sin'(0)=1\\
 f_j''(0)&=&f_j''(\pi)=\sin''(0)=0\\
&...&\\
f_j^{(k)}(0)&=&f_j^{(k)}(\pi)=\sin^{(k)}(0).
\end{eqnarray}
To have $C^k$ convergence, these derivatives must converge 
$$
f_j^{(i)}(0)\to f_\infty^{(i)}(0) \textrm{ for } i=0,...,k.
$$

\subsection{Convergence of Distances between Points}

Distances are preserved under $C^0$ convergence of  $M_j \to M_\infty$.   To see this, first we show 
that for any curve $C:[0,1]\to M$ such that $C(0)=p_\infty$ and $C(1)=q_\infty$ we have
$C_j=C\circ H_j^{-1}:[0,1]\to M_j$ with $C_j(0)=p_j$ and $C_j(1)=q_j$.  Since
\be
L_{g_j}(C_j)=\int_0^1 g_j(C_j'(s),C_j'(s))^{1/2}\,ds=\int_0^1 (H_j^*g_j)(C'(s),C'(s))^{1/2}\, ds \to L_{g_\infty}(C).
\ee
we see that 
$
| L_{g_j}(C_j)-L_{g_\infty}(C)|<\epsilon_j^{1/2}.
$
Since the distance $d_j$ is the infimum of the length, $L_j$, of all curves from $p_j$ to $q_j$, we have
\be
d_j(p_\infty,q_\infty) \le L_{g_j}(C_j)\le L_{g_\infty}(C)+\epsilon_j^{1/2}.
\ee
Since the distance, $d_\infty$, is the infimum of the length, $L_\infty$, over all curves, $C$, as above
\be
d_j(p_\infty,q_\infty) \le d_\infty(p_\infty,q_\infty) +\epsilon_j^{1/2}
\ee
Applying this entire argument in reverse we get the desired convergence:
\be
|d_\infty(p_\infty,q_\infty) -d_j(p_j,q_j)|\le\epsilon_j^{1/2} \to 0.
\ee

\subsection{Lipschitz Convergence}\label{Lip-Conv}

Lipschitz convergence of Riemannian manifolds is defined by viewing a Riemannian manifold, $(M,g)$,
as a metric space $(M, d_g)$ where $d_g$ is the Riemannian distance.   We say two Riemannian Manifolds
are biLipschitz if there exists a bijection $H:M_1 \to M_2$ such that
\begin{eqnarray}
\Lip(H)\,\,\,&=&\sup_{x\neq y} d_{2}(H(x),H(y))/d_{1}(x,y) <\infty \\
\Lip(H^{-1})&=&\sup_{x\neq y} d_{1}(H^{-1}(x),H^{-1}(y))/d_{2}(x,y) <\infty.
\end{eqnarray}
The Lipschitz distance between them is:
\be
d_{Lip}((M_1,d_1),(M_2,d_2))= \log \min\{\Lip(H), \Lip(H^{-1})\, : \, H:M_1\to M_2\}.
\ee
It is easy to verify that $M_j\to M_\infty$ in the $C^0$ sense implies $d_{Lip}(M_j,M_\infty)\to 0$.

Consider the class of smooth Riemannian manifolds which are biLipschitz to a standard sphere, ${\mathbb S}^2$ with Lipschitz bounds in both directions $\le 2$.   By the Arzela-Ascoli Theorem, this class is precompact with respect to Lipschitz convergence.  However, the limit spaces need not be a smooth Riemannian manifolds.   

In fact $M_\infty$ could have a conical singularity.   This can be seen by considering the following sequence of functions:
\be
g_j=\,d\varphi^2 + f_j(\varphi)^2\,d\theta^2
\ee
where $f_j$ are smooth functions such that 
\begin{eqnarray}
f_j(\varphi)&=&\sin(\varphi) \textrm{ on } [0,\delta_j]\cup [\pi-\delta_j,\pi]\\
f_j(\varphi)&=& \sin(\varphi)/2  \textrm{ on } [0,5\delta_j]\cup [\pi-5\delta_j,\pi]
\end{eqnarray}
and monotone on $[\delta_j,5\delta_j]\cup[\pi-5\delta_j,\pi-\delta_j]$.
Taking $\delta_j \to 0$, the Lipschitz limit is an American-football-shaped $C^0$ Riemannian manifold with conical singularities at the poles where $\varphi=0$ and $\varphi=\pi$.

This sequence has Ricci and sectional curvature bounded from below by $0$ but not from above and has no lower bound on the injectivity radius.   However we can still apply Cheeger-Yau's theorem in (\ref{Cheeger-Yau-heat}) to say that
the heat kernels of these manifolds are bounded from below by the heat kernel on Euclidean space.   On balls which avoid the poles, the heat kernels are bounded on both sides.  In fact as we shall soon see, the heat kernels converge to a function on the limit space which has many of the properties of a heat kernel.

\subsection{$C^2$ convergence preserves $\Delta$ and $h_t(x,y)$}

If $(M_j, g_j) \to (M_\infty,g_\infty)$ in the $C^2$ sense, then it is easy to see that their curvatures and Laplacians 
are preserved,
\be
|\, \Delta^{M_j} (h\circ H_j^{-1}) - \Delta^{M_\infty} h \,|_{C^0} \to 0 \textrm{ for any } C^2 \textrm{ function }h: M_\infty \to {\mathbb R}.
\ee 
because Laplacians depend only on first and second derivatives of the metric tensor.

Since solutions to the heat equation depend smoothly on the Laplace operator, the heat kernel on $M_j$ can be shown to converge to the heat kernel on $M_\infty$,
\be
| h^{M_j}_t(H_j^{-1}(x), H_j^{-1}(y))- h^{M_\infty}_t(x,y) |_{C^2} \to 0,
\ee
and similarly the eigenfunctions can be shown to converge
\be
|\phi^{M_j}_k(H_j^{-1}(x))-\phi^{M_\infty}_k(x)|_{C^2}\to 0 \textrm{ and }|\lambda_k(M_j)-\lambda_k(M_\infty)|\to 0
\ee
where as usual the $\lambda_k$ are counted with multiplicity.   

What is quite surprising is that heat kernels and eigenvalues can converge even with far weaker notions of convergence.

\subsection{Diffeomorphic Compactness Theorems}

With the above definitions of convergence of manifolds we require that the sequence of manifolds be diffeomorphic in order to compare their metric tensors.   So if we are studying a class or family of Riemannian manifolds, ${\mathcal M}$,
and wish to prove a compactness theorem stating that given a sequence $M_j\in {\mathcal M}$ a subsequence converges to $M_\infty \in {\mathcal M}$, we first need to show a subsequence is pairwise diffeomorphic to one another.   This can be achieved using the pigeonhole principle if we know there are only finitely many diffeomorphism types in the
class $\mathcal{M}$.

In his 1967 doctoral dissertation \cite{Cheeger:1967}\cite{Cheeger:1970}, Cheeger observed that there were only finitely many diffeomorphism types in
\be
{\mathcal M}_{K_1, V}^{K_2, D}=\{M^n\, |\, K_1\le \sect\le K_2, \vol(M)\ge V, \diam(M)\le D\}
\ee
by proving a uniform lower bound on injectivity radius based on $K_1$, $K_2$, $V$, and $D$.   He also found other classes of Riemannian manifolds with similar controls.   

Note that the tori 
\be
{\mathbb S}^1_A\times {\mathbb S}^1_B \in {\mathcal M}_{K_1, V}^{K_2, D} \textrm{ for any } K_1\le 0\le K_2
\ee
if and only if
\be \label{ABch}
4\pi^2AB \ge V \textrm{ and }(\pi A)^2+ (\pi B)^2 \le D^2
\ee
because they are flat with sectional curvature $=0$.   Notice that $(\ref{ABch})$ implies $A, B \le D/\pi$
which implies $A,B \ge V/(4\pi D)>0$.   So we do not have arbitrarily thin tori with $B_j \to 0$ in this class.

Once Cheeger had a uniform lower bound on injectivity radius, he proved 
\be
{\mathcal M}_{K_1, V}^{K_2, D} \subset {\mathcal M}_{K_1 i_0}^{K_2, D}
\ee
where 
\be
{\mathcal M}_{K_1 i_0}^{K_2, D}=\{M^n\, |\, K_1\le \sect\le K_2, injrad\ge i_0, \diam(M)\le D\}.
\ee
Note that any manifold $M^n$ in this second family
can be covered by at most $N_{K_1 i_0}^{K_2, D}$ charts as follows.   About any point $p\in M^n$, we can take a chart defined by the exponential map, $\exp_p$, which covers a ball about $p$ of radius $r=i_0$.   To choose a finite collection, we take a maximal collection of disjoint balls of radius $r/2$, so that the balls of radius $r$ cover $M^n$.   Let
\be
N(M,r)=\textrm{ max number of disjoint balls of radius $r/2$ in $M$}.
\ee
To find a uniform bound on this number, 
note that $Vol(B(p,r/2))\ge V^n_{K_2}(r/2)$ which is the volume of a ball is a simply connected space of constant sectional curvature $K_2$, so
\be
N(M,r) V^n_{K_2}(r/2) \le \sum_{i=1}^{N(M,r)} \vol(B(p_i,r/2)) \le \vol(M) \le V^n_{K_1}(D).
\ee
Thus
\be
N(M,r)\le N(r)=V^n_{K_1}(D)/V^n_{K_2}(r/2).
\ee
Note that the lower bound on the $Vol(B(p,r/2))$ here strongly uses the fact that $r\le i_0\le \injrad(p)$.   Clearly such a uniform lower bound on volume does not hold on increasingly thin tori with $B_j \to 0$.

Once we have a uniform bound on the number of charts, we see that we have only finitely many homeomorphism types bounded by all the possible ways these charts can overlap.   If we have a sequence of $M_j$ in a class ${\mathcal M}$ with a uniform bound $N(M,r)\le N(r)$ then we can apply the pigeon hole principle to choose a subsequence, $M_j$, such that they are all homeomorphic and such that $N(M_j,r)$ is constant and such that the same pairs of balls of radius $r$ intersect.   
To complete the proof of Cheeger's Compactness Theorem, we would then only need to show uniform bounds on the sequence of metric tensors within these charts.   We recommend Petersen's textbook for a discussion of the many ways that 
metric tensors can be bounded within the charts \cite{Petersen-text}.   See also Gromov's textbook \cite{Gromov-text}.

\subsection{Injectivity Radius, Heat Kernels, and Eigenvalues}
We have seen above in Sections~\ref{sect-est-heat} and ~\ref{sect-est-evalue} there are entire classes of manifolds which share the same bounds on their heat kernels and on their eigenvalues.  
 In particular there
was the class
\be
{\mathcal M}^n_H,i_0=\{M^n\, | \,\Ricci_M\ge (n-1)H \textrm{ and } injrad_M \ge i_0\}.
\ee
The injectivity radius bound here is essential.  
The Ricci curvature bound in also a key step.   

Consider for example the rescaled circles ${\mathbb S}^1_R$ that we studied 
in Section~\ref{sect-recircle-heat}
and~\ref{sect-recircle-evalue}.   Their injectivity radii 
\be
\injrad_{{\mathbb S}^1_R}=\pi R \to 0 \textrm{ as }R \to 0.
\ee
We saw in Section~\ref{sect-recircle-heat} that their heat kernels,
\be
h^R_t(\theta_1,\theta_2)= R^{-1} \,h_{t/R^2}(\theta_1,\theta_2)
\ee
are not uniformly bounded above for sequences as $R\to 0$.
We saw in Section~\ref{sect-recircle-evalue} that their eigenvalues diverge
\be
 \lambda_{2j}=\lambda_{2j-1}=j^2/R^2 \to \infty \textrm{ as } R\to 0.
\ee
Sequences of rescaled spheres, ${\mathbb S}_R^2$, have nonnegative Ricci curvature
and the same problems with their injectivity radii, heat kernels, and eigenvalues
as the ${\mathbb S}_R^1$ as $R\to 0$.

\subsection{Harmonic Coordinates and Harmonic Radius}

In \cite{DeTurck-Kazdan:1981},
DeTurck and Kazdan proved that whenever an atlas exists covering $M$ with charts such that 
$h\circ\varphi_i \subset C^k(U_i)$, then the same is true if we chose the atlas of {\bf \em harmonic coordinate charts}.
A coordinate chart $\varphi: U\subset {\mathbb R}^m \to M$ is {\em harmonic} if the components of $\varphi^{-1}$ are harmonic functions:
\be
\varphi^{-1}(p)=(x_1(p),x_2(p),...,x_n(p)) \textrm{ where } \Delta x_i=0.
\ee
A manifold $M$ has {\bf {\em harmonic radius}} $r_h>0$ if it can be covered by a collection of harmonic coordinate charts, each of which covers a ball of radius $r_h$. 

In 1990, Anderson \cite{Anderson:1990} proved a $C^{1,\alpha}$ precompactness theorem for
\be
{\mathcal{M}}^D_{H, i_0}=\{ M^n\,:\, |\Ricci|\le (n-1)H, \injrad(M)\ge i_0>0, \diam(M)\le D\}
\ee
using harmonic coordinate charts by finding a uniform lower bound on the harmonic radius, which
is a lower bound on the radius of harmonic coordinate charts.   Building upon this work, 
Anderson-Cheeger proved a $C^{\alpha}$
compactness theorem for a larger class of manifolds in \cite{Anderson-Cheeger:1992}.

The harmonic radius can also be used to control eigenfunctions as is seen for example in the
work of Portegies  \cite{Portegies:2016} that we will present at the end of this paper.

\subsection{Gromov-Hausdorff Convergence} \label{sect-GH}

In 1983 Gromov introduced the notion of the Gromov-Hausdorff distance
between compact metric spaces \cite{Gromov-text}.   Since a 
compact Riemannian manifold, $(M,g)$
 can be viewed 
as metric space, $(M, d_g)$, his notion allowed him to define the Gromov-Hausdorff convergence
of Riemannian manifolds to limits which may only be metric spaces.    The notion
allowed him to study sequences of collapsing manifolds (like increasingly thin tori) whose limits have lower dimension
as well as sequences of noncollapsing manifolds.

Before we define Gromov-Hausdorff distance
we review that the {\bf Hausdorff distance} between two compact sets, $K_1$ and $K_2$, in a metric space, $Z$,
is defined
\be\label{Hausdorff}
d_H^Z(K_1, K_2) =\inf\{R:\, K_1\subset T_R(K_2) \textrm{ and } K_2\subset T_R(K_1)
\ee
where $T_R(K_i)$ is the tubular neighborhood about $K_i$:
\be
T_R(K_i)=\{z\in Z:\, \exists x\in K_i \textrm{ s.t. } d_Z(x,z)<R\}.
\ee
   This is an
extrinsic way of measuring the distance between $K_1$ and $K_2$ which depends very
much on the extrinsic space, $Z$, in which they both sit.

To define an intrinsic distance between two compact metric spaces, $(X_1, d_1)$ and $(X_2, d_2)$,
which are not sitting in a common metric space, Gromov considered all distance preserving maps,
\be
F_i: X_i \to Z \textrm{ such that } d_Z(F_i(p), F_i(q))=d_{X_i}(p,q) \qquad \forall p,q\in X_i
\ee
and all possible compact metric spaces, $(Z, d_Z)$.   He defined the
{\bf Gromov-Hausdorff distance} by taking the infimum over all such possibilities:
\be
d_{GH}(X_1, X_2))=\inf d_H^Z(F_1(X_1), F_2(X_2)).
\ee
He called this an intrinsic Hausdorff distance because it depends only on 
intrinsic information about the pair of metric spaces, $(X_1, d_1)$ and $(X_2, d_2)$,
and is defined using the Hausdorff distance.

Gromov proved that if a sequence $(X_j, d_j \to (X_\infty,d_\infty)$ in the Gromov-Hausdorff sense
then there exists a common compact metric space, $Z$, and distance preserving maps
$F_j: X_j \to Z$ such that 
\be
d_H^Z(F_j(X_j), F_\infty(X_\infty)) \to 0.
\ee
In particular, for all $p_\infty\in X_\infty$ there exists $p_j\in X_j$ such that
$F_j(p_j)\to F_\infty(p_\infty)$.   Although this is not uniquely assigned, this is
how  
\be\label{conv-pts}
p_j\in X_j \textrm{ converges to }p_\infty\in X_\infty
\ee
is defined.

\noindent
{\bf Gromov's  Compactness Theorem}: {\em A sequence of Riemannian
manifolds 
\be
M_j^n\in {\mathcal M}^n_{H,D}=\{M^n\, | \,\,\Ricci_{M^n}\ge (n-1)H,\, \diam(M^n)\le D \}.
\ee
has a subsequence which converges in the Gromov-Hausdorff sense to a compact
metric space. }

Gromov proved this in \cite{Gromov-text}  by first proving the following theorem:

\noindent
{\bf Bishop-Gromov Volume Comparison
Theorem:}  {\em for any $M^n\in {\mathcal M}^n_{H,D}$, and any
$p\in M^n$ and any $R>r>0$, the volumes of balls satisfy the following:
\be
\vol(B(p,r))/\vol(B(p,R)) \ge V^n_H(r) / V^n_H(R)
\ee
where $V^n_H(s)$ is the volume of a ball of radius $s$ in an $n$ dimensional simply connected
manifold of constant sectional curvature, $H$. }  

This allowed Gromov to uniformly count the number, $N(r)$, of balls of any
given radius in $M_j^n\in {\mathcal M}^n_{H,D}$
space.  He applied this to create uniformly dense countable collection points, $p_{j,i}\in M_j$
and matrices of distances between these points, 
\be
d_{j,i,k}=d_{M_j}(p_i,p_k) \in [0,D].
\ee
For each, $(i,k)$, a subsequence of $j$ (that we still denote by $j$) has
\be
d_{j,i,k}\to d_{\infty,j,k} \in[0,D].
\ee
Gromov diagonalized the sequence to obtain this convergence for all pairs of points
and then built a limit space out of the countable collection of points assigning them
\be
X=\{p_j\,:\,j\in {\mathbb N}\} \textrm{ and }d_{\infty}(p_j,p_k)=d_{\infty,j,k}.
\ee
He used the triangle inequality and the uniform density of the points to take a metric completion,
$\bar{X}$, and prove $M_j\to \bar{X}$ in the Gromov-Hausdorff sense for this subsequence.
See \cite{Gromov-text} and see also the text of Burago-Burago-Ivanov \cite{BBI-text}).

\subsection{Increasingly Thin Tori GH-Converge to Circles}
An example of a collapsing sequence of $M_j\subset  {\mathcal M}^n_{H,D}$
which converges in the Gromov-Hausdorff sense
is the sequence of increasingly thin tori
\be
{\mathbb S}^1_A \times {\mathbb S}^1_{1/j} \GHto {\mathbb S}^1_A.
\ee
This can be seen by observing that for any fixed $j$, we can take 
\be
Z={\mathbb S}^1_A \times {\mathbb S}^1_{1/j}
\ee
and distance preserving maps
\be
F_1: {\mathbb S}^1_A \times {\mathbb S}^1_{1/j} \to Z \textrm{ to be } F_2(\theta,\varphi)=(\theta,\varphi).
\ee
and 
\be
F_2: {\mathbb S}^1_A\to Z \textrm{ to be } F_2(\theta)=(\theta,0).
\ee
The Hausdorff distance between 
\be
F_2( {\mathbb S}^1_A)=\left\{(\theta,0):\theta\in [0,2\pi]\right\} \subset Z \textrm{ and } F_1( {\mathbb S}^1_A \times {\mathbb S}^1_{1/j})=Z
\ee
is easily seen to be $\pi/j$.   Thus
\be
d_{GH}\left({\mathbb S}^1_A \times {\mathbb S}^1_{1/j}, {\mathbb S}^1_A\right)<\pi/j \to 0.
\ee

In Section~\ref{sect-thin-evalue}, we saw that for each $i$ the $i^{th}$ eigenfunctions
of increasingly thin tori ${\mathbb S}^1_A\times {\mathbb S}^1_B$ converge to the
$i^{th}$ eigenfunctions of the circle, ${\mathbb S}^1_A$.
Thus the truncated heat kernels of thin tori ${\mathbb S}^1_A\times {\mathbb S}^1_B$ 
rescaled by $(2\pi B)$ converge as  $B\to 0$
 to
$$
\tfrac{1}{2\pi} 
+\, \tfrac{1}{\pi}\sum_{j=1}^{N} e^{-j^2 t} 
(\sin(j\theta_1) \sin(j\theta_2)  +\cos(j\theta_1)\cos(j\theta_2)) 
$$
which
is the truncated heat kernel on a circle ${\mathbb S}^1_A$.   Indeed the full heat kernel of 
${\mathbb S}^1_A\times {\mathbb S}^1_B$
rescaled by $(2\pi B)$ converges
to the heat kernel of  ${\mathbb S}^1_A$ as $B\to 0$.

\subsection{Eigenvalues and Measured Gromov-Hausdorff Convergence}\label{sect-mGH}

This nice behavior of the eigenvalues and heat kernels is not true in general.  
Fukaya found examples of sequences of Riemannian manifolds which
converge in the Gromov-Hausdorff sense to lower dimensional manifolds such that the eigenfunctions
did not converge the eigenvalues of the limit manifold  \cite{Fukaya:1987aa}.   More precisely, the
eigenvalues, $\lambda_k(M_j)$, do not converge
to eigenvalues, $\lambda_k(M_\infty)$ of the standard Laplace operator, $\Delta_{M_\infty}$, on the limit manifold,
$M_\infty$.

Fukaya studied collapsing
warped tori:
\be
(M_j,g_j)=({\mathbb S}^1 \times {\mathbb S}^1, d\theta^2 + f^2(\theta)/j^2 d\phi^2
\ee
where $f(\theta)$ is a nonconstant positive periodic warping function, which converge in the Gromov-Hausdorff
sense to the standard $({\mathbb S}^1, d\theta^2)$.   Since the Laplacian on $(M_j,g_j)$ has the form
\be
\Delta_{M_j}=1\frac{1}{f(\theta)}\frac{\partial}{\partial \theta} f(\theta) \frac{\partial}{\partial \theta}
-\frac{1}{j^2f^2(\theta)} \frac{\partial^2}{\partial \varphi^2},
\ee
Fukaya explained that their eigenvalues, $\lambda_k=\lambda_k(M_j)$ converge as $j \to \infty$, to the eigenvalues, 
$\lambda_k'$, of the operator
\be
1\frac{1}{f(\theta)}\frac{d}{d \theta} f(\theta) \frac{d}{d \theta} 
\ee
on ${\mathbb{S}}^1$.    This operator is not equal to the Laplace operator, $\Delta_{{\mathbb S}^1}=d^2/d\theta^2$,
and has different eigenvalues, $ \lambda_k'\neq \lambda_k({\mathbb S}^1)$.

In \cite{Fukaya:1987aa}, Fukaya proposed the notion of {\bf measured Gromov-Hausdorff convergence} of
metric measure spaces:
\be
(X_j, d_j, \mu_j) \to (X_\infty, d_\infty, \mu_\infty)
\ee
where $(X_j, d_j) \to (X_\infty, d_\infty)$ in the Gromov-Hausdorff sense and if
$p_j\in X_j$ converges to $p_\infty\in X_\infty$ as in (\ref{conv-pts}) then the volumes
of balls about them converge:
\be
\mu_j(B(p_j,r))\to \mu_\infty(B(p_\infty,r)) \textrm{ for all } r>0.
\ee

Fukaya proved that with this notion of metric measure convergence, and a natural kind of Laplacian defined in a way
that depends on the measure, there is semicontinuity of
the eigenvalues:
\be
\limsup_{j\to \infty} \lambda_k(X_j, d_j, \mu_j) \le \lambda_k(X_\infty, d_\infty, \mu_\infty).
\ee
Fukaya conjectured that eigenvalues and other spectral properties
will behave well under metric measure convergence if the sequence is in ${\mathcal M}^n_{H,D}$
and he proved this conjecture in the case with two sided sectional curvature bounds.  He found
counter examples without the lower bound on Ricci curvature.

In \cite{ChCo-Part3}, Cheeger and Colding proved Fukaya's conjecture for sequences of 
$M_j \in {\mathcal M}^n_{H,D}$.   In fact, they proved that
any sequence $M_j \in {\mathcal M}^n_{H,D}$ endowed with measures defined using rescaled
volumes has a subsequence which converges in the metric measure sense to a
metric measure space with a measure that satisfies the Bishop-Gromov Inequality.   They defined
a natural notion of eigenfunction on this limit space and proved the eigenfunctions converge as well.
Building on their work, Yu Ding proved the heat kernels converge as well \cite{Ding:02}.

More recently Sturm has defined the metric measure distance between Riemannian manifolds is defined as the
infimum over all compact metric spaces $Z$ and all distance preserving maps $F_i: M_i \to Z$
of the Wasserstein distance between the push forwards of their measures into $Z$
\cite{Sturm-D}.  This notion is useful for studying $CD(H,n)$ spaces.  
These $CD(H,n)$ spaces were first introduced by
Lott and Villani and by Sturm using optimal transport \cite{Lott-Villani:2009}\cite{Sturm-D}
and were shown to include  the limits of 
$M_j \in {\mathcal M}^n_{H,D}$.
See also Villani's text \cite{Villani-text}.  
A stronger class of spaces called
$\mathrm{RCD}^*(H,n)$ spaces which were
 introduced by Ambrosio-Gigli-Savare in \cite{Ambrosio-Gigli-Savare:2014} have all the 
 properties of $CD(H,n)$ spaces and also satisfy a splitting theorem.   They
 proved all limits of manifolds $M_j \in {\mathcal M}^n_{H,D}$ are included in this stronger class as well.

It is also of interest to study sequences of manifolds which do not have uniform bounds on their Ricci or sectional curvatures.   

\subsection{Barbells}

In \cite{Fukaya:1987aa}, Fukaya considered sequences of {\bf barbells} which are pairs of spheres of radius $R_1$ and $R_2$ attached
by a cylinders of length $L$ and radius $r\le \min\{R_1,R_2\}$ that can be described by as warped products
over a sphere:
\be
(M_{R_1,R_2,L,r}, g_{R_1,R_2,L,r})=({\mathbb S}^2, ds^2+f^2(s)d\theta^2)
\ee
where
\begin{eqnarray*} \label{barbell-eq}
f(s) &=&R_1\sin((s+s_{R_1})/R_1) \textrm{ for } s\in [-s_{R_1},-L/2] \\
f(s) &=& r  \textrm{ for } s\in [-L/2, L/2] \\
f(s) &=& R_2\sin( (s_{R_2}-s)/R_2) \textrm{ for } s\in [L/2, s_{R_2}]]
\end{eqnarray*}
where $s_{R_i}$ are chosen so that we have a $C^0$ metric:
\be
R_1\sin((L/2+s_{R_1})/R_1)=r=R_2\sin( (s_{R_2}-L/2)/R_2).
\ee
These are $C^0$ Riemannian manifolds, but can easily be smoothed to be $C^\infty$ Riemannian manifolds
with metric tensors arbitrarily close in the $C^0$ sense to the given metric tensor in (\ref{barbell-eq}).   However the Ricci curvature of the
sequence approaching this decreases to negative infinity for any such sequence because these spaces
do not satisfy the Bishop-Gromov inequality.

\begin{figure}[h] %  figure placement: here, top, bottom, or page
   \centering
   \includegraphics[width=3in]{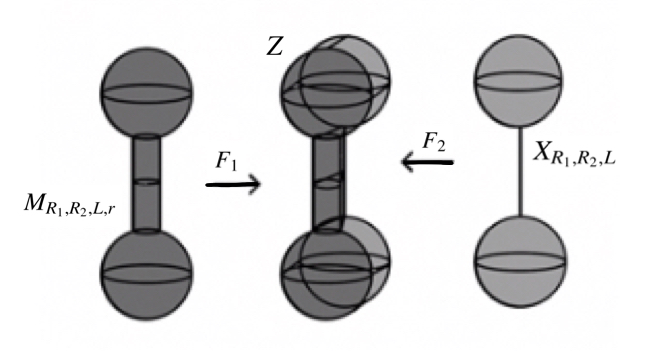} 
   \caption{Here we see distance preserving maps $F_1:M_{R_1,R_2,L,r}\to Z$
   and $F_2:X_{R_1,R_2,L}\to Z$ that can be used to estimate the Gromov-Hausdorff distance,
   $d_{GH}(M_{R_1,R_2,L,r}, X_{R_1,R_2,L})\approx \pi r$.
   }
   \label{fig:GH-barbell}
\end{figure}

Fukaya fixed $R_2=R_1$ and fixed $L>0$ and took $r\to 0$.  He explained that the sequence
of such barbells converges in the measured Gromov-Haudorff sense to a metric measure space which
is a pair of standard spheres joined by a line segment of length $L$:
\be
M_{R_1,R_2,L,r}
\mGHto X_{R_1,R_2,L}={\mathbb S}_{R_1}^1\disjointunion|_{\sim} [0,L]\disjointunion|_{\sim} {\mathbb S}_{R_2}^1 \textrm{ as } r\to 0
\ee   
where one pole ${\mathbb S}_{R_1}$ is identified to $0\in [0,L]$
and one pole in ${\mathbb S}_{R_2}$ is identified to $L\in [0,L]$.   See Figure~\ref{fig:GH-barbell}.

Fukaya noted measure on this limit space is the standard measure on each sphere and is $0$ on the interval.
He also studied sequences where the metric on the cylinder was increasingly bumpy so its
area would not converge to $0$ and the measure on the interval would be nonzero.
Furthermore he proved he had strong enough control on this sequence, that the eigenvalues 
of the limit space are found by taking all the eigenvalues of ${\mathbb S}_{R_1}^1$ and
all the eigenvalues of ${\mathbb S}_{R_2}^1$ and
all the eigenvalues of the interval $[0,L]$ and listing them in order with multiplicity.   Despite the lack of uniform
control on Ricci curvature, he stated that the eigenvalues of the barbells converge to these eigenvalues.  It would be interesting to investigate what the eigenfunctions and heat kernels look like on this sequence.  

Fukaya also studied that case where $R_1$ is fixed and $L\to 0$ and $r\to 0$ and $R_2>r\to 0$ so that
\be
M_{R_1,R_2,r,L}\GHto {\mathbb S}^2_{R_1}  \textrm{ as } r\to 0, \,\,L\to 0,\,\, R_2\to 0.
\ee
By taking
$r\to 0$ faster than the other two, he can use his previous result to estimate the eigenvalues.
In particular for $r$ sufficiently small, 
\be
M_{R_1,R_2,r,L}\approx X_{R_1,R_2,L}
\ee
and then 
\be
X_{R_1,R_2,L} \GHto {\mathbb S}^2_{R_1}  \textrm{ as } r\to 0, \,\,L\to 0,\,\, R_2\to 0.
\ee
He studied the eigenvalues in this setting and it would be interesting to study the eigenfunctions and heat kernels as well.

\subsection{Intrinsic Flat Convergence}\label{sect-SWIF}

The Sormani-Wenger intrinsic flat distance between a pair of oriented Riemannian manifolds,
$(M_1^n,g_1)$ and $(M_2^n, g_2)$, is essentially the filling volume between them.
It is defined by taking the infimum over all complete metric spaces $Z$ and all distance preserving maps,
$F_i: (M_i, d_{g_i})\to (Z,d_z)$  of the flat distance between their images:
\be
d_{SWIF}(M_1, M_2))=\inf d_F^Z(F_{1\#}[[M_1]], F_{2\#}[[M_2]]).
\ee
As it requires Geometric Measure Theory and Ambrosio-Kirchheim Theory to define
the flat distance between two submanifolds, we instead state a simple useful estimate:
\be
d_F^Z(F_{1\#}[[M_1]], F_{2\#}[[M_2]])\le \vol_n(A^n) +\vol_{n+1}(B^{n+1})
\ee
where $A^n$ is an $n$-dimensional submanifold of $Z$ and the filling manifold, $B^{n+1}$, is an $(n+1)$-dimensional submanifold of $Z$
such that the boundary of $B$ is $F_{1\#}[[M_1]]$, $F_{2\#}[[M_2]])$ and $A$ keeping track of the orientation:
\be \label{AB-filling}
\partial B=F_{1\#}[[M_1]]-F_{2\#}[[M_2]])-A.
\ee
More
generally, $A$ and $B$ are integral currents satisfying (\ref{AB-filling} and their volumes are replaced by their Ambrosio-Kirchheim masses.   In fact the intrinsic flat distance is defined between pairs of {\bf integral current spaces} which are rectifiable metric
spaces with rectifiable boundaries that have a weighted orientation called an integral current structure.
The intrinsic flat limits of sequences of oriented Riemannian manifolds are integral current spaces included possibly the $0$ space when the sequence collapses \cite{Sormani-Wenger:2011}.

This notion of convergence was defined to study the noncollapsed limits of sequences of Riemannian manifolds
with nonnegative scalar curvature.   It is designed so that thin wells disappear in the limit.

\begin{figure}[h] %  figure placement: here, top, bottom, or page
   \centering
   \includegraphics[width=5in]{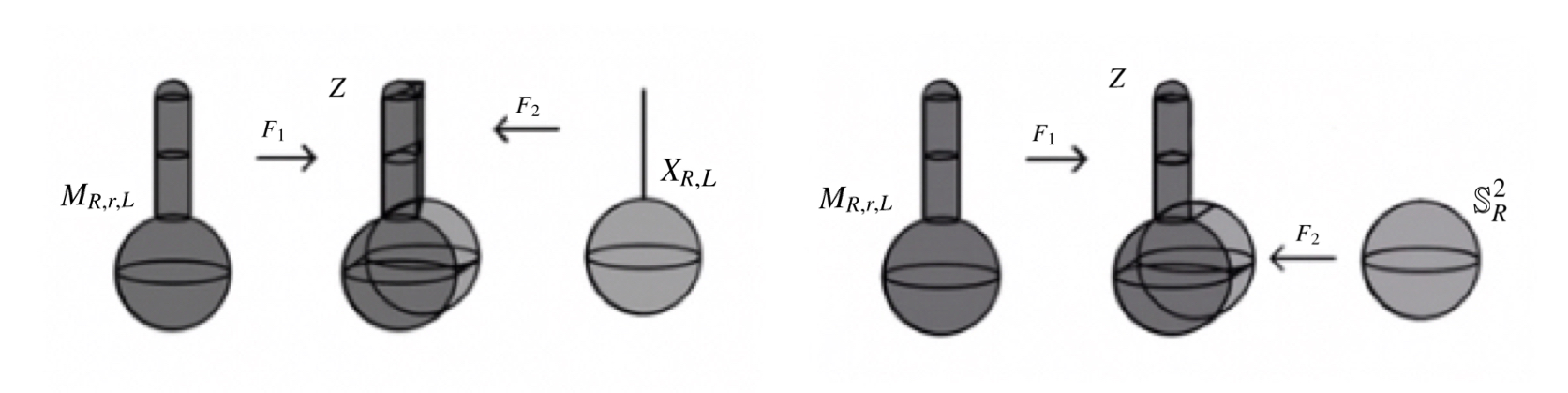} 
   \caption{On the left we see distance preserving maps $F_1:M_{R,L,r}\to Z$
   and $F_2:X_{R,L}\to Z$ that can be used to estimate the Gromov-Hausdorff distance,
   $d_{GH}(M_{R,L,r}, X_{R,L})\approx \pi r$.    On the right we see the same maps can be used
   to show the intrinsic flat distance between  $M_{R,L,r}$ and ${\mathbb S}^2_R$ is small.  }
   \label{fig:GH-SWIF-lollipop}
\end{figure}

Sormani-Wenger studied spheres with wells that our student research team liked to call {\bf lollipops}.  These consist of
a single sphere of radius $R$ attached to a capped cylinder of length $L$ and radius $r\le R$ in \cite{Sormani-Wenger:2011}.   These manifolds
\be
(M_{R,r,L}, g_{R,r,L}) =({\mathbb S}^2, ds^2+f^2(s)d\theta^2)
\ee
are defined using the following $C^0$ warping functions:
\begin{eqnarray*} \label{lollipop-eq}
f(s) &=&R\sin((s+s_{R})/R) \textrm{ for } s\in [-s_{R},-L/2] \\
f(s) &=& r  \textrm{ for } s\in [-L/2, L/2] \\
f(s) &=& r\sin( (s_{r}-s)/R_2) \textrm{ for } s\in [L/2, s_{r}]]
\end{eqnarray*}
as in (\ref{barbell-eq}) with $R_1=R$ and $R_2=r$.
If we fix $R=r>0$ and take $L\to 0$ these converge in the $C^0$ sense to a rescaled sphere of radius $R$:
\be
M_{R,r,L}\to{\mathbb S}^2_R \textrm{ as }  L\to 0.
\ee
If we fix $R, L>0$ and take $r_j\to 0$, these have no $C^0$ limit, but converge in the $mGH$ sense to
a sphere of radius $R$ attached to an interval of length $L$:
\be
M_{R,r_j,L}\mGHto X_{R,L}={\mathbb S}^1\disjointunion [0,L] \textrm{ as }  r_j\to 0.
\ee
See Figure~\ref{fig:GH-SWIF-lollipop}.
The measure on this limit space is the standard measure on the sphere and is $0$ on the interval.
Sormani-Wenger showed intrinsic flat limit is just the sphere by explicitly constructing distance preserving maps into a sequence of complete spaces $Z_j$ and constructing fillings $B_j$ of one dimension higher between the sphere and the lollipop that has very small volume.    
Portegies studied the eigenvalues
of this example in \cite{Portegies-evalues} and it would be interesting to know how the eigenfunctions and heat kernels behave.

If the warping functions on the sequence are increasingly bumpy on the cylindrical part instead of
being constant, Sormani-Wenger showed the measured Gromov-Hausdorff limit has a nonzero measure 
on the interval but the intrinsic flat limit is still only the sphere \cite{Sormani-Wenger:2011}.  

Ilmanen created a sequence of spheres with increasingly many increasingly thin capped cylinders with no Gromov-Hausdorff limit because such a sequence cannot be embedded with distance preserving maps into a compact metric space, $Z$.   Sormani-Wenger showed that if the volume of the
sequence is kept bounded then the filling volumes between these manifolds and the sphere converges to $0$, so
that the intrinsic flat limit of Ilmanen's Example is a standard sphere.   The Intrinsic Flat convergence was designed to handle
sequences like these that have no Gromov-Hausdorff limit.   Wenger's Compactness Theorem states that any
sequence of manifolds with $\vol(M_j)\le V$ and $\diam(M_j)\le D$ has a subsequence which converges in the
intrinsic flat sense to an integral current space.

In \cite{Portegies-evalues}, Portegies introduced  {\bf volume preserving intrinsic flat convergence},
$M_j\to M_\infty$ iff 
\be
 M_j\to M_\infty \textrm{ in the intrinsic flat sense and }\vol(M_j) \to \vol(M_\infty).
\ee   
Applying the work of Fukaya,
he proved that under this convergence 
\be
\limsup_{j\to \infty} \lambda_k(M_j) \le \lambda_k(M_\infty).
\ee
It would be interesting to study what happens to the heat kernels in this setting as well. 

\subsection{VADB Convergence}
Most recently the notion of VADB Convergence has been defined by Allen-Perales-Sormani
in \cite{APS-JDG} and they proved VADB convergence implies Intrinsic Flat Convergence.
Ilmanen's Example converges to a sphere under VADB convergence.   

A sequence of diffeomorphic manifolds $(M,g_j)$ converge to $(M, g_0)$ in the {\bf VADB} sense
if $\vol_{g_j}(M)\to \vol_{g_0}(M)$, $\diam_{g_j}(M)\le D$, and $g_j \ge (1-1/j)g_0$.   

This paper is building on prior work of All-Sormani concerning the 
$L^p$ convergence of the
metric tensors which is not enough to achieve uniform, Gromov-Hausdorff, or intrinsic flat convergence
unless there are also uniform bounds on
the metric tensors from above and below as shown by Allen-Sormani in \cite{Allen-Sormani:2019} and \cite{Allen-Sormani:2020}.  

A number of conjectures concerning VADB convergence as well as a survey of results applying VADB
convergence appears in a survey by the second author \cite{Sormani-IAS}.   It would be interesting to study how the eigenvalues, eigenfunctions, and heat kernels behave under VADB convergence.

%%%%%%%%%%%%%%
\section{\bf Embedding Riemannian Manifolds via the Heat Kernel}\label{sect-embed-map}

In 1994 B\'erard, Besson, and Gallot introduced the idea of embedding a Riemannian manifold into $\ell^2$ using the heat kernel \cite{BBG:1994}.   Recall that $\ell^2$ is the vector space of square summable sequences,
\be
\ell^2=\left\{ \, \{v_j\}_{j=0}^\infty=\{v_0,v_1,v_2,v_3,...\}\, :\, \sum_{i=0}^\infty v_i^2 <\infty \right\},
\ee
with the $\ell^2$-inner product and $\ell^2$-norm,
\be
<\{v_j\}_{j=0}^\infty, \{\tilde{v}_j\}_{j=0}^\infty>_{\ell^2}\,=\,\sum_{j=0}^\infty v_j\cdot \tilde{v}_j
\quad \textrm{ and } \quad ||\{v_j\}_{j=0}^\infty||_{\ell^2}^2 =\,\sum_{j=0}^\infty v_j^2.
\ee
Given an $L^2(M)$ orthonormal basis of eigenfunctions 
\be
\{\phi_0, \phi_1, \phi_2, ....\} \textrm{ such that } \Delta \phi_j=-\lambda_j\phi_j
\textrm{ with }0=\lambda_0<\lambda_1\le \lambda_2\le ...,
\ee 
B\'erard-Besson-Gallot define the family of maps
$
\Psi_t: M \to \ell^2
$
for $t>0$ as
\be\label{BBG-map}
\Psi_t(x)=\left\{e^{-\lambda_j t/2 }\phi_j(x)\right\}_{j=1}^\infty.
\ee
They use the fact that
\begin{eqnarray*}
||\Psi_{t}(x)-\Psi_t(y)||^2_{\ell^2}&=&\sum_{j=1}^\infty e^{-\lambda_j t } (\phi_j(x)-\phi_j(y))^2
=\sum_{j=0}^\infty e^{-\lambda_j t } (\phi_j(x)-\phi_j(y))^2\\
&=&\sum_{j=0}^\infty e^{-\lambda_j t } (\phi_j(x))^2
+\sum_{j=0}^\infty e^{-\lambda_j t } (\phi_j(y))^2
-2 \sum_{j=0}^\infty e^{-\lambda_j t } \phi_j(x)\phi_j(y)\\
&=&h_{t}(x,x)+h_{t}(y,y)-2 h_{t}(x,y) 
\end{eqnarray*}
%and furthermore that
%\be
%|\Psi_{t_i}(x_i)-\Psi_t(x)|=h_{t_i}(x_i,x_i)+h_t(x,x)-2 h_{(t_i+t)/2}(x_i,x)
%\ee
to prove the map is continuous and thus the image is compact.   

To prove this is an embedding, B\'erard-Besson-Gallot
 use the fact that any $L^2(M)$ basis of orthonormal eigenfunctions {\em separates points}:
\be
\forall p, q\in M \,\, \exists \phi_j \textrm{ such that } \phi_j(p)\neq \phi_j(q).
\ee
This is also used to prove that the $L^2(M)$ basis of orthonormal eigenfunctions spans $L^2(M)$. 

Note that the the eigenvalues, $0=\lambda_0<\lambda_1\le \lambda_2\le...$
and the eigenspaces of these eigenvalues are uniquely determined by $M$ and 
these eigenspaces are $L^2(M)$-orthogonal to
one another.  There is a choice being made when
selecting the $L^2(M)$-orthonormal basis for each eigenspace, and thus the
$L^2(M)$-orthonormal eigenfunctions depend on this choice.  For this reason
B\'erard-Besson-Gallot label their eigenfunctions with an $a$.

\subsection{Rescaled Heat Kernel Embeddings}\label{sect-re-map}

B\'erard, Besson, and Gallot \cite{BBG:1994} defined rescaled heat kernel embeddings
and proved the following theorem:

\begin{thm}
Given an $n$-dimensional closed Riemannian manifold,  let the map
$
\psi_t: M^n \to \ell^2
$
for $t>0$ be defined by
\be
\psi_t(x)=
 \sqrt{2}(4\pi)^{n/4} t^{(n+2)/4} \Psi_t(x).
\ee
Then the pulled back metric tensor satisfies
\be \label{rBBG-map}
(\psi_t)^*g_{\ell^2} = g+ \tfrac{t}{3}(\tfrac{1}{2}\Scal_g \cdot g - \Ricci_g)  + \textrm{O}(t^2) \textrm{ when }t\to 0^+
\ee
where the pulled back metric tensor is defined on a vector $V_x\in TM_x$ by
\be
(\psi_t)^*g_{\ell^2} (V_x,V_x)= || d\psi_{t}V_x||^2_{\ell^2}=
\left(\sqrt{2}(4\pi)^{n/2} t^{(n+2)/4}\right)^2\sum_{j=1}^\infty e^{-\lambda_jt} |d\phi_j V_x|^2.
\ee
\end{thm}

B\'erard-Besson-Gallot  proved (\ref{rBBG-map}) by
observing that 
\be
\sum_{j=1}^\infty e^{-\lambda_jt} |d\phi_j V_x|^2 
= \partial_{y_1}\partial_{y_2} h_t(y_1, y_2)|_{y_1=y_2=x}.
\ee
They were able to prove their theorem taking these partial derivatives on the heat kernel using the 
Minakshisundaram-Pleijel expansion.

Minakshisundaram and Pleijel studied the heat kernel on an manifold $M$ 
within a normal neighborhood about a point $x$ using the method of analytic continuation applied to the time variable \cite{Minak-Pleijel:1949}.  They required that the metric tensor $g$ of $M$ is an analytic function of the normal coordinates.
They prove the heat kernel is the $C^{\omega}$ limit
\be \label{Minak-Pleijel-1}
h_t(x,y)=\lim_{P\to \infty} {\tilde h}_t^P(x,y)\\
\ee
where
\be\label{Minak-Pleijel-2}
{\tilde h}_t^P(x,y)=  \tfrac{1}{(4\pi t)^{n/2} } \, e^{-d_M(x,y)^2/(4t)}   \sum_{i=0}^P U_i(x,y) \, t^i
\ee
where $U_i(x,y)$ satisfy an iterative set of partial differential equations (each depending on the value of
the previous $U_i(x,y)$ and the metric tensor $g$). They proved the $U_i(x,y)$ are uniquely
determined by the requirement that $U_0(x,x)=1$ and $U_i(x,x)$ is finite.     One advantage of
Minakshisundaram-Pleijel's expansion is that a $C^{\omega}$ limit allows us to differentiate as many times as we wish and the derivatives will converge to the derivatives of the limit.   The
proof of convergence in \cite{Minak-Pleijel:1949} involves the decay of the eigenvalues.  The disadvantage of this formula is that it requires that $x,y$ lie in a normal neighborhood.

Abdalla later applied Minakshisundaram and Pleijel to study $C^{\omega}$ parametrized families of Riemannian manifolds embedded via heat kernel maps in \cite{Abdalla:2012}.

\subsection{Truncated Heat Kernel Embeddings}\label{sect-trunc-map}

It is much simpler to consider the truncated heat kernel embeddings into Euclidean space
$
\psi^N_t: M^n \to {\mathbb E}^N
$
for $t>0$ such that
\be\label{BBG-trunc}
\psi^N_t(x)= \sqrt{2}(4\pi)^{n/4} t^{(n+2)/4} \{e^{-\lambda_j t/2 }\phi_j(x)\}_{j=1}^N \in {\mathbb E}^N.
\ee
Note that we easily obtain continuity for the truncated heat kernel maps and they are one-to-one as long as 
$\{\phi_1,...\phi_N\}$ separate points.   
 
We can then rewrite B\'erard-Besson-Gallot's Theorem stated in (\ref{rBBG-map}) as 
\be\label{metric-lim}
\lim_{t\to 0^+}\lim_{N\to \infty} (\psi_t^N)^*g_{{\mathbb E}^N} (V_x,V_x)= g(V_x,V_x),
\ee
and 
\be
\lim_{t\to 0^+}(3/t)  \left(\lim_{N\to \infty}(\psi_t^N)^*g_{{\mathbb E}^N} (V_x,V_x) - g(V_x,V_x)\right)=
\tfrac{1}{2}\Scal_g \cdot g - \Ricci_g
\ee
where
\be \label{this-sum}
(\psi_t^N)^*g_{{\mathbb E}^N} (V_x,V_x)
=\left(\sqrt{2}(4\pi)^{n/4} t^{(n+2)/4}\right)^2\sum_{j=1}^N e^{-\lambda_jt} |d\phi_j V_x|^2.
\ee
Note that these limits must be taken in the correct order.   If we try to take the limit as $t\to 0^+$ first in (\ref{metric-lim}),
then
\be
\lim_{t\to 0^+}(\psi_t^N)^*g_{{\mathbb E}^N} (V_x,V_x)
=\lim_{t\to 0^+}\left(\sqrt{2}(4\pi)^{n/4} t^{(n+2)/4}\right)^2\sum_{j=1}^N e^{-\lambda_jt} |d\phi_j V_x|^2
= 0 
\ee
because $\lim_{t\to 0+} t^{(n+2)/4} e^{-\lambda t} =0$.   

We can define a BBG time $T_\epsilon$ such that
\be\label{BBG-time}
|(\psi_t)^*g_{\ell^2} (V_x,V_x)-g(V_x,V_x)|<\epsilon g(V_x,V_x) \quad \forall t\in (0,T_\epsilon)
\ee
and a BBG number $N_{t,\varepsilon}$
such that
\be \label{BBG-N}
|(\psi_t^N)^*g_{{\mathbb E}^N} (V_x,V_x)- (\psi_t)^*g_{\ell^2}(V_x,V_x)| <\varepsilon g(V_x,V_x) \quad \forall N\ge N_{t,\varepsilon}
\ee
If we choose $t<T_\epsilon$ first and then $N>N_{t,\varepsilon}$ we have
\be
|(\psi_t^N)^*g_{{\mathbb E}^N} (V_x,V_x)-g(V_x,V_x)|<(\epsilon+\varepsilon) g(V_x,V_x).
\ee
These estimates are explored further in the work of Portegies \cite{Portegies:2016} which we will discuss below
in Section~\ref{sect-Portegies}.

Note that there is some lack of uniqueness when defining the truncated embeddings.  If $N$ has
$\lambda_N =\lambda_{N+1}$, then $\psi_t^N$ includes $\phi_N$ in the truncated embedding but doesn't include $\phi_{N+1}$
although both are in the same eigenspace and could have been listed in any order when first finding the 
$L^2(M)$ basis of eigenfunctions.   If $N$ has
$\lambda_N <\lambda_{N+1}$, then the truncated embedding is unique up to an orthogonal transformation
of ${\mathbb E}^N$, because we would be including a complete $L^2(M)$ basis for each eigenspace with $j\le N$
and choosing a different basis for the $m_j$ dimensional eigenspace of $\lambda_j$
only transforms the basis by an $O(m_j)$ action and thus transforms the image in ${\mathbb E}^N$ by
the same action.

\subsection{Embedding Circles with a Truncated Heat Kernel }\label{sect-circle-map}

Let us consider 
$\psi_t: {\mathbb S}^1 \to \ell^2$
for $t>0$ using the eigenfunctions (which come in pairs) that we found in the Section~\ref{sect-circle-evalue} .   The truncated embedding is
$$
\psi_t^{2N}(x)= \sqrt{2}(4\pi)^{n/4} t^{(n+2)/4} \left\{e^{-j^2 t/2 }\pi^{-\tfrac{1}{2}\,\,}\sin(j\theta),e^{-j^2 t/2} \pi^{-\tfrac{1}{2}\,\,}\cos(j\theta)) \right\}_{j=1}^N \subset {\mathbb E}^{2N}.
$$

\begin{figure}[h] %  figure placement: here, top, bottom, or page
   \centering
   \includegraphics[width=5in]{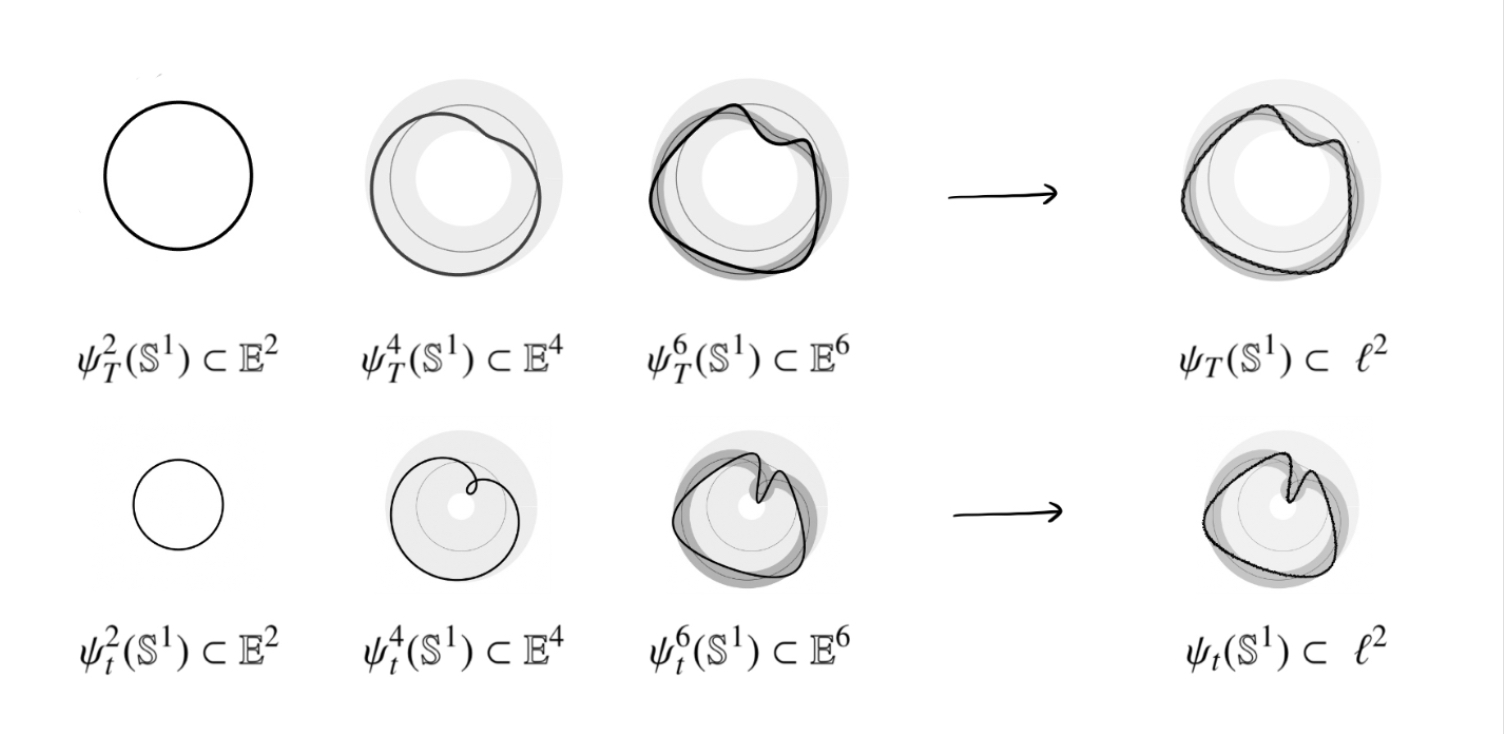} 
   \caption{Here we look at the images of the heat kernel embeddings of the circle
   for $T>t$.   Image credit: Pen Chang (Student, Stony Brook).
   }
   \label{fig-circle-embeddings}
\end{figure}

In Figure~\ref{fig-circle-embeddings}, we have drawn some of the images of the truncated heat kernel embeddings with $T>t$.   Note that already with $N=2$ we have an embedding, but the circle is too short.  With $N=4$ we wind around $\psi_{T}^2({\mathbb{S}}^1)$ twice, and with $N=6$ we wind around $\psi_{T}^4({\mathbb{S}}^1)$ thrice.   Taking limit as $N\to \infty$ for fixed time $T$, we moving to the right across the first row
of Figure~\ref{fig-circle-embeddings}.   The limit, $\Psi_T({\mathbb S}^1)$ is a circle winding infinitely many times around itself.   The same happens in the second row.

Now consider the columns in Figure~\ref{fig-circle-embeddings}.  If we take $t$ small, for fixed $N$, we see that the radius of winding is closer to the original radius, because
$e^{-j^2t}$ is close to $1$ for both $j=1$ and $j=2$.  However, the whole image is scaled down by $\sqrt{2}(4\pi)^{n/2} t^{(n+2)/4}$.
So for fixed $N$ taking $t\to 0$ the images contract to a point. 

To find the lengths of these curves, we use the homogeneity of ${\mathbb{S}}^1$ to see that
\be
L\left(\psi_t^{2N}({\mathbb S}^1)\right)= \int_0^{2\pi} (\psi_t^{{2N}})^*g_{{\mathbb E}^{2N}} (\partial_ \theta,\partial_ \theta)^{\tfrac{1}{2}\,} \, d\theta=2\pi (\psi_t^{{2N}})^*g_{{\mathbb E}^{2N}} (\partial_ \theta,\partial_ \theta)^{\tfrac{1}{2}\,}.
\ee
We can evaluate this
using (\ref{this-sum}) at any $\theta$.  So lets take $\theta=0$
where our eigenfunctions, $\phi_{2j}(0)=0$ and $\phi_{2j-1}(0)=\pi^{-\tfrac{1}{2}\,\,}$, to see that
 \be
(\psi_t^{{2N}})^*g_{{\mathbb E}^{2}N} (\partial_\theta,\partial_\theta)=
2(4\pi)^{\tfrac{1}{2}\,} t^{3/2}\sum_{j=1}^{2N} e^{-\lambda_jt} \lambda_j |\phi_j(0)|^2=
4\pi^{-\tfrac{1}{2}\,\,} t^{3/2}\sum_{j=1}^{N} e^{-j^2t} j^2. 
\ee
So the lengths of the images of ${\mathbb S}^1$ are 
\be \label{length}
L\left(\psi_t^{2N}({\mathbb S}^1)\right)=2\pi \left(  4\pi^{-\tfrac{1}{2}\,\,}  t^{3/2}\sum_{j=1}^{N} e^{-j^2t} j^2    \right)^{\tfrac{1}{2}\,}
\ee
and the lengths of their limits as $N\to \infty$ increase to
\be
L\left(\psi_t({\mathbb S}^1)\right)=2\pi \left(4\pi^{-\tfrac{1}{2}\,\,}  t^{3/2}\sum_{j=1}^{\infty} e^{-j^2t} j^2 \right)^{\tfrac{1}{2}\,}
\ee
and by (\ref{metric-lim}) we have
\be\label{verified}
\lim_{t\to 0^+} L\left(\psi_t({\mathbb S}^1)\right)=2\pi.
\ee
However if we take $t\to 0$, before $N\to\infty$, we get
$$
\lim_{t\to 0^+} L\left(\psi_t^{2N}({\mathbb S}^1)\right)=\lim_{t\to 0^+} 2\pi \left(  4\pi^{-\tfrac{1}{2}\,\,}  t^{3/2}\sum_{j=1}^{N} e^{-j^2t} j^2    \right)^{\tfrac{1}{2}\,}
=0(1^2+2^2+\cdots +N^2)=0.
$$

Note that (\ref{verified}) can be independently verified using a Theta function as seen 
in Edwards' text \cite{Edwards-text}.
Consider 
\be
\Theta(x) = \sum_{n=1}^{\infty} e^{-\pi n^2 x}.
\ee
Then $\Theta(x)$ satisfies the equation
\be
\frac{1+2\Theta(x)}{1+2\Theta(1/x)} = \frac{1}{\sqrt{x}}.
\ee
Differentiating this equation, we obtain the following:
\be
-4 x^{3/2} \Theta'(x) = 1 + 2 \Theta(1/x) + \frac{4}{x} \Theta'(1/x)
\ee
That is, 
\be
4 \pi x^{3/2} \sum_{n=1}^{\infty}  n^2 e^{-\pi n^2 x} = 1 + 2 \Theta(1/x) + \frac{4}{x} \Theta'(1/x)
\ee
Setting $t = \pi x$, we get 
\be
4 \pi^{-1/2} t^{3/2} \sum_{n=1}^{\infty}  n^2 e^{- n^2 t} = 1 + 2 \Theta(\pi/t) + \frac{4\pi}{t} \Theta'(\pi/t)
\ee
Note that the right hand side approaches $1$ as $t \to 0$ due to exponential decay of the Theta function. Hence, we have 
\be
\lim_{t\to 0^+} L\left(\psi_t({\mathbb S}^1)\right)=2\pi.
\ee

\subsection{Embedding Spheres with a Truncated Heat Kernel } \label{sect-sphere-map}

For a round sphere, ${\mathbb S}^2$ of dimension $2$, the first $3$ nonconstant eigenfunctions all share the same eigenvalue (as we saw in \ref{sect-sphere-evalue}).   Due to the symmetry, the truncated heat kernel $\psi^3_{t}({\mathbb S}^2)$ is an embedding whose image is a round sphere in
${\mathbb E}^{3}$ but its radius is smaller than that of the original sphere and decreases to $0$ as $t\to 0$.   In fact, for a round sphere, ${\mathbb S}^n$ of dimension $n$, the first $n+1$ nonconstant eigenfunctions all share the same eigenvalue.  Thus for $N=n+1$, the truncated heat kernel,
$\psi^{N}_{t}({\mathbb S}^n)$, is an embedding whose image is a round sphere in ${\mathbb E}^{n+1}$.

If a manifold $(M^n,g)$ is $C^2$ close to a round sphere ${\mathbb S}^n$, then its eigenvalues and eigenvectors will be close to that of the sphere.   If they are sufficiently close then $\psi^N_t: M^n \to {\mathbb E}^{N}$ will be an embedding for $N=n+1$ as well.   

In fact, if $(M_j^n,g_j)$ have
$\Ricci \ge (n-1)H$ and $M^n_j \GHto {\mathbb S}^n$ then by Cheeger-Colding \cite{ChCo-Part3}, their eigenfunctions and
eigenvalues converge.  So perhaps, for $j$ sufficiently large, $\psi^N_{j\,t}: M_j \to {\mathbb E}^N$ is an embedding for $N=n+1$.

\subsection{Embedding Tori with a Truncated Heat Kernel Map} \label{sect-tori-map}

For a flat torus of the form ${\mathbb S}^1 \times {\mathbb S}^1$, we saw in Section~\ref{sect-tori-evalue} that
the first four nonconstant eigenfunctions are
\be
A_{1,0}\sin(\theta), \,\,\, A_{1,0}\cos(\theta),\,\,\, A_{0,1}\cos(\varphi), \textrm{ and } A_{0,1}\sin(\varphi)
\ee
where $A_{1,0}=A_{0,1}=1/(\pi \sqrt{2})$ with eigenvalue $\lambda_{1,0}=\lambda_{0,1}=1$
because $A=B=1$.
So the truncated truncated heat kernel map, $\psi_t^{4}$,
maps the torus onto 
\be
{\mathbb S}_R^1 \times {\mathbb S}_R^1\subset {\mathbb E}^2 \times {\mathbb E}^2
\ee 
where $R=1/(\pi \sqrt{2})$.
Thus the pullback of the Euclidean metric tensor, $\psi_t^{4*}g_{{\mathbb E}^4}$,
is just a rescaling of the flat 
torus metric tensor, $g_{{\mathbb S}^1 \times {\mathbb S}^1}$.   

For a flat torus of the form ${\mathbb S}^1 \times {\mathbb S}_{1/10}^1$ we saw in Section~\ref{sect-thin-evalue} that the first twenty eigenfunctions
only depend upon on the first parameter and so $\psi_t^{4}$ even $\psi_t^{19}$ maps into a parametrized circle lying in ${\mathbb E}^{19}$
exactly as in the image of a single ${\mathbb S}^1$ although it will be scaled differently.    The pull back of the metric tensor, $\psi_t^{19*}g_{{\mathbb E}^{19}}$ will not be a positive definite metric tensor on the torus.  On a thin torus, we need to take $N$ very large before the pull back becomes positive definite.

\subsection{Embedding Tori with Selected Eigenfunctions } \label{sect-tori-smap}

Here we embed the thin torus, ${\mathbb S}^1 \times {\mathbb S}_{1/10}^1$,
with its metric tensor of the
form,
\be
g_{{\mathbb S}^1 \times {\mathbb S}_{1/10}^1}= d\theta^2 + (1/10)^2 d\varphi^2,
\ee
into ${\mathbb E}^4$ using only four eigenfunctions.   To do this, we select the four eigenfunctions corresponding to 
\be
\lambda_{0,1}=\lambda_{10,0}=100
\ee
and define the map, 
\be
\psi^{\lambda=100}_t(\theta,\phi)=\left(\,\tfrac{\sqrt{5}}{\pi} \sin(10\,\theta),
\,\tfrac{\sqrt{5}}{\pi} \cos(10\,\theta),
\,\tfrac{\sqrt{5}}{\pi} \sin(\varphi),\,\tfrac{\sqrt{5}}{\pi} \cos(\varphi)\right).
\ee
We claim that
\begin{eqnarray}
\left(\psi_t^{\lambda=100}\right)^*g_{{\mathbb E}^4}&=&\left(\psi_t^{\lambda=100}\right)^*(dx_1^2+dx_2^2+dx_3^2+dx_4^2)\\
&=& \left(10 \,{\sqrt{5}}/{\pi} \right)^2(\cos^2(10\,\theta)+\sin^2(10\,\theta)) \, d\theta^2\\
&&\qquad + \,\, \left({\sqrt{5}}/{\pi} \right)^2(\cos^2(\varphi)+\sin^2(\varphi)) \, d\varphi^2.
\end{eqnarray}
To see this, observe that
\begin{eqnarray*}
\left(\left(\psi_t^{\lambda=100}\right)^*(dx_i)\right)(\partial_{\theta})
&=&dx_1\left((\psi^{\lambda=100}_t)_*\partial_\theta\right)
=\left((\psi_t^{\lambda=100})_*\partial_\theta\right) x_i\\
&=&\left((\psi_t^{\lambda=100})_*\partial_\theta\right) x_i
= \partial_\theta \left(x_i\circ \Psi^{\lambda=100}_t\right).
\end{eqnarray*}
So we have
\begin{eqnarray}
\left(\left(\psi_t^{\lambda=100}\right)^*(dx_1)\right)(\partial_{\theta})&=&
\partial_\theta \left(\tfrac{\sqrt{5}}{\pi} \sin(10\,\theta)\right)=
10 \,\tfrac{\sqrt{5}}{\pi} \cos(10\,\theta)\\
 \left((\psi_t^{\lambda=100})^*(dx_2))(\partial_{\theta}\right) 
 &=& \partial_\theta \left(\tfrac{\sqrt{5}}{\pi} \cos(10\,\theta)\right)=
-10 \,\tfrac{\sqrt{5}}{\pi} \sin(10\,\theta)\\
 \left((\psi_t^{\lambda=100}\right)^*(dx_3))(\partial_{\phi})
 &=& \partial_\varphi \left(\tfrac{\sqrt{5}}{\pi} \sin(\varphi)\right)=
\tfrac{\sqrt{5}}{\pi} \cos(\varphi)\\
 \left((\psi_t^{\lambda=100}\right)^*(dx_4))(\partial_{\phi})
 &=& \partial_\varphi \left(\tfrac{\sqrt{5}}{\pi} \cos(\varphi)\right)=
- \tfrac{\sqrt{5}}{\pi} \sin(\varphi)
\end{eqnarray}
and the other pullbacks of $dx_i$ acting on $\partial_{\theta}$ and $\partial_{\phi}$ are zero.
Thus
\be
(\psi^{\lambda=100})^*g_{{\mathbb E}^4}= \left(10 \,\tfrac{\sqrt{5}}{\pi} \right)^2 d\theta^2+
\left(\tfrac{\sqrt{5}}{\pi} \right)^2 d\varphi^2
\ee
is a rescaling of 
\be
g_{{\mathbb S}^1 \times {\mathbb S}_{1/10}^1}= d\theta^2 + (1/10)^2 d\varphi^2
\ee

This also works to embed ${\mathbb S}^1 \times {\mathbb S}_{{1}/{2}\,}^1$
into ${\mathbb E}^4$.   We can take the four eigenfunctions corresponding to $\lambda_{0,1}=\lambda_{4,0}=4$
to define a {\bf {\em selective embedding}}:
\be
\psi^{\lambda=4}_t(\theta,\phi)=e^{-4t/2}\left(A_{2,0}\sin(2\theta),A_{2,0}\cos(2\theta),A_{0,1}\sin(\varphi),A_{0,1}\cos(\varphi)\right).
\ee
Note that a {\bf {\em selective embedding}} as described above
is a well chosen projection of the heat kernel embedding in carefully
selected directions.   This is our own notion introduced vaguely here
and is not mentioned in the literature elsewhere.

This will not work for arbitrary ${\mathbb S}^1_A \times {\mathbb S}_B^1$ with $A>B$ as $\lambda=(1/B)^2$
might only have multiplicity $2$ with a pair of eigenfunctions that only depend on $\varphi$.   In that case we might take a pair of eigenvectors depending on $\theta$ with
an eigenvalue {\em near} $(1/B)^2$ to define an  {\bf {\em approximately selective embedding}}
\be
\psi_t^{\lambda \textrm{ near } (1/B)^2}: {\mathbb S}^1_A \times {\mathbb S}_B^1\to {\mathbb E}^4.
\ee
Although this would not be a rescaling, it would be biLipschitz close to one.
  
It would be easier to find
a pair of eigenvectors depending on $\theta$ with
an eigenvalue nearer to $(k/B)^2$ for large $k$.  We might then use these in conjunction with a pair of eigenvectors depending on $\varphi$ with $\lambda=(k/B)^2$ to define an {\bf {\em approximately selective embedding}}
\be
\psi_t^{\lambda \textrm{ near } (k/B)^2}: {\mathbb S}^1_A \times {\mathbb S}_B^1\to {\mathbb E}^4.
\ee   
We might even conjecture that
that
with appropriate rescaling by a well chosen sequence, $\zeta_k\in (0,\infty)$, the
pullbacks of these {\bf {\em approximately selective embedding}} converge to the standard flat metric
on the torus:
\be
\lim_{k\to \infty} \zeta_k (\psi^{\lambda \textrm{ near } (k/B)^2})^*g_{{\mathbb E}^4} =
g_{{\mathbb S}^1_A \times {\mathbb S}_B^1}\,.
\ee
This would be straightforward to check for tori.
Could something like this work for more general classes of manifolds?   See Conjecture~\ref{conj-selective}
at the end of this paper.

\subsection{Estimating Lengths vs Extrinsic Distances}\label{sect-length-vs-distance}

Note that B\'erard-Besson-Gallot heat kernel embeddings can be used to approximate the
lengths of curves:
\be
L_g(C)=\int_0^1 g(C'(s),C'(s))^{1/2\,} \, ds.
\ee
If we choose $t<T_\epsilon$ as in (\ref{BBG-time})
then
\be
|L_{(\psi_t)^*g_{\ell^2} }(C)-L_g(C)|< \epsilon L_g(C).
\ee
So
\be
(1-\epsilon)\, L_g(C) <L_{(\psi_t)^*g_{\ell^2} }(C) <(1+\epsilon) \,L_g(C) 
\ee
The Riemannian distance between points $x,y\in M$ is found by taking
\be
d_M(x,y)=\inf\{L(C)\,:\, C(0)=x,\,C(1)=y, \, C:[0,1]\to M\}.
\ee
So
\be
(1-\epsilon)\, d_M(x,y) <d_{\psi_t(M)}(\psi_t(x),\psi_t(y)) <(1+\epsilon) \,d_M(x,y)
\ee
So $\psi_t$ is $(1+\epsilon)$-biLipschitz with respect to the intrinsic metric of the image $\psi_t(M)$
but not with respect to the metric $d_{\ell^2}$.   

To make this more clear, lets use the truncated heat kernel embeddings.  Taking
$N>N_{t,\varepsilon}$ so that we have (\ref{BBG-N}), we can repeat the above to see that
\be
(1-\epsilon-\varepsilon) L_g(C) <L_{(\psi^N_t)^*g_{{\mathbb E}^N} }(C) <(1+\epsilon+\varepsilon) L_g(C) 
\ee
and
\be
(1-\epsilon-\varepsilon) \,d_M(x,y) <d_{\psi^N_t(M)}(\psi_t(x),\psi_t(y)) <(1+\epsilon+\varepsilon) \,d_M(x,y)
\ee
So $\psi_t$ is $(1+\epsilon+\varepsilon)$-biLipschitz with respect to the intrinsic distance within the image 
$\psi^N_t(M)$
but not with respect to the extrinsic Euclidean distance.

In fact,
the extrinsic Euclidean distance between the points in the image are converging to $0$.  This can be seen 
as follows:
\be
d_{{\mathbb E}^N}( \psi_t(x),\psi_t(y))\le  d_{{\mathbb E}^N}( \psi_t(x),0)+d_{{\mathbb E}^N}( 0,\psi_t(y))
\ee
and
\begin{eqnarray}
d_{{\mathbb E}^N}( \psi_t(x),0)^2 &=&
\left(\sqrt{2}(4\pi)^{n/4} t^{(n+2)/4}\right)^2 \sum_{j=1}^N e^{-\lambda_j t}(\phi_j(x))^2\\
&=& 
\left(\sqrt{2}(4\pi)^{n/4} t^{(n+2)/4}\right)^2 h_t^N(x,x)
\end{eqnarray}
Taking $N\to \infty$ first we get
\be
d_{\ell^2}( \psi_t(x),\psi_t(x))^2=\left(\sqrt{2}(4\pi)^{n/4} t^{(n+2)/4}\right)^2 h_t(x,x).
\ee
By the Minakshisundaram-Pleijel expansion in (\ref{Minak-Pleijel-1})-(\ref{Minak-Pleijel-2}), we can see that,
\be
\lim_{t\to 0} \left(\sqrt{2}(4\pi)^{n/4} t^{(n+2)/4}\right)^2 h_t(x,x)
=\lim_{t\to 0} 2 (4\pi)^{n/2} t^{n/2+1} \frac{e^{-d_M^2(x,x)/(4t)}}{(4\pi t)^{n/2} } =0.
\ee
Thus  
\be \label{Diff-Distance-to-0}
\lim_{t\to 0} \lim_{N\to \infty} d_{{\mathbb E}^N}( \psi_t^N(x),\psi_t^N(y))
=\lim_{t\to 0} d_{\ell^2}( \psi_t(x),\psi_t(y))=0.
\ee

\subsection{Other Rescalings of the Heat Kernel Embedding}

B\'erard-Besson-Gallot  defined an embedding where they rescaled by volume
\be\label{BBG-It}
I_t(x)=\sqrt{\vol(M)} \Psi_t(x),
\ee
in Part III of \cite{BBG:1994}. 

They also defined an
embedding into an infinite dimensional sphere in Part VIII of \cite{BBG:1994}.
\be
K_t(x)=\Psi_t(x)/|\Psi_t(x)|
\ee
They proved the pull back satisfies
\be
K_{t}^*g_{\ell^2}=1/(2t) \left(g-(t/3)\Ricci+O(t^2)\right) \textrm{ as }t\to 0^+.
\ee

\subsection{Spectral Convergence of Riemannian Manifolds}\label{sect-spectral-conv}

Inspired by the work of Gromov \cite{Gromov-text} and of Fukaya \cite{Fukaya:1987aa},
 B\'erard-Besson-Gallot defined a spectral distance between Riemannian manifolds.

 To define a {\bf {\em spectral distance}}
between two Riemannian manifolds, $M_1$ and $M_2$, they
took the Hausdorff distance between $I_t(M_1)$ and $I_t(M_2)$ of (\ref{BBG-It})
and considered
a minmax over various reorderings of the eigenfunctions in Part IV of \cite{BBG:1994}. 
 They proved a compactness theorem
for their spectral distance and they prove that when sequences of manifolds converge with respect to this spectral distance their eigenvalues converge.  

Kasue and Kumura defined a different spectral convergence  in \cite{Kasue-Kumura:1994}  
\cite{Kasue_Kumura:2002}
\cite{Kasue:2002}.  Their notion is related to this one but is defined in a completely different way using maps
that almost preserve heat kernels rather than Hausdorff distances between embeddings.   They also proved that when sequences of manifolds converge with respect to this spectral distance their eigenvalues converge.

\section{\bf Applications to Dimension Reduction of High Dimensional Data Sets}

\begin{figure}[h] %  figure placement: here, top, bottom, or page
   \centering
   \includegraphics[width=2in]{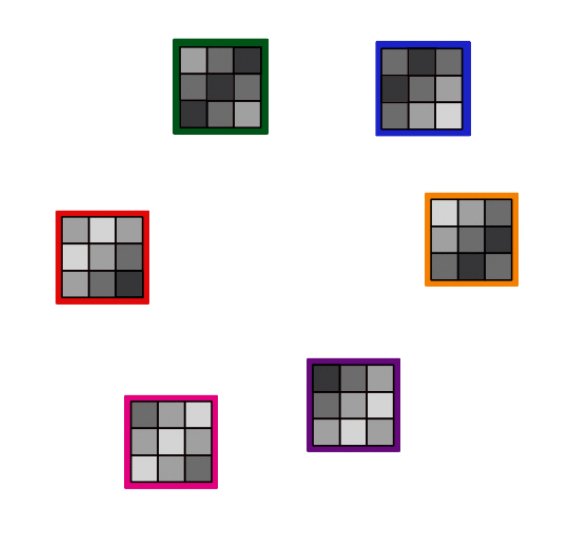} 
   \caption{Here we have $6$ photos each of which has $9$ pixels, 
   so they may be viewed as points in ${\mathbb R}^9$.  
   For example, the leftmost photo is $(2,1,2,1,2,3,2,3,4)$ and the 
one above it is $(2,3,4,3,4,3,4,3,2)$.   
 The distance between these photos measured using the Euclidean distance is
   then $\sqrt{24}$.
   }
   \label{fig:photos1}
\end{figure}

Heat kernel embeddings have been applied to achieve {\bf{\em dimension reduction}}: the mapping of data sets of points in high dimensional spaces to lower dimensional spaces.  An example of such a data set would be a collection of photographs.  Each photograph would be a point.  See for example Figure~\ref{fig:photos1}.  If the photographs are in grey scale with $1,000\times1,000$ pixels, then their points lie in 
$X={\mathbb R}^{1,000,000}$.   In general, the data set is described by $\{x_1,...x_k\}\in X$, where
$X$ is the high dimensional space endowed with a distance, $d_X$, that is meaningful to the application.   

\begin{figure}[h] %  figure placement: here, top, bottom, or page
   \centering
  \includegraphics[width=2in]{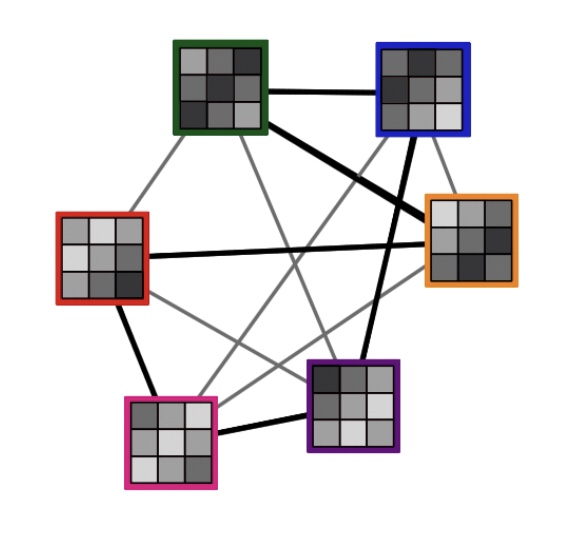} 
   \caption{We have created a graph with a photo at each vertex.   Here we have connected photos that are within an 
   $\epsilon=5$, so that there are edges between pairs of photos that are a distance less than $5$ apart.  Note it turns out there are 4 edges per vertex, so this is also a kNN graph with $k=4$.   Note that the thick lines connect vertices that are a distance $3$ apart and the thinner lines are for the vertices that are further apart.  
    }
   \label{fig:photos2}
\end{figure}

The distances between the points can be stored as a {\bf distance matrix}, $d_X(i,j)$.   Sometimes the distances are just the extrinsic Euclidean distances.   A neighborhood graph can be created by joining some nearby points with edges.
One method is to select the $k$ nearest neighbors (kNN) for each point to be connected with an edge so the resulting graph has degree $k$: $k$ edges per vertex.  Another method is to connect each point with every point in an $\epsilon$ neighborhood about it.  See Figure~\ref{fig:photos2}.   A distance of shortest paths can then be found using an algorithm by Dijkstra \cite{Dijkstra:1959} or we can use the extrinsic Euclidean distance.  

In 2000, Tenenbaum-deSilva-Langford \cite{Tenenbaum_deSilva_Langford:2000} created the ISOmap algorithm combining these ideas with Multidimensional Scaling (MDS) which is a technique from Linear Algebra.  Roweis-Saul \cite{Roweis_Saul:2000} introduced Nonlinear Dimensionality Reduction using the Locally Linear Embedding (LLE) algorithm which eliminated the need to estimate pairwise distances between widely separated points.   These initial papers introduced algorithms and demonstrated their effectiveness in studying sample data sets but did not include proofs.

When we assume the data set of points lies on a submanifold of the high dimensional space, the process of dimension reduction is called {\bf {\em manifold learning}}.   For example, we might have a collection of photos from a continuous film, in which case they lie on a curve which is a 1 dimensional manifold.   Or they might be photos of the same scene taken from many directions at the same distances, in which case the underlying submanifold would be $SO(3)$.   Note that the submanifold is not known: it is just presumed to exist, and then if the right technique is used, its properties allow us to reduce the dimension.   

\begin{figure}[h] %  figure placement: here, top, bottom, or page
   \centering
  \includegraphics[height=2in]{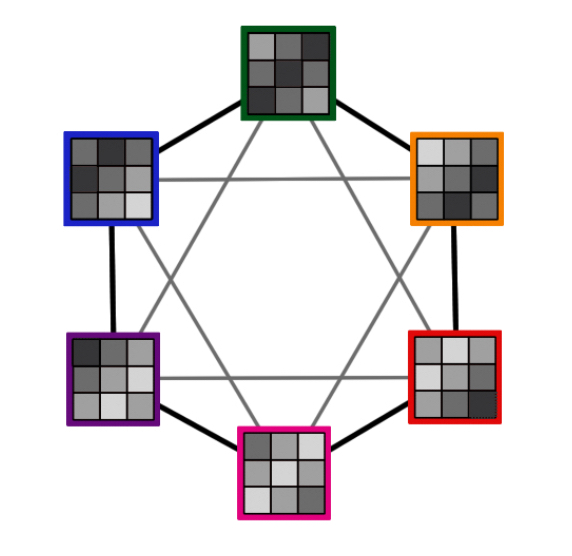}   \includegraphics[height=2.1in]{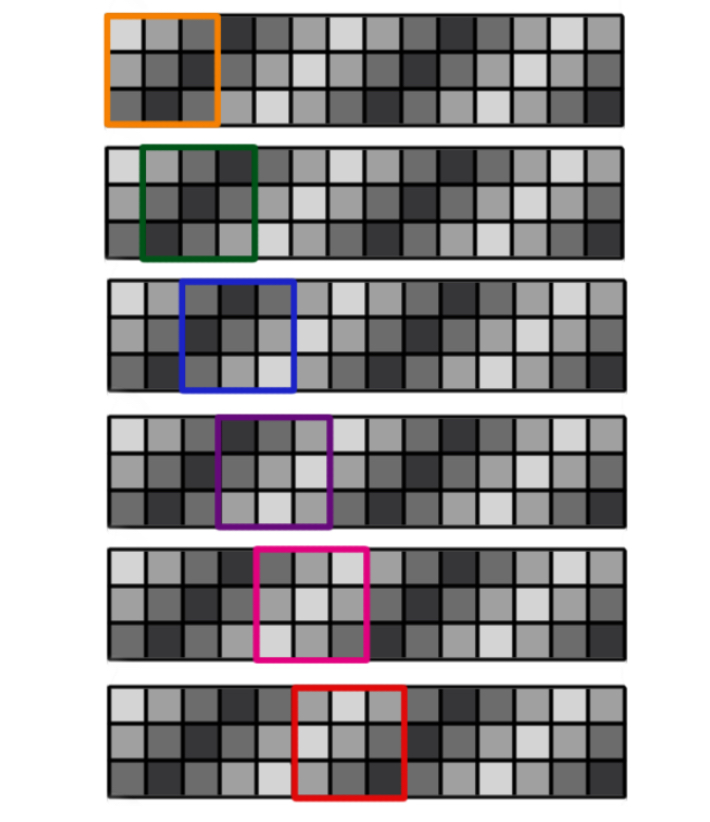} 
   \caption{On the left is the image of the data set from Figure~\ref{fig:photos1} after it has been mapped into ${\mathbb R}^2$ by an eigenmap.  We have placed each photo on top of its image under the eigenmap
   of Coifman-Lafon and we then drew the edges connecting them. We see that the data is now arranged symmetrically into a circle.   On the right, we have used the circular symmetry to glue the photos together to form
 a panoramic image which is periodic.  This symmetry was not visible in the original disorganized data set.}
   \label{fig:photos3}
\end{figure}

The goal is to map the data points into a lower dimensional space and thus discover the underlying structure.   Ideally we might even preserve the distances under this map (just as in a heat kernel embedding), although even just obtaining an embedding is of interest.  The heat kernel is also useful because it is robust to noise. In 2001, Belkin and Niyogi introduced the idea of embedding data sets using {\em eigenmaps} whose components are eigenfunctions \cite{Belkin-Niyogi:2001}.  In this section we describe their approach and then later we describe the 
{\em diffusion maps} of Coifman-Lafon whose components are weighted eigenfunctions \cite{Coifman-Lafon:2006}.
See Figure~\ref{fig:photos3}.

\subsection{Belkin-Niyogi Eigenmaps of Data Sets} \label{sect-Belkin-Niyogi}

In 2001, Belkin and Niyogi introduced {\em Laplacian eigenmaps} \cite{Belkin-Niyogi:2001} as a method to provide dimension reduction for data sets that lie on submanifolds in a high dimensional space, $X$.   Their algorithm has three steps as follows:
%In a series of papers published in 2001-2006 Belkin and Niyogi introduced {\em Laplacian eigenmaps} \cite{Belkin-Niyogi:2001} \cite{Belkin-Niyogi:2003} \cite{Belkin:2003} \cite{Belkin-Niyogi:2005}  \cite{Belkin-Niyogi:2006convergence}.    

\vspace{.2cm}
\noindent{\bf Algorithm:}

\vspace{.2cm}
\noindent
{\bf Step 1:} {\em Given $k$ points, 
$\{x_1,...,x_k\}\subset X$,  they constructed an adjacency graph connecting neighboring points, $x_i \sim x_j$, either using an $\epsilon$-neighborhood approach or the $k$ nearest neighbors (kNN).

\vspace{.2cm}
\noindent
{\bf Step 2:} They created an $k\times k$ {\bf \em heat kernel matrix},
\be\label{Wij}
W^t_{i,j}= e^{-|x_i-x_j|^2/t} \textrm{ for } x_i \sim x_j \textrm{ and } W^t_{i,j}=0 \textrm{ otherwise.}
\ee
In particular, $W^t_{i,i}=e^{-0^2/t}=1$ because $x_i\sim x_i$.   
They define a {\bf {\em diagonal degree matrix}} 
\be\label{degree-matrix}
D_{i,i}=\sum_{j=1}^n W_{i,j} \textrm{ where } W_{i,j}=W^t_{i,j},
\ee
a {\bf \em graph Laplacian matrix}, $L=D-W$.  
 
 \vspace{.2cm}
\noindent{\bf Step 3:}  Belkin-Niyogi define an {\bf {\em eigenmap}}: $V_t^N: \{x_1,...,x_k\} \to {\mathbb R}^N$, for any
$N\le k-1$ by
\be\label{BN-emap}
V_t^N(x_i)=(v_1(i),v_2(i),...,v_N(i)) \in {\mathbb R}^N \textrm{ where } v_j(i) \textrm{ is the } i^{th} \textrm{ component of } v_j.
\ee
where
$
0=\mu_0<\mu_1\le \cdots \le \mu_{k-1} 
$
are the eigenvalues and 
\be \label{BN-evectors}
v_0,v_1,..., v_{k-1} \in {\mathbb R}^k
\ee
is an orthonormal collection of {\bf \em eigenvectors} of the {\bf \em normalized graph Laplacian} $D^{-1}L$:
\be \label{BN-e}
D^{-1}Lv_j=\mu_j v_j. 
\ee
}
 \vspace{.2cm}

Belkin and Niyogi leave out the $0^{th}$ eigenvector $v_0$ in the definition of their eigenmap $V_t^N$ because all of its entries are the same.   This can easily be seen by observing
that the sum of the entries in each row of $D^{-1}L$ is $0$:
\be
\sum_{j=1}^k (D^{-1} L)_{i,j}=\sum_{j=1}^k D^{-1}_i (D-W)_{i,j}=\sum_{j=1}^k I_{i,j} -D^{-1}_i \sum_{j=1}^k W_{i,j}=1-D^{-1}_iD_i=0.
\ee

In \cite{Belkin-Niyogi:2001}, Belkin and Niyogi explain that intuitively if $\{x_1,...x_k\}\subset X$ lie as a uniformly distributed collection of points in a submanifold, $M^n$, of small dimension, $n$, in the high dimensional Euclidean space, $X$, then they can use eigenfunctions, $\phi_j: M \to {\mathbb R}$, of the Laplacian, $\Delta$, on $M$ to
define a {\bf \em smooth eigenmap},
\be\label{BN-smooth-eigenmap}
\Phi^N: M^n \to {\mathbb R}^N \textrm{ as } \Phi^N(x_i)=(\phi_1(x_i),\phi_2(x_i),...,\phi_N(x_i)).
\ee
Note that this {\bf \em smooth eigenmap}, $\Phi^N(x)$, is a truncation of an unweighted version of B\'erard-Besson-Gallot's heat kernel
embedding, $\Psi_t(x)$, that we saw in (\ref{BBG-map}) \cite{BBG:1994}.   

Belkin and Niyogi's three step algorithm defines an eigenmap, $V_t^N$,  on the data set, $\{x_1,...,x_k\}$, without knowing the ambient manifold, $M$, such that
\be\label{BN-Phi-approx-V}
\Phi^N(x_i)=(\psi_1(x_i),\psi_2(x_i),...,\psi_N(x_i))\approx (v_1(i),v_2(i),...,v_N(i))=V_t^N(x_i)
\ee  
The approximation $\Phi^N ( x_i ) \approx V_t^N(x_i)$ was not proven rigorously in Belkin and Niyogi's early papers however they did provide some of intuition leading to their algorithm.

The intuitive reasoning behind Belkin and Niyogi's choice of matrix in (\ref{Wij}) was explained by Belkin roughly as follows
in his doctoral dissertation \cite{Belkin:2003}.   Recall that given a function, $f: M \to {\mathbb R}$,
then 
\be
u(x,t)=\int_M h_t(x,y) f(y) \,dvol_y
\ee
 satisfies 
 \be
 \partial_t u=\Delta u \textrm{ and } \lim_{t\to 0^+}u(x,t)=f(x). 
 \ee
 Thus
\begin{eqnarray}
(\Delta f) (x) &=& (\Delta u)(x,0)= (\partial_t u)(x,0)= \lim_{t\to 0^+} \tfrac{-1}{t} \left(u(x,0)-u(x,t)\right)\\
&=&\lim_{t\to 0^+} \tfrac{-1}{t} \left( f(x)-\int_M h_t(x,y) f(y) \, dvol_y \right)\\
\end{eqnarray}
Since they are studying a data set, Belkin and Niyogi ignored the limit as $t\to 0^+$ and consider integration to be summation:
\be 
(\Delta f) (x_i)\approx f(x_i) - \sum_{j=1}^k h_t(x_i,x_j)f(x_j)
\ee
They decided to choose a weight matrix $W_{i,j}$ so that 
\be \label{DinvWapproxHeat}
D^{-1}_{i,i}W_{i,j} \approx h_t(x_i,x_j) 
\ee
so that
\be 
(\Delta f) (x_i)\approx f(x_i) - \sum_{j=1}^k (D^{-1}W)_{i,j} f(x_j) = (I-D^{-1}W)_{i,j} f(x_j).
\ee
In this way, the eigenfunctions $\phi_j(x)$ of the Laplace Beltrami operator, $\Delta \phi_j =-\lambda_j \phi_j$,
seem to correspond with the eigenvectors, $v_j$, of $I-D^{-1}W=D^{-1}L$ as in (\ref{BN-e}) so that
\be \label{BN-Phi-approx-V-2}
\phi_j(x_i) \quad \approx \quad v_j(i) \qquad \forall i,j\in\{1,..,k\} %\textrm{ and possibly } \lambda_j \approx \mu_j 
\ee
which intuitively justifies (\ref{BN-Phi-approx-V}).
Note that they made no claim that the eigenvalues match
and in practice the eigenvalues, $\lambda_j$ and $\mu_j$ 
are not close even as $t\to 0^+$ as we will see in upcoming sections.   

Note that by (\ref{DinvWapproxHeat}) and the definition of the degree matrix (\ref{degree-matrix}).
they achieve the normalization of the heat kernel:
\be
1=\int_M h_t(x,y) \, dy \quad \approx \quad \sum_{j=1}^k h_t(x_i,x_j) =\sum_{j=1}^k D_i^{-1\,\,}W_{i,j}=1.
\ee
Finally Belkin-Niyogi approximated the heat kernel as we've seen throughout by
\be
h_t(x_i,x_j) \approx \tfrac{1}{(4\pi t)^{n/2} } \,  e^{-|x_i - y_i |^2/(4t) }
\ee
to justify their choice of $W_{i,j}$ in (\ref{Wij}) although their choice has division by t instead of $4t$ in the exponent.   
Since this just rescales time, it does not cause any issues for them.

A more rigorous discussion appears in the 2006 papers of Belkin-Niyogi using a probabilistic approach \cite{Belkin-Niyogi:2005} \cite{Belkin-Niyogi:2006convergence} .   See also Kondor-Lafferty work in \cite{Kondor-Lafferty:2002}.
% where they also multiply their matrix $W$ by $(1/n)1/(4\pi t)^{dim(M)/2}$.   
See also their joint work with Sinddhwani in
\cite{Belkin-Niyogi-Sinddhwani:2006}.

\subsection{A Brief History of Graph Laplacians starting with Dodziuk}

It should be noted that the graph Laplacian was first defined by Dodziuk in 1984 \cite{Dodziuk:1984}.   His graph Laplacian was on a simple graph without weighted edges.  He used the {\bf difference Laplacian},
\be
``\Delta'' f (x_i) =\sum_{x_j \sim x_i} (f(x_j)-f(x_i))= -D_i f(x_i) + \sum_{x_j \sim x_i} f(x_j)
\ee
that had been studied on square lattices in Euclidean space as early as 1928 by Courant, Friedrichs and Lewy
\cite{CFL:1928}.    Here $D_i$ is the degree of vertex $x_i$ and is equal to the number of vertices $x_j$ such that 
$x_j \sim x_i$.   Defining $W$ to be the adjacency matrix such that $W_{i,j}=1$ when $x_i=x_j$ or there
is an edge between $x_i$ and $x_j$, and $0$ otherwisse, we see that Dodziuk's graph Laplacian satisfies,
\be
``\Delta'' f (x_i) =\sum_{x_j \sim x_i} (f(x_j)-f(x_i))=  \sum_{j=1}^k (W-D)_{i,j}f(x_j)   
\ee
where $D$ is the diagonal degree matrix $D_{i,i}=D_i$.   Note the opposite sign on this version.
In his paper, Dodziuk proved the maximum principle, the Harnack inequality, and Cheeger's bound for the lowest eigenvalue all extend to the graph setting in \cite{Dodziuk:1984}.     In 1988, he developed a Green's formula with Karp in \cite{Dodziuk-Karp:1988}.   For early applications of these ideas to Markov chains see the paper of Lawler and Sokol \cite{Lawler-Sokol:1988}.

\subsection{Continuing the History of Graph Laplacians with work of Chung}
A deep analysis of graph Laplacians  appears in the work of Chung in  \cite{Chung:1989} where
she studies weighted graphs with a symmetric weight matrix, $W$.  The degree of the vertex is then
defined to be
$
D_i=\sum_{j=1}^k W_{i,j} 
$
Chung's Graph Laplacian matrix is defined $L=D-W$ with the opposite sign used by Dodziuk.   She defines a normalized graph Laplacian matrix, $D^{-1}L$, which might not be symmetric.  So
she defines a normalized symmetrized graph Laplacian matrix,
\be\label{graph-laplacian}
{\mathcal L}=D^{1/2\,\,} (D^{-1}L )\,D^{-1/2\,\,}=D^{-1/2\,\,} L D^{-1/2\,\,}= I -D^{-1/2\,\,} W D^{-1/2\,\,}, 
\ee
which is a symmetric matrix.   
Taking a function $f: \{x_1,...,x_k\} \to {\mathbb R}$,
she defines $({\mathcal L} f):  \{x_1,...,x_k\}  \to {\mathbb R}$ as follows:
 \be
({\mathcal L} f )(x_i) = \sum_{j=1}^k (I -D^{-1/2\,\,} W D^{-1/2\,\,})_{i,j} \,f(x_j),
\ee
so that
 \be
({\mathcal L} f )(x_i) 
\,\,\,=\,\,\, f(x_i) -\sum_{j=1}^k D_i^{-1/2\,\,}W_{i,j}\,D_j^{-1/2\,\,}f (x_j).
\ee

Note that the eigenvalues $\mu_j$ of the normalized graph Laplacian, $D^{-1}L$,  are also the eigenvalues of 
${\mathcal L}$ because
\begin{eqnarray}
{\mathcal L}(D^{1/2}v_j)&=&D^{-1/2\,\,} L D^{-1/2\,\,}(D^{1/2}v_j)
=D^{1/2\,\,}D^{-1} Lv_j\\
&=&D^{1/2\,\,}(-\mu_j v_j)=-\mu_j (D^{1/2}v_j).
\end{eqnarray}
Because of this well known relationship, some people redefine Belkin-Niyogi's eigenmaps using the eigenvectors of
the symmetrized normalized Laplacian $\mathcal L$.

In 1995-1997 Chung completed 
joint work with Yau on Graph eigenvalues and Sobolev inequalities \cite{Chung-Yau:1995} , 
joint work with Grigor'yan and Yau proving upper bounds on eigenvalues of finite graphs \cite{Chung-Grigorian-Yau:1996}, and joint work with Yau on the heat kernel and its trace for graphs  \cite{Chung-Yau:1997trace}.
She has applied these ideas in many directions including developing a pagerank for web search engines  \cite{Chung:pagerank}.   

Note that Belkin-Niyogi's graph Laplacian is defined using their own weight
matrix $W$ that was not studied by Chung.  It does not agree with the adjacency matrix of Dodziuk even
when all edges have the same weight.   Their graph Laplacian only fits with Chung's work in the sense that it uses her notion of a weighted
graph Laplacian (as in 1.4 of her textbook \cite{Chung:text}).   It is not the heat kernel matrix defined in Chung's work with Yau in  \cite{Chung-Yau:1997trace} and explored further in the later chapters of her textbook \cite{Chung:text}.   

\subsection{Graph Laplacians Approximating Smooth Laplacians}

Hein, Audibert, and Luxburg provide a rigorous discussion of the relationship between graph Laplacians to manifold Laplacians in \cite{Hein_Audibert_Luxburg:2005}.
Singer improves the convergence rate from graph to manifold Laplacians in \cite{Singer:2006} by also assuming uniformly distributed points.

In \cite{WMKG-2007}, Wardetzky, Mathur, Kalberer, and Grinspun study five different discrete Laplacian and discuss which properties each shares with the continuous Laplacian.  They conclude that no discrete Laplacian can have all the properties but that different ones may be used to achieve different properties.

Most recently, Burago-Ivanov-Kurylev have shown in \cite{Burago_Ivanov_Kurylev:2015} that eigenvalues and eigenfunctions of the Laplace-Beltrami operator on a Riemannian manifold are approximated by eigenvalues and eigenvectors of a (suitably weighted) graph Laplace operator of a proximity graph on an epsilon-net.   Earlier work
by Fujiwara in \cite{Fujiwara:1995} had also considered $\epsilon$-nets but he only approximated the eigenvalues and his estimates were not as strong.  There is also unpublished work of Aubry from 2014  providing estimates that depend only on the upper bounds on diameter and sectional curvature of $M$ and lower bounds on the injectivity radius.   

\subsection{An Eigenmap of a Hexagon}

As a simple case let us work out an eigenmap of $6$ points, $\{x_1, x_2, x_3, x_4, x_5, x_6\}$, where
\be
x_j=(\cos(2\pi j/6), \sin(2\pi j/6))
\ee  lie on the unit circle.   
We can consider $kNN$ with $k=2$
nearest points or $\epsilon$ slightly $>1$ to determine the same collection of unit length edges, so that
our graph is the standard hexagon with two edges of length one attached at every vertex.
Then by (\ref{Wij}), we have
\be
W= \begin{pmatrix}
1 & e^{-1/t} & 0 & 0& 0& e^{-1/t}  \\
e^{-1/t} & 1 & e^{-1/t}  & 0 & 0 & 0\\
0  & e^{-1/t}   & 1& e^{-1/t}  &0 &0 \\
0  &  0& e^{-1/t}  & 1 &e^{-1/t}  & 0\\
0& 0  &  0& e^{-1/t}  & 1 &e^{-1/t}  \\
e^{-1/t}&0& 0  &  0& e^{-1/t}  & 1 \\
\end{pmatrix}.
\ee
and our diagonal matrix $D$ has $D_{i,i}=\sum_{j=1}^n W_{i,j}=1+2e^{-1/t}$.  So
\be
L=D-W=
\begin{pmatrix}
2e^{-1/t} & -e^{-1/t} & 0 & 0& 0& -e^{-1/t}  \\
-e^{-1/t} & 2e^{-1/t} & -e^{-1/t}  & 0 & 0 & 0\\
0  & -e^{-1/t}   & 2e^{-1/t}& -e^{-1/t}  &0 &0 \\
0  &  0& -e^{-1/t}  & 2e^{-1/t} &-e^{-1/t}  & 0\\
0& 0  &  0& -e^{-1/t}  & 2e^{-1/t} &-e^{-1/t}  \\
-e^{-1/t}&0& 0  &  0& -e^{-1/t}  & 2e^{-1/t} \\
\end{pmatrix}.
\ee
Here $D$ is a multiple of the identity, so we have
${\mathcal L}=D^{-1/2\,\,}L D^{-1/2\,\,}=D^{-1}L$, and
\be
{\mathcal L}=D^{-1}L=\frac{2e^{-1/t}}{1+2e^{-1/t}}
\begin{pmatrix}
1 & -\tfrac{1}{2}\,\,& 0 & 0& 0& -\tfrac{1}{2}\,\,  \\
-\tfrac{1}{2}\,\, & 1 & -\tfrac{1}{2}\,\, & 0 & 0 & 0\\
0  & -\tfrac{1}{2}\,\,  & 1& -\tfrac{1}{2}\,\, &0 &0 \\
0  &  0& -\tfrac{1}{2}\,\, & 1 & -\tfrac{1}{2}\,\,  & 0\\
0& 0  &  0& -\tfrac{1}{2}\,\, & 1 & -\tfrac{1}{2}\,\,  \\
-\tfrac{1}{2}\,\,&0& 0  &  0& -\tfrac{1}{2}\,\,  & 1 \\
\end{pmatrix}
\ee
Thus we will get the same eigenmaps, $V_t^N:\{x_1,...,x_k\}\to {\mathbb R}^N$, defined using
Belkin-Niyogi's $D^{-1}L$ as we get using the symmetrized $\mathcal L$ to compute the eigenvectors.

We know the $\lambda_0=0$ eigenfunction is constant and quickly confirm the same for eigenvectors of ${\mathcal L}$:
\be
\frac{2e^{-1/t}}{1+2e^{-1/t}}
\begin{pmatrix}
1 & -\tfrac{1}{2}\,\,& 0 & 0& 0& -\tfrac{1}{2}\,\,  \\
-\tfrac{1}{2}\,\, & 1 & -\tfrac{1}{2}\,\, & 0 & 0 & 0\\
0  & -\tfrac{1}{2}\,\,  & 1& -\tfrac{1}{2}\,\, &0 &0 \\
0  &  0& -\tfrac{1}{2}\,\, & 1 & -\tfrac{1}{2}\,\,  & 0\\
0& 0  &  0& -\tfrac{1}{2}\,\, & 1 & -\tfrac{1}{2}\,\,  \\
-\tfrac{1}{2}\,\,&0& 0  &  0& -\tfrac{1}{2}\,\,  & 1 \\
\end{pmatrix}
\,
\begin{pmatrix}
1\\
1\\
1\\
1\\
1\\
1\\
\end{pmatrix}
=
\begin{pmatrix}
0\\
0\\
0\\
0\\
0\\
0\\
\end{pmatrix}
=
0
\begin{pmatrix}
1\\
1\\
1\\
1\\
1\\
1\\
\end{pmatrix}
\ee
Ordinarily we would next compute the eigenvectors numerically, but here we can quickly guess them
using Belkin-Niyogi's intuition as described in (\ref{BN-Phi-approx-V}): $v_j(i)\approx \phi_j(x_i)$.

In Section~\ref{sect-circle-map}, we found the eigenfunctions
\be
\phi_1(\theta_i)=A_1\cos(\theta_i)\textrm{ and } \phi_2(\theta_i))= A_1\sin(\theta_i).
\ee
By (\ref{BN-Phi-approx-V}), we can guess that the next two eigenvectors, $v_1$ and $v_2$, are
$$
\begin{pmatrix}
v_1(1)\\
v_1(2)\\
v_1(3)\\
v_1(4)\\
v_1(5)\\
v_1(6)\\
\end{pmatrix}
=
\begin{pmatrix}
\cos(2\pi/6)\\
\cos(2\pi 2/6)\\
\cos(2\pi 3/6)\\
\cos(2\pi 4/6)\\
\cos(2\pi 5/6)\\
\cos(2\pi 6/6)\\
\end{pmatrix}
=
\begin{pmatrix}
\tfrac{1}{2}\,\\
-\tfrac{1}{2}\,\,\\
-1\\
-\tfrac{1}{2}\,\,\\
\tfrac{1}{2}\,\\
1
\end{pmatrix}
\textrm{ and }
\begin{pmatrix}
v_2(1)\\
v_2(2)\\
v_2(3)\\
v_2(4)\\
v_2(5)\\
v_2(6)\\
\end{pmatrix}
=
\begin{pmatrix}
\sin(2\pi/6)\\
\sin(2\pi 2/6)\\
\sin(2\pi 3/6)\\
\sin(2\pi 4/6)\\
\sin(2\pi 5/6)\\
\sin(2\pi 6/6)\\
\end{pmatrix}
=
\begin{pmatrix}
\tfrac{\sqrt{3}}{2}\\
\tfrac{\sqrt{3}}{2}\\
0\\
-\tfrac{\sqrt{3}}{2}\\
-\tfrac{\sqrt{3}}{2}\\
0
\end{pmatrix}
$$
where $\sqrt{3}$ is needed to make them orthonormal.
Checking the first we see:
$$
\begin{pmatrix}
1 & -\tfrac{1}{2}\,\,& 0 & 0& 0& -\tfrac{1}{2}\,\,  \\
-\tfrac{1}{2}\,\, & 1 & -\tfrac{1}{2}\,\, & 0 & 0 & 0\\
0  & -\tfrac{1}{2}\,\,  & 1& -\tfrac{1}{2}\,\, &0 &0 \\
0  &  0& -\tfrac{1}{2}\,\, & 1 & -\tfrac{1}{2}\,\,  & 0\\
0& 0  &  0& -\tfrac{1}{2}\,\, & 1 & -\tfrac{1}{2}\,\,  \\
-\tfrac{1}{2}\,\,&0& 0  &  0& -\tfrac{1}{2}\,\,  & 1 \\
\end{pmatrix}
\,
\begin{pmatrix}
\tfrac{1}{2}\,\\
-\tfrac{1}{2}\,\,\\
-1\\
-\tfrac{1}{2}\,\,\\
\tfrac{1}{2}\,\\
1
\end{pmatrix}
=
\begin{pmatrix}
\tfrac{1}{2}\, + \tfrac{1}{4}\, -\tfrac{1}{2}\,\,\\
-\tfrac{1}{4}\,-\tfrac{1}{2}\,\,+\tfrac{1}{2}\,\\
\tfrac{1}{4}\,-1+\tfrac{1}{4}\,\\
\tfrac{1}{2}\,-\tfrac{1}{2}\,\,-\tfrac{1}{4}\,\\
\tfrac{1}{4}\,+\tfrac{1}{2}\,-\tfrac{1}{2}\,\,\\
-\tfrac{1}{4}\,-\tfrac{1}{4}\,+1\\
\end{pmatrix}
=
\begin{pmatrix}
\tfrac{1}{4}\,\\
-\tfrac{1}{4}\,\\
-\tfrac{1}{2}\,\,\\
-\tfrac{1}{4}\,\\
\tfrac{1}{4}\,\\
\tfrac{1}{2}\,\\
\end{pmatrix}
=
\frac{1}{2}
\begin{pmatrix}
\tfrac{1}{2}\,\\
-\tfrac{1}{2}\,\,\\
-1\\
-\tfrac{1}{2}\,\,\\
\tfrac{1}{2}\,\\
1
\end{pmatrix}
$$
so that 
\be
{\mathcal L} \,\,v_1\,= \,\frac{2e^{-1/t}}{\,1+2e^{-1/t}\,} \,\frac{1}{2}\, v_1.
\ee
Note that the eigenvalue itself does not match the eigenvalue of the Laplacian
even though the eigenvector matches the eigenfunction so well.

Next we see that the truncated eigenmap, with $N=2$, is just a rescaling of the heat kernel embedding, with $N=2$:
\be
V_t^2(x_i)=(v_1(i),v_2(i))=(\cos(2\pi i/6),\sin(2 \pi i/6)) =\, e^t\, \Psi^2_t(x_i ).
\ee
The fact that we have equality instead of an approximate equality is due to the symmetry of the circle and the points lying on the circle.  %I had computed the next eigenvector which can be found in TWAS9

\subsection{Finding the Eigenmap of a Data Set of $3 \times 3$ Photos}\label{sect-photos}

\begin{figure}[h] %  figure placement: here, top, bottom, or page
   \centering
\includegraphics[width=2in]{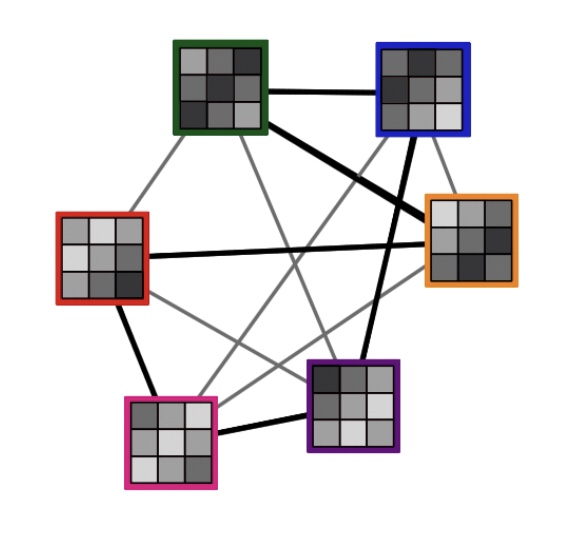}\qquad 
\qquad \includegraphics[width=2in]{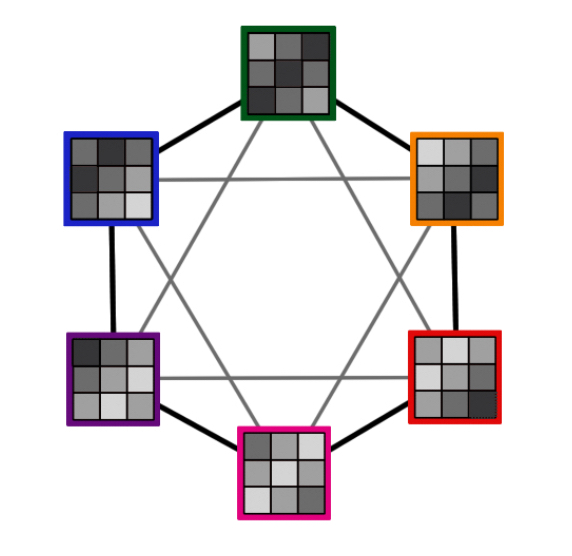} 
   \caption{On the left is the data set in ${\mathbb R}^9$ with the points marked by their $3\times 3$ photos
   and on the right is the image under
    the truncated eigenmap, $V_t^2$, with the image points marked by their $3\times 3$ photos lying 
    on ${\mathbb S}^1\subset {\mathbb R}^2$.}     
   \label{fig:photos23}
\end{figure}

At the beginning of this section we described  the process of dimension reduction using a collection of six $3\times 3$ pixel photos.   See Figure~\ref{fig:photos23}.  We will now
describe exactly how this is done using a Belkin-Niyogi eigenmap.

We label the original collection of  six $3\times 3$ pixel photos on the left side of Figure~\ref{fig:photos23}
starting with the rightmost and proceeding counter clockwise:
$$
\textcolor{orange}{p_1=
\begin{pmatrix}
1\\2\\3\\2\\3\\4\\3\\4\\3\\
\end{pmatrix}}
\quad
\textcolor{blue}{p_2=
\begin{pmatrix}
3\\4\\3\\4\\3\\2\\3\\2\\1\\
\end{pmatrix}}
\quad
\textcolor{teal}{p_3=
\begin{pmatrix}
2\\3\\4\\3\\4\\3\\4\\3\\2\\
\end{pmatrix}}
\quad
\textcolor{red}{p_4=
\begin{pmatrix}
2\\1\\2\\1\\2\\3\\2\\3\\4\\
\end{pmatrix}}
\quad
\textcolor{magenta}{p_5=
\begin{pmatrix}
3\\2\\1\\2\\1\\2\\1\\2\\3\\
\end{pmatrix}}
\quad
\textcolor{violet}{p_6=
\begin{pmatrix}
4\\3\\2\\3\\2\\1\\2\\1\\2\\
\end{pmatrix}}
\quad
$$

Evaluating the distances between these points we see that some have distance $\sqrt{24}$ and
others have every pixel off by 1 so their distance is $\sqrt{9}=3$.   Taking $t=1$ and observing that
$e^{-24}\approx 0$ and $e^{-3}\approx .05$ we get a $6\times 6$ matrix
\be
W=
\begin{pmatrix}
1&0&.05&.05&0&0\\
0&1&.05&0&0&.05\\
.05&.05&1&0&0&0\\
.05&0&0&1&.05&0\\
0&0&0&.05&1&.05\\
0&.05&0&0&.05&1\\
\end{pmatrix}
\ee
Thus our diagonal matrix has $D_{i,i}=1+.05+.05=1.1=11/10$ so
\be
L=D-W=
\begin{pmatrix}
.1&0&-.05&-.05&0&0\\
0&.1&-.05&0&0&-.05\\
-.05&-.05&.1&0&0&0\\
-.05&0&0&.1&-.05&0\\
0&0&0&-.05&.1&-.05\\
0&-.05&0&0&-.05&.1\\
\end{pmatrix}
\ee
and since $D^{-\tfrac{1}{2}\,\,}_{i,i}= (11/10)^{-1/2}=\sqrt{10/11}$ we have
\be
{\mathcal L}=D^{-\tfrac{1}{2}\,\,}L D^{-\tfrac{1}{2}\,\,}= \frac{10}{11}
\begin{pmatrix}
1&0&-1/2&-1/2&0&0\\
0&1&-1/2&0&0&-1/2\\
-1/2&-1/2&1&0&0&0\\
-1/2&0&0&1&-1/2&0\\
0&0&0&-1/2&1&-1/2\\
0&-1/2&0&0&-1/2&1\\
\end{pmatrix}
\ee
It is easy to check that
\be
{v}_1=
\begin{pmatrix}
\tfrac{\sqrt{3}}{2} \\
-\tfrac{\sqrt{3}}{2} \\
0 \\
\tfrac{\sqrt{3}}{2} \\
0 \\
-\tfrac{\sqrt{3}}{2} \\
 \end{pmatrix}
 \textrm{ and }
{v}_2=
\begin{pmatrix}
1/2 \\
1/2 \\
1 \\
-1/2 \\
-1 \\
-1/2 \\
 \end{pmatrix}
\ee
are perpendicular eigenfunctions with eigenvalue $=-10/22$.  Thus up to scaling,
Belkin-Niyogi's eigenmap, $V^2_t(p_i)=(v_1(p_i),v_2(p_i))$, takes the following values:
\begin{eqnarray*}
\textcolor{orange}{V_2(p_1)=(\sqrt{3}/2,1/2)} & 
\textcolor{blue}{V_2(p_2)=(-\sqrt{3}/2,1/2) } &
\textcolor{teal}{V_2(p_3)=(0,1) } \\
\textcolor{red}{V_2(p_4)=(\sqrt{3}/2,-1/2) } &
\textcolor{magenta}{V_2(p_5)=(0,-1) } &
\textcolor{violet}{V_2(p_6)=(-\sqrt{3}/2,-1/2) } 
\end{eqnarray*}
which reveals the circular structure of the image seen on the right in Figure~\ref{fig:photos23}.

%removed subsection about eigenvalue computation

\subsection{An Eigenmap of a Sphere}\label{sect-emap-sphere}

Our 2020 student research team (Maziar Farahzad, Julinda Pillati Mujo, and Esteban Alcantara) tried out the Belkin-Niyogi method of eigenmaps using MATLAB.  Recall that in Section~\ref{sect-sphere-evalue} we 
observed that the first three nonconstant eigenfunctions of ${\mathbb S}^2$  up to scaling are
\be
\phi_1(x)= X(x) \qquad \phi_2(x)=Y(x) \qquad \phi_3(x)= Z(x)
\ee
so the image of the heat kernel embedding, $\Psi_t^3: {\mathbb S}^2 \to {\mathbb E}^3$, is a round sphere.    So this was a good example to test their code on an evenly distributed
data set of points in ${\mathbb S}^2$ to produce an eigenmap $V^3_t$ to ${\mathbb E}^3$.

\begin{figure}[h] %  figure placement: here, top, bottom, or page
   \centering
  \includegraphics[width=3in]{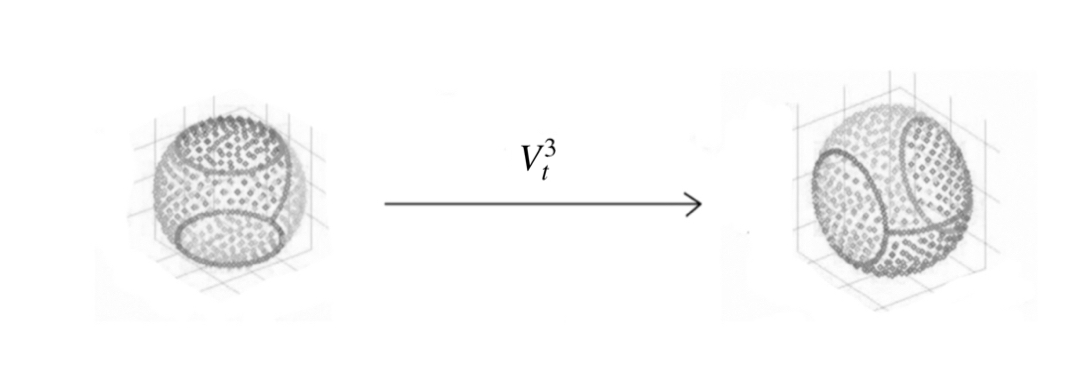} 
   \caption{Here we see a data set of evenly distributed points on a sphere mapped to a sphere by a truncated eigenmap,
   $V_t^3$, and we see that it agrees up to scaling with the image under a heat kernel embedding, 
   $\psi_t^3: {\mathbb S}^2\to {\mathbb E}^3$.   MATLAB computation by Maziar Farahzad, Julinda Pillati Mujo, and Esteban Alcantara.
   }
   \label{fig:sphere-eigenmap}
\end{figure}

To get the image to look like a sphere, the points had to be evenly distributed.  Ultimately the team devised a collection of biLipschitz disjoint charts to cover a sphere, used those charts to map a triangular lattice from a plane into the sphere, and then added extra points along the edges of the charts, to distribute the points as evenly as possible.  

In Figure~\ref{fig:sphere-eigenmap} we see the evenly distributed points on the sphere as our data set on the left, and then we see the image under $V_t^3$ on the right.   Since $N=3$ this agrees up to scaling with the image under the truncation of the sphere's heat kernel embedding, $\psi_t^3$, that we studied in Section~\ref{sect-sphere-evalue}.   It has been rotated because the eigenvectors were selected in a different order by MATLAB so that up to scaling,
\be
V_t^3(x_j)=(v_1(i),v_2(i),v_3(i)) \quad \leftrightarrow \quad \psi_t^3(x_j)=(-Z(x_j),-X(x_j), Y(x_j)).
\ee
This captures the concern B\'erard-Besson-Gallot mentioned regarding choice of eigenfunctions leading to rotations of the image of the heat kernel embeddings by elements of SO(N) that we discussed at the end of Section~\ref{sect-trunc-map}.

\subsection{Eigenmaps of Barbells} \label{sect-barbell}
\begin{figure}[h] %  figure placement: here, top, bottom, or page
   \centering
  \includegraphics[width=3in]{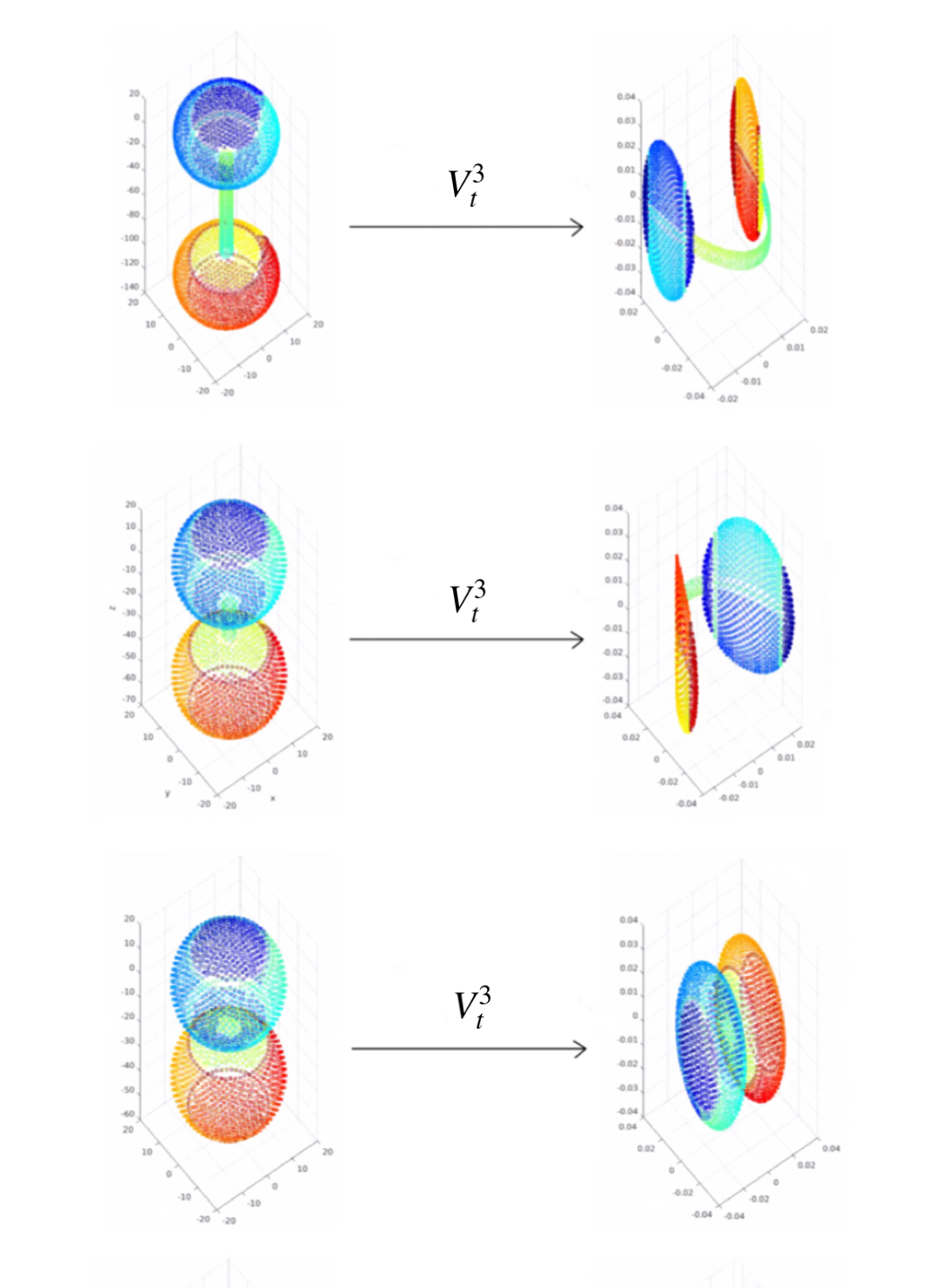} 
   \caption{
   Here we see a data set of evenly distributed points on various barbells mapped by an eigenmap,
   $V_t^3$, which gives some possible insight into the image under a heat kernel embedding, 
  $\Psi_t^3$ of these barbells.    MATLAB computation by Maziar Farahzad, Julinda Pillati Mujo, and Esteban Alcantara.
   }
   \label{fig:barbells}
\end{figure}

Our 2020 student research team (Maziar Farahzad, Julinda Pillati Mujo, and Esteban Alcantara) next applied the the method of eigenmaps to study what the images of heat kernel embeddings of deformed spheres shaped like barbells might look like.   In Figure~\ref{fig:barbells}, there are evenly distributed points on various barbells: {\em pairs of spheres of radius, $20$, with caps removed and cylinders of radius, $5$, of various lengths, $L=100$, $L=10$ and $L=1$, glued between them}.   

When the cylinder is long, $L=100$ as on the top of Figure~\ref{fig:barbells},
there is a first nonzero eigenvalue corresponding to the full height of the manifold which has an eigenfunction that looks something like $cos(z)$.  The next two eigenfunctions seem to correspond to $Z$ and $X$.   That is the points $x_i$ are mapped to $V_t^3(x_i)=(\cos(Z(x_i)), Z(x_i), X(x_i))$ and so the image of $\psi_t^3$ is probably not an embedding.  It is flattened and then bent around a cosine curve.   The students called this image a {\em headset }.  

In the second row of Figure~\ref{fig:barbells}, the cylinder only has $L=10$ and yet the same phenomenon has occurred.  The image is once again not an embedding and takes the form of a headset.   The headset is rotated 
the points $x_i$ seem to be mapped to 
\be
V_t^3(x_i)=(-\cos(Z(x_i)), -Z(x_i), X(x_i))
\ee
up to scaling because different eigenvectors were selected first by MATLAB. 

In the bottom row of Figure~\ref{fig:barbells}, the cylinder is finally short enough that the eigenfunctions correspond to something like 
$X$, $Y$ and $Z$, so the image is not flattened, and the eigenmap is an embedding.   

Please note that these are images under the Belkin-Noyogi eigenmaps, $V_t^3$, not the images under the smooth
eigenmaps, $\Phi_t^3$, nor the heat embeddings, $\Psi_t^3$ of B\'erard-Besson-Gallot.  Nevertheless there is an intuitive
relationship between these concepts.

It should be noted that, to save on computation time, the students
did not construct an adjacency graph but used the Euclidean distances to define their matrix $W$ (since points that were not close lead to very small entries in $W$ anyway).   Note also that shortest paths cutting across lattices of points are not
approximating the lengths of direct line segments, so it is not clear that these are better than Riemannian lengths when trying to approximate the heat kernel of a Laplacian.

\subsection{Eigenmaps of Lollipops} \label{sect-lollipop}

\begin{figure}[h] %  figure placement: here, top, bottom, or page
   \centering
       \includegraphics[width=2.8in]{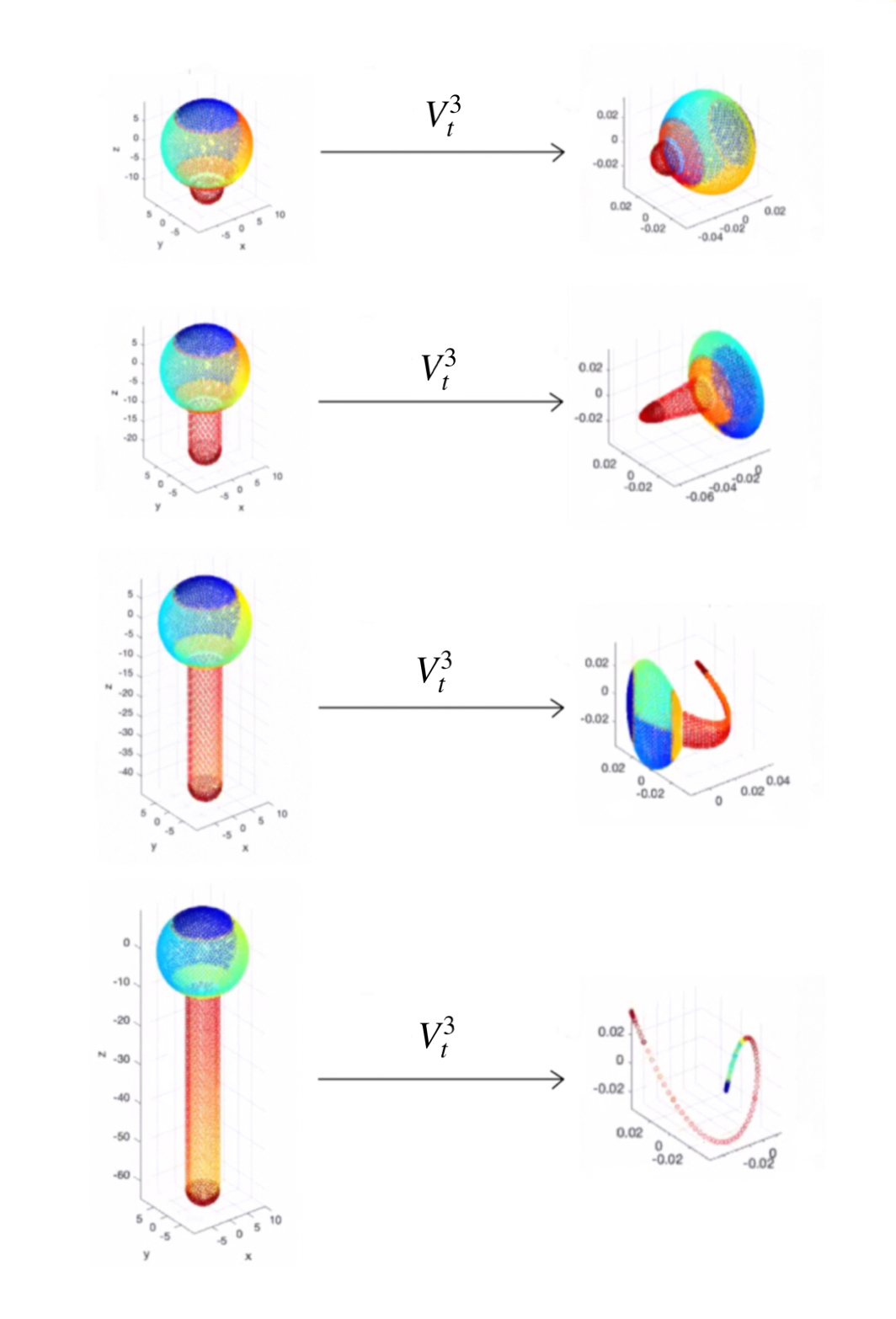} 
   \caption{Here we see a data set of evenly distributed points on various lollipops mapped by an eigenmap,
   $V_t^3$, which gives some possible insight into the image under a heat kernel embedding, 
   $\psi_t^3$ of these lollipops.    MATLAB computation by Maziar Farahzad, Julinda Pillati Mujo, and Esteban Alcantara.
      }
   \label{fig:lollipops}
\end{figure}

Our 2020 student research team (Maziar Farahzad, Julinda Pillati Mujo, and Esteban Alcantara) next applied the the method of eigenmaps to study what the heat kernel embeddings of deformed spheres shaped like lollipops might look like.   In Figure~\ref{fig:lollipops}, they have evenly distributed points on various barbells: {\em a sphere of radius, $10$, with a cap removed and cylinders of radius, $5$, of various lengths, $L=1$, $L=10$, $L=30$, and $L=50$, glued in, and then capped off with a hemisphere of radius $5$}.   

At the top, $L=1$, so the data set lies on a shape which is almost a sphere and the first three nonconstant eigenfunctions are close to the spherical eigenfunctions, $(X,Y,Z)$ that
we discussed in Section~\ref{sect-sphere-evalue}.  The eigenmap seems to be close to $(-Z(x_i), X(x_i), -Y(x_i), )$ giving a rotated image of the original shape just as in Section~\ref{sect-emap-sphere}.   In the second row, $L=10$, we still see an embedding but one direction is becoming more deformed.   The rotationally symmetry is preserved but changing the X axis to a Z axis.

In the third row, $L=30$, and as with the long barbell, there is a first nonzero eigenvalue corresponding to the full height of the manifold which has an eigenfunction that looks something like $\cos(Z)$.  The next two eigenfunctions seem to correspond to $z$ and $x$.   That is the points $x_i$ are mapped to $(\cos(Z(x_i)), Z(x_i), X(x_i))$ and so the image is not an embedding.  It is flattened and then bent around a cosine curve.   

In the bottom row, $L=100$ and now the fact that there is a sphere on one end seems to be undetectible.   The first three eigenfunctions seem to all be functions of $z$ alone causing the image under the eigenmap to be a parametrized curve.   This is similar to the effect that we saw with the thin torus in Section~\ref{sect-thin-evalue}.

Please note again that these are images under the Belkin-Niyogi   eigenmaps, $V_t^3$ not the images under the smooth
eigenmaps, $\Phi_t^3$, nor the heat embeddings, $\Psi_t^3$ of B\'erard-Besson-Gallot.  It should be noted again that, to save on computation time, the students
did not construct an adjacency graph but used the Euclidean distances to define their matrix $W$ (since points that were not close lead to very small entries in $W$ anyway).

\subsection{Coifman and Lafon's Diffusion Maps of Data Sets}\label{sect-Coifman-Lafon}

In 2004-2006 Coifman and Lafon extended this idea by introducing the notion 
of {\em diffusion maps} (DM) of high dimensional data sets 
\cite{Lafon:2004}\cite{Coifman-Lafon:2005} \cite{Coifman-Lafon:2006}.    They allowed for data sets which are not uniformly distributed and discuss various kernel operators, $P$, based on the Riemannian Laplacian, based on the weighted graph Laplacian, and based on other linear operators.  They studied the eigenvalues, $\tilde{\lambda}_j$, of these kernel's operators.   

In particular, Coifman and Lafon studied the heat kernel of the Riemannian Laplacian,
which they denoted as, $P_\tau=e^{-\tau\Delta}$.  To understand this notation, we may use Taylor series,
\be
P_\tau=e^{-\tau\Delta}=\sum_{n=0}^\infty\frac{(-\tau \Delta)^n}{n!} = I-\tau \Delta+\frac{(\tau\Delta)^2}{2!}-\frac{(\tau\Delta)^3}{3!}+\cdots.
\ee
Note that $e^{-\tau\Delta}$ has the same eigenfunctions, $\phi_j$, as the Laplacian, $\Delta \phi_j=-\lambda_j \phi_j$, as seen by:
\be
e^{-\tau\Delta}\phi_j=\sum_{n=0}^\infty\frac{(-\tau \Delta)^n}{n!} \phi_j= \sum_{n=0}^\infty\frac{(-\tau \lambda_j)^n}{n!} \phi_j
=e^{-\lambda_j\tau} \phi_j.
\ee
Coifman and Lafon wrote their eigenvalues of $P_\tau$ as $\tilde{\lambda}_j^\tau$ decreasing from $\tilde{\lambda}_0^\tau=1$ so that
\be
P_{\tau} \phi_j=\tilde{\lambda}^\tau_j \phi_j. 
\ee 
When $P_{\tau} =e^{-\tau\Delta}$ we have
\be\label{tilde-lambda}
P_{\tau} \phi_j=\tilde{\lambda}^\tau_j \phi_j \iff 
\Delta \phi_j=-\lambda_j \phi_j  \textrm{ so } \tilde{\lambda}_j=e^{-\lambda_j}
\ee

Coifman-Lafon's {\bf \em diffusion map} is a weighted eigenmap into Euclidean space:
\be
\tilde{\Psi}^N_\tau(x)=\{\tilde{\lambda}_j^{\tau}\phi_j(x)\}_{j=1}^N \subset {\mathbb E}^N
\ee
Recall the {\em heat kernel embedding} of B\'erard-Besson-Gallot we saw in (\ref{BBG-map}) 
and observe that by taking $\tau=t/2$ and applying (\ref{tilde-lambda}), 
Coifman-Lafon's diffusion map is a truncation
of the heat kernel embedding of B\'erard-Besson-Gallot:
\be\label{BBG-CL}
\tilde{\Psi}^N_{t/2}(x)=\{e^{-\lambda_j t/2 }\phi_j(x)\}_{j=1}^N=\Psi^N_t(x).
\ee
They did not directly cite \cite{BBG:1994} but they do cite other papers and textbooks that cite the original
article.  

\begin{figure}[h] %  figure placement: here, top, bottom, or page
   \centering
   \includegraphics[width=4in]{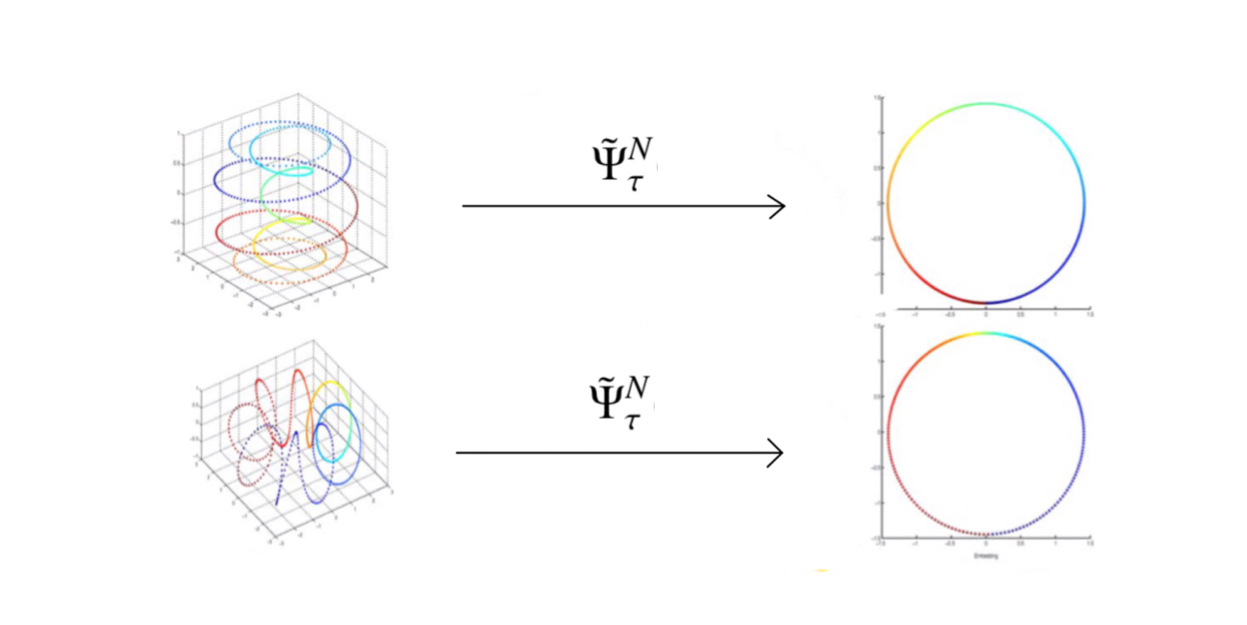} 
   \caption{Coifman-Lafon demonstrate that the images of data sets of evenly distributed points lying on closed curves under a truncated diffusion map, $\tilde{\Psi}^N_\tau$, with N=2 is round circle.  This embedded circle recovers the intrinsic Riemannian geometry of the 1 dimensional submanifolds containing the point sets.  Image Credit: Excerpt from
   Figure 4 in \cite{Coifman-Lafon:2006}.
   }
   \label{fig-Coifman-Lafon}
\end{figure}

Coifman-Lafon demonstrated with graphics as in Figure~\ref{fig-Coifman-Lafon} that this weighted eigenmap captures the intrinsic geometry of a data set lying on a curve.  In particular they demonstrated that collections of points lying on a long twisted and knotted curve in the original space maps onto a perfectly round circle in ${\mathbb E}^2$ capturing the intrinsic geometry of the original curve.   In \cite{Coifman-Lafon:2006}, Coifman and Lafon proved rigorously how $N$ must be large after $\tau$ is chosen to be small to justify the good behavior of their method for any member of their class of kernels.   

\subsection{Diffusion Distance}
In \cite{Coifman-Lafon:2006}, Coifman and Lafon coined the term {\em\bf Diffusion Distance} 
\be
D_\tau(x,y)^2=\sum_{j=1}^\infty \tilde{\lambda}_1^{2\tau}(\phi_j(x)-\phi_j(y))^2
=\sum_{j=1}^\infty e^{-\lambda_j t}(\phi_j(x)-\phi_j(y))^2
\ee
which agrees with pull back of the $\ell^2$ distance of B\'erard-Besson-Gallot when applied to the heat kernel:
\be
D_{t/2}(x,y)= || \Psi_t(x)-\Psi_t(y)||_{\ell^2}=(\Psi_t^*d_{\ell^2})(x,y).
\ee
While we have provided the relationship here between their work and that of B\'erard-Besson-Gallot, remember that
Coifman and Lafon's work applies to a larger class of kernels.   Once again we apologize for not including probability references in this survey that inspired and helped justify the work of Coifman and Lafon.

It is important to note that the Diffusion Distance is an extrinsic distance measured using direct line segments in $\ell^2$.
Even when  studying the heat kernel, 
\be
D_{t/2}(x,y)^2=\sum_{j=1}^\infty e^{-\lambda_j t}(\phi_j(x)-\phi_j(y))^2.
\ee
is an extrinsic distance between points and does not approximate $d_g(x,y)$.  Even with appropriate scaling and taking $t\to 0$, the heat kernel map is not distance preserving, but is only length preserving as we discussed in Section~\ref{sect-length-vs-distance}.   Nevertheless $D_{t/2}(x,y)$ is an interesting distance function
on the original manifold $M$. It induces rectifiable lengths of curves, $C:[0,1]\to M$,
\be
L_t(C)=\sup\left\{ \sum_{i=1}^n D_{t/2}(C(s_{i-1},C(s_i))\,: \,0\le s_0<s_1<\cdots <s_n=1\right\},
\ee
where the sup is taken over all partitions, $0\le s_0<s_1<\cdots <s_n=1$.   These lengths
agree with the lengths of curves defined using the pullback of the Euclidean metric tensor.
By the work of B\'erard-Besson-Gallot in \cite{BBG:1994}, with appropriate rescaling, 
these lengths converge 
\be
\lim_{t\to 0} L_t(C)= L_g(C)
\ee
With that rescaling we have (\ref{Diff-Distance-to-0}).

Memoli and Sapiro computed the
intrinsic distance between points in data sets lying on a manifold using related ideas in \cite{Memoli_Sapiro:2005}

\section{{\bf Spectral Embeddings and Embedding Dimensions}}

In this section we present various results on {\bf \em spectral embeddings} of Riemannian manifolds where
we refer to an embedding as spectral if it is defined using eigenfunctions, eigenvalues, and/or heat kernels.   
There has been significant research activity in this area inspired by the applications of spectral
embeddings to study high dimensional data sets.    It is of particular interest to study the
smallest dimension, $N$, such that $M^n$ can be embedded into ${\mathbb R}^N$
and to find uniform bounds on this dimension $N$.   This is referred to as the { \bf \em embedding dimension}
of $M^n$.   

Note a truncated heat kernel map, $\psi^N_t$, as defined in (\ref{BBG-trunc}) is an embedding iff
the corresponding smooth eigenmap, $\Phi^N$, as defined in (\ref{BN-smooth-eigenmap})
is also an embedding.   This is because the only distinction between 
these maps is the scaling of the various entries using eigenvalues.   Other spectral maps discussed above
have similar uniform bounds on $N$.

\subsection{Deforming Manifolds and their Embeddings}\label{sect-deform}

As we saw in Section~\ref{sect-embed-map}, B\'erard, Besson, and Gallot 
proved that for any $\epsilon>0$ there is a $t_\epsilon>0$ sufficiently small
such that for any $t\in (0,t_\epsilon)$, there is an $N_{t,\epsilon}$ sufficiently large
that for any $N>N_{t,\epsilon}$ we have
\be
| (\psi_t^N)^*g_{{\mathbb E}^N} -g|<\epsilon
\ee
where
 \be
\psi^N_t(x)= \sqrt{2}(4\pi)^{n/4} t^{(n+2)/4} \{e^{-\lambda_j t/2 }\phi_j(x)\}_{j=1}^N \in {\mathbb E}^N.
\ee
In particular for $t$ sufficiently small, and $N$ sufficiently large, they can guarantee that
$(\psi_t^N)^*g_{{\mathbb E}^N}$ is close enough to $g$ to be positive definite.  In fact
for sufficiently large $N$ the truncated eigenmap, $\psi^N_t$, is an embedding.
We saw that for a thin torus, we would need to take $N$ very large to achieve this in Section~\ref{sect-length-vs-distance}.
We also saw that the extrinsic distances were not well controlled as $t\to 0$ in Section~\ref{sect-length-vs-distance}.

If the manifold has computable eigenfunctions, like a sphere or a torus, then we can compute $t_\epsilon$
and $N_{t,\epsilon}$ for that specific manifold.   In fact, we know we only need three eigenfunctions on a sphere
and four on a torus to achieve an embedding.  Naturally it is of interest to have similarly strong controls on manifolds that are close to these special manifolds.

In fact B\'erard, Besson, and Gallot proved their embeddings are stable under $C^0$ perturbations of the metric with a uniform lower bound on Ricci curvature (See Theorem 21 in \cite{BBG:1994}):

\begin{thm}\cite{BBG:1994}
Let $(M,g)$ be a closed $n$-dimensional Riemannian manifold, $\epsilon_0 >0$ and $N_0$ be a positive integer. Let $g'$ be any metric on $M$ 
such that 
\be
(1-\epsilon)g \leq g' \leq (1+\epsilon)g, \epsilon < \epsilon_0.
\ee
 We assume that all metrics under consideration satisfy  $\mathrm{Ric}_{(M,g')} \geq -(n-1)H$ 
for some constant $H\ge 0$. Then there exist constants $\eta_{g,i,H}(\epsilon)$, $1 \leq i \leq N_0$, which go to $0$ with $\epsilon$, such that to any orthonormal basis $\{\phi'_i \}$ of 
eigenfunctions of $\Delta_{g'}$ we can associate an orthonormal basis $\{ \phi_i \}$ of eigenfunctions of $\Delta_g$ satisfying 
\be
\| \phi_i - \phi'_i \|_{L^{\infty}} \leq \eta_{g,i,H}(\epsilon)
\ee
for all $i \leq N_0$. 
\end{thm}

As a consequence, suppose we consider the standard sphere ${\mathbb S}^2$.  We saw in Section~\ref{sect-sphere-map} that its first three nonconstant eigenfunctions map it to a rescaled sphere in ${\mathbb E}^3$.   If
$(M,g)$ is $C^2$ close to ${\mathbb S}^2$ then it's eigenfunctions are close to that of the standard sphere
and its truncated eigenmap will map it close to a sphere.   If
$(M,g)$ is only $C^0$ close to ${\mathbb S}^2$ but has Ricci curvature bounded below, then this works as well.
This works for lollipops with fixed cylinder radius, $r$, and length $L_i\to 0$.
The increasingly thin lollipops with $r_i\to 0$ and $L$ fixed are not $C^0$ converging to a standard sphere
and do not have Ricci curvature bounded from below uniformly, so we cannot apply the B\'erard, Besson, and Gallot 
$C^0$ stability theorem to study them.   

\subsection{Local Embeddings of $M^n$ into ${\mathbb E}^n$}

In 2008, Jones, Maggioni, and Schul studied local properties of embeddings using eigenfunctions. 
They showed that any smooth $n$-dimensional manifold admitting charts in which the metric is $C^{\alpha}$ can be locally embedded in
$\mathbb{E}^n$ by eigenfunctions of the Laplace operator \cite{Jones_Maggioni_Schul:2008, Jones-Maggioni-Schul:2010}. 
These coordinates are bi-Lipschitz on embedded balls of the domain or manifold, with distortion constants that depend only on natural 
geometric properties of the domain or manifold.  

\subsection{Embedding Dimensions of Spectral Embeddings}

  In 2014 Bates considered the class
  \be\label{Bates-class}
\mathcal{M}_{H, i_0}^V = \{ (M,g) : \dim M = n, \mathrm{Ric} \geq -(n-1)H, \mathrm{inj}(M) \geq i_0, \mathrm{Vol}(M) \leq V \}
\ee
in \cite{Bates:2014}. 
Recall this is the class Anderson-Cheeger proved was $C^\alpha$ compact in \cite{Anderson-Cheeger:1992}.

Bates proved  that any closed, connected Riemannian manifold $M$ can be embedded 
  using a smooth truncated eigenmap where the number of required eigenfunctions
  is uniform on the class:

\begin{thm}\cite{Bates:2014}
Given any $H\in {\mathbb R}$, any $i_0>0$ and any $V>0$, there
 is a positive integer $N=N_{H,n,i_0,V}$ and constant $\epsilon_{H,n,i_0,V} >0$, such that, for any 
$M^n \in \mathcal{M}_{H, i_0}^V$
and for all $z\in M$,
\be
\Phi^N: M \to \mathbb{E}^N \textrm{ where}
x \mapsto (\phi_1(x), \cdots, \phi_N(x)) 
\ee
is a smooth embedding.   
\end{thm}

As a consequence there is a uniform upper bound on the embedding dimension of $M$ on this class
using eigenmaps as embeddings.   Note that 
in his proof, Bates 
applied the work of Jones, Maggioni, and Schul that we mentioned above.
  
\subsection{Uniform Bounds on Spectral Embeddings by Portegies}\label{sect-Portegies}

In 2016, Portegies considered various spectral embeddings on the same class $\mathcal{M}_{H, i_0}^V$.   In particular he studied the truncated heat kernel
map of B\'erard-Besson-Gallot, 
\be
\psi^N_t(x)= \sqrt{2}(4\pi)^{n/4} t^{(n+2)/4} \{e^{-\lambda_j t/2 }\phi_j(x)\}_{j=1}^N \in {\mathbb E}^N,
\ee
that we reviewed in Section~\ref{sect-trunc-map}.   Recall that 
 B\'erard-Besson-Gallot had proven that
 the pull back of the metric tensor under their map is close to the original metric tensor
 under the following limits:
\be
\lim_{t\to 0^+}\lim_{N\to \infty} (\psi_t^N)^*g_{{\mathbb E}^N} (V_x,V_x)= g(V_x,V_x).
\ee
Portegies proved this convergence is uniform on $\mathcal{M}_{H, i_0}^V$ in \cite{Portegies:2016}:

\begin{thm}
Let $\epsilon >0$. Then there exists a $t_0 = t_0(n,\epsilon, H, i_0 )$ such that for all $0 < t \leq t_0$ there exists an $N_0 = N_0(n, \epsilon, t, H, i_0, V)$ such that if $N \geq N_0$, for all 
$(M,g) \in \mathcal{M}_{H, i_0}^V$, the map 
\be
\psi^N_t(x) = (2t)^{(n+2)/4} \sqrt{2}(4 \pi)^{n/4} (e^{-\lambda_1 t} \phi_1(x), \cdot, e^{-\lambda_Nt} \phi_N(x))
\ee
is an embedding of $M$ into $\mathbb{E}^N$ such that 
\be
1 - \epsilon < |(d\psi^N_t)_p | < 1 + \epsilon.
\ee
\end{thm}

In particular, the sequence of Riemannian manifolds $(M,  (\psi^N_t)^*g_{{\mathbb E}^N})$
converges in the $C^0$ sense as $N\to \infty$ to the original manifold $(M,g)$.   The speed of this
convergence is the same for any $(M,g) \in \mathcal{M}_{H, i_0}^V$.  Note that
the sequence of Riemannian manifolds $(M,  (\Phi^N)^*g_{{\mathbb E}^N})$ defined by pulling back metric tensors using eigenmaps diverges.  We need the rescaling by the eigenvalues prevent the
lengths from stretching to infinity in the limit.

Portegies proved this theorem by first obtaining a quantitative version of the estimate on the harmonic radius found by Anderson-Cheeger in \cite{Anderson-Cheeger:1992}.   He then studied the images under
the map $\psi^N_t$ restricted to these harmonic charts.   In addition to proving the above theorem he
proved bounds on other truncated spectral embeddings defined using heat kernels \cite{Portegies:2016}.

\subsection{Spectral Embeddings depending on the Connection Laplacian}

There are other spectral embeddings based on other types of Laplacian operators, such as the connection Laplacian. The connection Laplacian is a Laplace operator acting on the various tensor bundles of a manifold. Inspired by the cryo-electron microscopy image problem for reconstructing macromolecular structures, Singer-Wu \cite{Singer_Wu:2012} introduced the {\bf \em vector diffusion maps} (VDM). The VDM method helped them
with aligning images and finding similar images (nearest neighbors) more effectively than the diffusion maps
studied by Coifman-Lafon. 

The VDM is a spectral embedding  \cite{Lin_Wu:2018} that depends on the connection 
Laplacian associated with a possibly nontrivial bundle structure. 
Its general goal is to integrate local information and the
relationship between these pieces of local information in order to obtain the global information
of the dataset; for example, the ptychographic imaging problem \cite{Marchesini_Tu_Wu:2014}, the synchronization problem \cite{Bandeira_Singer_Spielman:2013}, the vector nonlocal
mean/median \cite{lin2018manifold}, the orientability problem \cite{Singer_Wu:2011}, etc. Numerically, the VDM depends on the
spectral study of the graph connection Laplacian  \cite{Chung_Zhao_Kempton:2013, ElKaroui_Wu:2015a, ElKaroui_Wu:2015b, Singer_Wu:2012, Singer_Wu:2016} , which is a
direct generalization of the graph Laplacian discussed in the spectral graph theory 
text by Chung \cite{Chung:text}.\\

In \cite{Lin_Wu:2018}, Lin and Wu investigated estimates of eigenvector fields and truncations of heat kernels and its derivatives. In the end, they showed that an embedding of a manifold can be constructed with finite eigenvector fields of the connection Laplacian for a class of closed smooth manifolds
\be
\mathcal{M}_{n, H, i_0, V} = \{ (M,g) : |\mathrm{Ric_g}| \leq H, \mathrm{inj}(M) \geq i_0, \mathrm{Vol}(M) \leq V \}.
\ee
Recall this is the class that Anderson proved was $C^{1,\alpha}$ precompact in \cite{Anderson:1990}.

In \cite{Lin:2022}, Lin proved that, in addition, we can ensure the pull back of the metric tensor is close to 
the original metric tensor under a limiting process and this can be done uniformly on the class
$\mathcal{M}^H_{n, H, i_0, V}$ obtaining the following theorem:

\begin{thm} \cite{Lin:2022}
Given $\epsilon >0$, there exists $t_0 = t_0(n,H,i_0, \epsilon)$ so that for all $t< t_0$ there exists $N_0 = N_0(n, H, i_0, \epsilon, t, V)$ such that for all 
$(M,g) \in \mathcal{M}_{n, H, i_0, V}$, there exist
points $p_1, \cdots, p_{N_0}$ on $M$ such that the map 
$\mathcal{H} : (M^n, g) \to \mathbb{E}^{N_0}$
defined by
\be
\mathcal{H}(p)= \frac{A(2t)^{(3n+2)/4}}{V_e} \left( \| K_{TM}(p,t; p_1)\|_{HS}^2, \cdots,  \| K_{TM}(p,t; p_{N_0})\|_{HS}^2\right)
\ee
is an embedding, where $A = A(n, H, i_0, \epsilon, t, V)$ is a constant,
 $K_{TM}$ denotes the heat kernel of the connection Laplacian on the tangent bundle $TM$, $\| \cdot \|_{HS}$ denotes the Hilbert-Schmidt norm, 
and 
$$
V_e = \left( \int_{\mathbb{R}^n} (\partial_{x_1} \| K^{\mathbb{E}^n}_{TM}(0, \frac 12; y)\|^2_{HS})^2 dy \right)^{\tfrac{1}{2}\,}.
$$
Here, $K^{\mathbb{E}^n}_{TM}$ is the standard Euclidean heat kernel on $\mathbb{E}^n$. In addition, 
$$
1 - \epsilon < |(d\mathcal{H})_p | < 1 + \epsilon \mbox{ for all $p\in M$.}
$$ 
The same statements hold with every heat kernel $K_{TM}$ replaced by the truncated version $K_{TM}^N$.   See \cite{Lin:2022}.
\end{thm}

The proof of this theorem also relies on a quantitative version of the estimate on the harmonic.   
Since Lin must control the tangent bundle $TM$
and the heat kernel of the connection Laplacian action on sections of the tangent bundle
 in \cite{Lin:2022}, it
is not surprising that Lin needs to study a class with $C^{1,\alpha}$ compactness.  

\subsection{$L^p$ Bounds on Pull Back Metrics}\label{sect-Lp}

In 2021, Ambrosio, Honda, Portegies and Tewodrose considered the class
\be
\mathcal{M}_{H,V}^{D}=\{ (M,g) : |\mathrm{Ric_g}| \leq H, \,\mathrm{Vol}(M) \geq V ,
\,\mathrm{diam}(M)\le D\}.
\ee
in \cite{AHPT:2021}. 
Note that this class has no uniform lower bound on injectivity radius.   This is the class
of noncollapsing manifolds with lower Ricci curvature bounds studied by Cheeger-Colding
that have subsequences which converge in the metric measure sense to limits
and have controls on their eigenfunctions \cite{ChCo-Part3}.  Yu Ding controlled the
heat kernels on such spaces in \cite{Ding:02}.  The limits of these spaces are in fact
$\mathrm{RCD}^*(H,N)$ spaces which were
first introduced by Ambrosio-Gigli-Savare \cite{Ambrosio-Gigli-Savare:2014}.

In fact Ambrosio, Honda, Portegies and Tewodrose considered
spectral maps defined on these $\mathrm{RCD}^*(H,N)$ spaces in \cite{AHPT:2021} extending many
known results to this class.   Theorem 6.8 and 6.9 in their paper \cite{AHPT:2021}
are also new results for smooth Riemannian manifolds in $\mathcal{M}_{H,V}^{D}$.
Theorem 6.9 concerns an eigenmap, so we include it
below as reformulated for us by Portegies for manifolds:  

\begin{thm}[Theorem 6.9 in \cite{AHPT:2021}]
For all $D > 0, V > 0, \epsilon > 0$ and $1 \leq p < \infty$ there exists \[t_0 := t_0(H, N, V, D, \epsilon, p) > 0\] such that for all $0 < t \leq t_0$ and every closed $n$-dimensional Riemannian manifold 
$M\in \mathcal{M}_{H,V}^{D}$, there exists an $N_0 = N_0(H, N, V, D, \epsilon, p, t)$ such that for all $N \geq N_0$, we have
\[
\| \, | \omega_n t^{(n+2)/2} g_t^N - c_n g|_{HS} \|_{L^p} \leq \epsilon.
\]
where
\[
g_t^N := \Psi_t^* g_{{\mathbb E}^N}=\sum_{i=1}^N e^{-2\lambda_i t}d \phi_i \otimes d\phi_i.
\]
\end{thm}

Observe that this theorem is stating that the pull back metric tensor is converging to the original tensor
but only with uniform $L^p$ bounds on the convergence.   In Conjectures~\ref{conj-GH-t-N} and~\ref{conj-SWIF-t-N} stated at the end of this paper, we ask whether we can achieve Gromov-Hausdorff or Intrinsic Flat Convergence with the right
choice of class ${\mathcal{M}}$ of Riemannian manifolds.   Note that $L^p$ convergence of the
metric tensor is not enough to achieve uniform, Gromov-Hausdorff, or intrinsic flat convergence
unless there are also uniform bounds on
the metric tensors from above and below as shown by Allen-Sormani in \cite{Allen-Sormani:2019} and \cite{Allen-Sormani:2020} or controls like those required 
by Allen-Perales-Sormani to achieve VADB convergence in \cite{APS-JDG}.

\section{\bf Open Questions}

Throughout this paper we have been exploring the properties of heat kernels, eigenfunctions, and
spectral maps for a variety of manifolds and classes of manifolds.   In this section, 
we vaguely state some conjectures related to these ideas.   We hope to bring together mathematicians
from around the world to study these problems together.

\subsection{Conjectures on Truncated Eigenmaps of Small Dimension}

In the work of Belkin-Niyogi reviewed in Section~\ref{sect-Belkin-Niyogi} and the
work of Coifman-Lafon reviewed in Section~\ref{sect-Coifman-Lafon}, truncated spectral maps were quite useful even when $N=2$ or $N=3$ and without taking $t\to 0$.   This is due in part to the underlying symmetries of the data sets that they are studying.   In Section~\ref{sect-sphere-map} we saw that the  symmetries of a sphere guarantee that it spectrally 
embeds into Euclidean space with only $N=3$ and does so symmetrically.  Data sets lying on surfaces
which are close to a sphere in the appropriate sense also have nearby eigenfunctions and thus also embed with $N=3$.  Other surfaces like the barbell and the lollipop studied by our students did not always spectrally embed with $N=3$
as seen in Sections~\ref{sect-barbell} and~\ref{sect-lollipop}.   Even the torus ${\mathbb S}^1\times {\mathbb S}^1$ needs $N=4$ and a thin torus needs a much
larger $N$ as discussed in Sections~\ref{sect-tori-map} and~\ref{sect-tori-smap}.

This leads naturally to the following conjecture:

\begin{conj}
Can we describe a class of manifolds $\mathcal{M}$  for which 
\be 
\Psi^N_t: M \to {\mathbb E}^N
\ee
is an embedding when $N=3$ or $N=4$?   
\end{conj}

Since eigenvalues and eigenfunctions vary smoothly with $C^2$ convergence of manifolds, we know
that this conjecture holds with $N=3$ for the class of manifolds which are sufficiently close to a standard sphere in the $C^2$ sense, and it holds with $N=4$ for the class of manifolds which are sufficiently close to a standard torus in the $C^2$ sense.  We saw in Section~\ref{sect-deform} that B\`erard-Besson-Gallot proved that if $M$ have Ricci curvature bounded uniformly from below and are $C^0$ close, then their eigenfunctions and eigenvalues are close, thus this conjecture holds with $N=3$ for classes of manifolds $C^0$ close to the
standard sphere with Ricci curvature bounded uniformly from below, and similarly for $N=4$ near the standard torus.   However there is no precise statement as to how close these manifolds must be.   It would also be interesting to define a class which does not mention the sphere or the torus.

\subsection{Conjectures on Selective Eigenmaps} \label{sect-conj-selective}

In Section~\ref{sect-tori-smap} we observed how shapes like a thin torus would not be embedded by
a truncated spectral map until $N$ is taken very large and then there has not been a very useful
reduction in dimension.   We suggested the idea of selecting eigenfunctions to define an embedding
rather than just truncating a spectral embedding.   
In particular,  we suggested only selecting eigenfunctions with eigenvalues near a given eigenvalue 
$\lambda$.   This leads to the following natural conjecture:

\begin{conj} \label{conj-selective}
If we are given a Riemannian manifold, $M^n$, is a class, $\mathcal{M}$, with $\injrad\ge i_0>0$
and we define a {\em selected spectral map} using only
eigenfunctions whose eigenvalues are close to $\lambda_0$ defined as a function of $i_0$,
can we prove this selected spectral map is an spectral embedding?
\end{conj}

Note that in the conjecture above we have been deliberately vague about the choice of $\lambda_0$
and how it depends on $i_0$.   We have also been vague about the class, $\mathcal{M}$.  It
might be natural to require the class have bounds on Ricci curvature, volume, or diameter.

\begin{conj}
Can we describe a class of manifolds $\mathcal{M}$  for which 
we can find {\em selected spectral map} that 
is an embedding when $N=3$ or $N=4$?   
\end{conj}

\subsection{Conjectures on Convergence of Pullbacks of Embeddings}

Recall in Section~\ref{sect-trunc-map} we reviewed the truncated heat kernel
map of B\'erard-Besson-Gallot, 
\be
\psi^N_t(x)= \sqrt{2}(4\pi)^{n/4} t^{(n+2)/4} \{e^{-\lambda_j t/2 }\phi_j(x)\}_{j=1}^N \in {\mathbb E}^N,
\ee
has the property that:
$
\forall \epsilon > 0, \,\exists t_\epsilon>0,\, \textrm{ s.t. } \forall t\in (0,t_\epsilon),
\,\exists N_{t,\epsilon} \textrm{ s.t. } \forall N \ge N_t \quad $
\be
| g_t^N - g| <\epsilon
\ee
where $g_t^N=(\psi_t^N)^*g_{{\mathbb E}^N}$.   In Section~\ref{sect-Portegies}, we
saw that Portegies proved the choice of $t_\epsilon$ and of $N_{t,\epsilon}$ can be
made uniform on the class $\mathcal{M}_{H, i_0}^V$.

Recalling the definition of the Lipschitz distance in Section~\ref{Lip-Conv}, we can rewrite this as
$
\forall \epsilon > 0, \,\exists t_\epsilon>0,\, \textrm{ s.t. } \forall t\in (0,t_\epsilon),
\,\exists N_{t,\epsilon} \textrm{ s.t. } \forall N \ge N_t \,\,$
\be
d_{Lip}\left( (M,d_{g_t^N}) - (M,d_g) \right) <\epsilon
\ee
where the choice of $t_\epsilon$ and of $N_{t,\epsilon}$ can be
made uniform on the class $\mathcal{M}_{H, i_0}^V$.

Recall the definition of Gromov-Hausdorff distance in Section~\ref{sect-GH}.

\begin{conj}\label{conj-GH-t-N}
Can we find a class of manifolds, ${\mathcal{M}}$ such that
$
\forall \epsilon > 0$, $\exists t_\epsilon>0,\, \textrm{ s.t. } \forall t\in (0,t_\epsilon),
\,\exists N_{t,\epsilon} \textrm{ s.t. } \forall N \ge N_t \,\,$
\be
d_{GH}\left( (M,d_{g_t^N}) - (M,d_g)\right) <\epsilon
\ee
where the choice of $t_\epsilon$ and of $N_{t,\epsilon}$ can be
made uniform on the class $\mathcal{M}$?
\end{conj}

Recall the definition of the intrinsic flat distance in Section~\ref{sect-SWIF}.

\begin{conj}\label{conj-SWIF-t-N}
Can we find a class of manifolds, ${\mathcal{M}}$ such that
$
\forall \epsilon > 0, \,\exists t_\epsilon>0,\, \textrm{ s.t. } \forall t\in (0,t_\epsilon),
\,\exists N_{t,\epsilon} \textrm{ s.t. } \forall N \ge N_t \,\, $
\be
d_{GH}\left( (M,d_{g_t^N}) - (M,d_g)\right) <\epsilon
\ee
where the choice of $t_\epsilon$ and of $N_{t,\epsilon}$ can be
made uniform on the class $\mathcal{M}$?
\end{conj}

Since Lipschitz convergence implies Gromov-Hausdorff and Intrinsic Flat convergence, both of
the above conjectures are immediately true for the class $\mathcal{M}_{H, i_0}^V$.  However they
should be provable without the injectivity radius bound.   See the discussion in Section~\ref{sect-Lp}.
%I think we can prove these immediately.

\subsection{Conjectures on Convergence of Manifolds and their Embeddings}

We saw in Section~\ref{sect-deform} that B\`erard-Besson-Gallot proved that if $M_j$ have Ricci curvature bounded uniformly from below and 
$M_j \to M_\infty$ in the $C^0$ sense, then $\psi_t^N(M_j) \to \psi_t^3(M_\infty)$ in 
the $C^0$ sense because the eigenfunctions converge.

Recall the definition of Hausdorff and Gromov-Hausdorff convergence in Sections~\ref{sect-GH}
and~\ref{sect-mGH}:

\begin{conj}
Suppose $M_j \subset \mathcal{M}$ for some class $\mathcal{M}$, and $M_j \to M_\infty$ in the measured Gromov-Hausdorff sense, does $\psi_t^N(M_j) \to \psi_t^N(M_\infty)$ in the Hausdorff sense as subsets of Euclidean Space, ${\mathbb E}^N$?  
\end{conj}

Recall the definition of the intrinsic flat convergence in Section~\ref{sect-SWIF}.

\begin{conj}
Suppose $M_j \subset \mathcal{M}$ for some class $\mathcal{M}$, and $M_j \to M_\infty$ in the SWIF sense, does $\psi_t^N(M_j) \to \psi_t^N(M_\infty)$ in the Flat sense
as submanifolds of Euclidean Space, ${\mathbb E}^N$?  
\end{conj}

These conjectures are quite deliberately vague as the specific class of manifolds has not been prescribed.
Both conjectures are probably easily proven to be true if we choose a class with a lower bound on injectivity radius and on Ricci curvature using the theorems surveyed within this article.   It would be far more challenging
to remove the injectivity radius bound.   

\vspace{.2cm}
{\em Please contact us if you would like to work on any of the conjectures in this survey.   We are happy to set up teams of collaborators to work together and might invite you to participate in workshops on the topic as well.}

\vspace{.2cm}

\section*{Acknowledgements}
 We would particularly like to thank our 2020 student research team members: Esteban Alcantara, Julinda Pillati Mujo, and Maziar Farahzad.    Professor Sormani would like to thank the organizers of the 2022 TWAS TYAN Virtual Workshop on Differential Geometry for the honor of being invited to moderate a session.  It was truly inspiring to meet Professor Romain Murenzi,: to hear him speak of his experience as a mathematics student and teacher from Burundi, and how important it can be to work in a young active subject.  He spoke of how Professor Alex Grossman adopted him into his research group and gave him the opportunity to speak internationally, and encouraged us all to do the same for young scientists from developing countries.   Here at the City University of New York, Professor Chen-Yun Lin and Professor Sormani often teach first and second generation students from developing countries and are delighted to build connections with faculty and students internationally.   

%\bibliography{database-doses}
\bibliographystyle{plain}
\bibliography{CYLIN-SORMANI-22.bib}

\strictpagecheck
\checkoddpage
\ifoddpage
\newpage{\ }\thispagestyle{empty}
\fi

\end{document}